\PassOptionsToPackage{pagebackref}{hyperref}
\documentclass[11pt]{article}
\usepackage[margin=1in]{geometry}
\usepackage{amsmath,amssymb,amsthm,mathtools,bm}
\usepackage{enumitem}
\usepackage{booktabs,array,longtable}
\usepackage{microtype}
\usepackage{hyperref}
\usepackage{url}
\usepackage{xcolor}
\usepackage{setspace}
\usepackage{bbm}
\usepackage{mathrsfs}
\usepackage{graphicx}
\usepackage{fancyhdr}
\usepackage{titlesec}
\usepackage{natbib}
\usepackage{algorithm}
\usepackage[pagebackref]{hyperref}

\renewcommand*{\backref}[1]{}
\renewcommand*{\backrefalt}[4]{%
  \ifcase #1\relax
  \or (Cited on page~#2.)%
  \else (Cited on pages~#2.)%
  \fi
}
\usepackage{algpseudocode}
\usepackage{titletoc}
\setcitestyle{authoryear,round}
\usepackage{placeins}
\hypersetup{colorlinks=true,linkcolor=blue,citecolor=blue,urlcolor=blue}
\newtheorem{theorem}{Theorem}[section]
\newtheorem{proposition}[theorem]{Proposition}
\newtheorem{lemma}[theorem]{Lemma}
\newtheorem{corollary}[theorem]{Corollary}
\newtheorem{assumption}[theorem]{Assumption}
\newtheorem{remark}[theorem]{Remark}

\newcommand{\PP}{\mathbb P}
\newcommand{\EE}{\mathbb E}
\newcommand{\RR}{\mathbb R}

\newcommand{\cO}{\mathcal O}
\newcommand{\cT}{\mathcal T}
\newcommand{\cH}{\mathcal H}
\newcommand{\cS}{\mathcal S}
\newcommand{\cU}{\mathcal U}

\newcommand{\eps}{\varepsilon}
\newcommand{\ind}{\mathbbm 1}
\newcommand{\tr}{\operatorname{tr}}
\newcommand{\rank}{\operatorname{rank}}

\newcommand{\argmax}{\operatorname*{arg\,max}}
\newcommand{\argmin}{\operatorname*{arg\,min}}

\newcommand{\weakto}{\rightsquigarrow}

\newcommand{\pto}{\xrightarrow{P}}

\newcommand{\norm}[1]{\left\lVert #1\right\rVert}
\newcommand{\abs}[1]{\left|#1\right|}

\newcolumntype{P}[1]{>{\raggedright\arraybackslash}p{#1}}

\title{\bfseries Weak Limits of Parameter Estimation \\
in Strongly Identifiable Mixture Models}
\author{}
\date{}

\begin{document}
%\maketitle

\begin{center}

{\bfseries\Large
Weak Limits of Parameter Estimation \\
in Strongly Identifiable Mixture Models
}

\vspace*{.2in}

{\large
\begin{tabular}{cccc}
Dung Le &
Trang Pham &
Alessandro Rinaldo &
Nhat Ho
\end{tabular}
}

\vspace*{.2in}

\begin{tabular}{c}
The University of Texas at Austin
\end{tabular}

\vspace*{.2in}
\today

\begin{abstract}
% Root-$n$ density asymptotics for finite mixtures do not by themselves determine the distribution of mixing-measure loss. We study global maximum likelihood estimation over mixtures with at most $k$ components when the true mixing measure has $p\le k$ distinct atoms. Under second-order strong identifiability, we establish two geometric results. First, we characterize the exact attainable density tangent set, including the allocation constraint $q+\sum_i\rank(C_i)\le k-p$ between nonlocal vanishing components and local splitting directions. Second, we prove a sequential first-order expansion of squared Wasserstein loss: every sequence attaining a given relative-density direction has the same leading transportation cost. Combining the two geometric results with a Gaussian likelihood argmax theorem yields joint weak convergence of the relative-density maximum likelihood estimator and $\sqrt n W_2^2(\widehat G_n,G_*)$, giving an estimator-level, generally non-Gaussian limit for $n^{1/4}W_2(\widehat G_n,G_*)$. As an inferential consequence, an $m$-out-of-$n$ parametric bootstrap consistently approximates the loss distribution when $m\to\infty$ and $m/n\to0$. A one-to-two-component Gaussian calculation illustrates the two competing tangent mechanisms. Numerical studies distinguish weak-limit scaling from finite-sample bootstrap and optimization errors, and a flow-cytometry application illustrates label-invariant uncertainty summaries under a working mixture model.
We study the asymptotic behavior of global maximum likelihood estimation in finite mixtures with an unknown number of components bounded by \(k\). Although convergence rates are well understood under strong identifiability, they do not determine the limiting distribution of mixing-measure loss. Under second-order strong identifiability, we establish two geometric results that connect density asymptotics to uncertainty about the mixing measure. First, we characterize the exact attainable density tangent set, accounting for the component budget shared by nonlocal vanishing components and local splitting directions. Second, we derive a sequential first-order expansion of squared Wasserstein loss, showing that its leading transportation cost is uniquely determined by the relative-density direction. Together with a Gaussian likelihood argmax theorem, these results yield joint weak convergence of the relative-density maximum likelihood estimator and \(\sqrt n W_2^2(\widehat G_n,G_*)\). The density limit maximizes a Gaussian likelihood process over the attainable tangent set, and the induced limit of \(n^{1/4}W_2(\widehat G_n,G_*)\) is generally non-Gaussian. We establish consistency of an \(m\)-out-of-\(n\) parametric bootstrap when \(m\to\infty\) and \(m/n\to0\). Numerical studies assess the approximation and distinguish bootstrap and optimization errors; a flow-cytometry application illustrates uncertainty summaries invariant to component labeling under a working mixture model.
\end{abstract}
\end{center}


\begingroup
%\color{red}
\section{Introduction}\label{sec:intro}

Finite mixture models combine flexible distributional modeling with
an interpretable description of population heterogeneity
\citep{Lindsay1995,McLachlanPeel2000}. With suitable kernels and
sufficiently many components, they can approximate arbitrary densities
in $L^1$ \citep{NguyenEtAl2023Approximation}. Their generative interpretation assigns each observation
to a latent component, providing a natural basis for clustering. Concretely, a finite mixture has density
\[
    p_G(x)
    := \int_{\Theta} f(x\mid\theta)\,dG(\theta)
    = \sum_{j=1}^{p} \pi_j f(x\mid\theta_j),
    \qquad x\in\mathbb R^d,
\]
where $\{f(\cdot\mid\theta):\theta\in\Theta\}$ is a family of
probability density kernels and $G=\sum_{j=1}^{p}\pi_j\delta_{\theta_j}$
is the \emph{mixing measure}, with distinct atoms
$\theta_1,\ldots,\theta_p\in\Theta$, positive weights
$\pi_1,\ldots,\pi_p$, and $\sum_{j=1}^{p}\pi_j=1$.
Equivalently, an observation is generated by first selecting a latent
component $Z$ with $\mathbb P(Z=j)=\pi_j$ and then drawing
$X\mid Z=j$ from $f(\cdot\mid\theta_j)$. The mixing measure thus
encodes both the component parameters and their relative prevalence
in the population.

Given independent observations $X_1,\ldots,X_n$ from $p_{G_*}$,
two related but distinct estimation problems arise. The first is
\emph{density estimation}, which seeks to recover the observable
distribution $p_{G_*}$. The second is \emph{mixing-measure estimation},
which seeks to recover the latent distribution $G_*$, including
the component parameters and their mixing proportions. The former
describes the distribution of the observed data, whereas the latter
describes the underlying population heterogeneity.

\subsection{Related literature}

Maximum likelihood estimation is a standard approach to estimating
both mixture densities and their mixing measures
\citep{KieferWolfowitz1956,Lindsay1995,McLachlanPeel2000}.
For a class $\mathscr G$ of mixing measures, the maximum likelihood
estimator is defined by
\[
    \widehat G_n\in
    \operatorname*{arg\,max}_{G\in\mathscr G}
    \sum_{\ell=1}^{n}\log p_G(X_\ell).
\]

When the model is restricted to a single component, estimation
reduces to ordinary parametric maximum likelihood. Under standard
regularity conditions, the estimator $\widehat\theta_n$ satisfies
\[
    \sqrt n(\widehat\theta_n-\theta_*)
    \rightsquigarrow
    N\bigl(0,I(\theta_*)^{-1}\bigr),
\]
where $I(\theta_*)$ is the Fisher information matrix
\citep{vdV1998}. For example, in a Gaussian location model with
known positive definite covariance $\Sigma$, the maximum
likelihood estimator is the sample mean, and
$\sqrt n(\widehat\mu_n-\mu_*)\rightsquigarrow N_d(0,\Sigma)$.
Asymptotic normality also holds for exactly specified
finite mixtures, after aligning component labels, under
suitable regularity conditions and nonsingular Fisher
information \citep[Theorem~3.1]{RednerWalker1984}.

A prominent framework for studying mixing-measure estimation,
developed by \citet{Nguyen2013} and extended by
\citet{HoNguyen2016,HoNguyen2016Weak}, connects density estimation
to the recovery of latent component parameters through inverse
inequalities. The density-estimation step uses empirical-process
arguments and bracketing entropy to control the complexity of the
mixture density class
\citep{vdVW1996,vanDeGeer2000}. These arguments yield bounds for Hellinger distance of the form
$d_{H}(p_{\widehat G_n},p_{G_*})=O_{P}(\varepsilon_n)$, 
where $\varepsilon_n$ is determined by
the relevant entropy bounds. To measure discrepancy between mixing measures, this framework uses
Wasserstein distances \citep{Villani2009}. For $r\geq1$, define
\[
    W_r(G,G')
    :=
    \left\{
        \inf_{\gamma\in\Gamma(G,G')}
        \int_{\Theta\times\Theta}
        \|\theta-\theta'\|^r\,d\gamma(\theta,\theta')
    \right\}^{1/r},
\]
where $\Gamma(G,G')$ denotes the set of probability measures with
marginals $G$ and $G'$. Inverse inequalities then relate density
discrepancies to transportation costs. In particular, under suitable
smoothness and strong identifiability conditions, one obtains locally
around a fixed $G_*$
\[
    d_{H}(p_G,p_{G_*})\geq c\,W_1(G,G_*)
    \qquad\text{in the exactly specified setting},
\]
and
\[
    d_{H}(p_G,p_{G_*})\geq c\,W_2^2(G,G_*)
    \qquad\text{in the over-specified setting},
\]
for some constant $c>0$
\citep{Nguyen2013,HoNguyen2016}.
The first inequality applies when the fitted and true numbers of
components agree, under first-order strong identifiability.
The second applies when additional components are permitted,
under second-order strong identifiability.

These inequalities translate a density rate $\varepsilon_n$ into
mixing-measure rates $\varepsilon_n$ under $W_1$ and
$\varepsilon_n^{1/2}$ under $W_2$, respectively. Thus, a root-$n$
density rate yields the familiar $n^{-1/2}$ and $n^{-1/4}$ rates;
logarithmic factors in the density bound propagate accordingly.
The change in transportation order reflects the possibility of
first-order cancellation: several fitted components near a true atom
may have weighted displacements that sum to zero, leaving a
second-order density discrepancy.

Subsequent work refines this framework in several directions.
\citet{ManoleHo2022} develop stronger loss functions that distinguish
the convergence rates of individual fitted components, revealing
heterogeneity that a single Wasserstein loss can obscure.
\citet{HeinrichKahn2018} use coarse-graining trees to organize
component configurations across multiple scales and establish uniform minimax rates. More recently, \citet{NguyenLeRinaldoHo2026} quantify
how separation between component parameters affects inverse
inequalities and estimation rates in finite Gaussian location
mixtures. Together, these developments clarify how component
configurations govern estimation accuracy.

Another prominent line of work studies mixing-measure estimation
through moment-based methods
\citep{WuYang2020,DossEtAl2023}.
By relating moment discrepancies to Wasserstein distance, these
methods yield minimax-optimal estimation rates for Gaussian location
mixtures in both univariate and high-dimensional settings without
requiring component separation.

\subsection{Contributions }
Existing work establishes asymptotic normality of likelihood
estimators for exactly specified mixtures under suitable
regularity conditions \citep{RednerWalker1984}, together with
extensive convergence-rate theory for overspecified mixtures
\citep{Chen1995,Nguyen2013,HoNguyen2016,ManoleHo2022}.
These results leave a distinct question: what is the limiting
distribution of the mixing-measure estimation error when redundant
components are permitted? Our goal is to answer this question for
the global maximum likelihood estimator
\[
    \widehat G_n\in
    \operatorname*{arg\,max}_{G\in\mathcal O_k(\Theta)}
    \ell_n(G),
    \qquad
    \ell_n(G):=\sum_{r=1}^{n}\log p_G(X_r),
\]
where $\mathcal O_k(\Theta)$ denotes the mixing measures with at most
$k$ atoms and the true mixing measure $G_*$ has $p\leq k$ distinct
atoms. In the second-order strongly identifiable regime, existing inverse
inequalities translate root-$n$ density convergence into the
fourth-root mixing-measure rate
$W_2(\widehat G_n,G_*)=O_{\mathbb P}(n^{-1/4})$
\citep{Nguyen2013,HoNguyen2016}, but do not determine its limiting
distribution. Under second-order strong identifiability and the stated
regularity conditions, we establish a weak limit for
$n^{1/4}W_2(\widehat G_n,G_*)$. Our analysis connects the local
likelihood geometry to the leading transportation cost through
two geometric results.

Our starting point is the scaled relative-density error
\[
    \widehat h_n
    :=\sqrt n\,\frac{p_{\widehat G_n}-p_*}{p_*},
    \qquad p_*:=p_{G_*}.
\]
To describe its possible limits, we characterize the \textit{attainable tangent set} 
\[
    \mathcal T_k
    :=
    \left\{
        h\in L^2(P_*):
        \begin{array}{l}
        \exists t_m\downarrow0
        \text{ and }G_m\in\mathcal O_k(\Theta):
        \dfrac{p_{G_m}-p_*}{t_m p_*}
        \longrightarrow h
        \text{ in }L^2(P_*)
        \end{array}
    \right\}.
\]
Writing $G_*=\sum_{i=1}^{p}\pi_i^*\delta_{\theta_i^*}$,
we show that every $h\in\mathcal T_k$ admits a unique representation
\begin{equation}
\label{eqn:dung_unique_presentation_of_h}
    h(x)
    =
    \frac{1}{p_*(x)}
    \left[
        \sum_{i=1}^{p}
        \left\{
            \alpha_i f_{\theta_i^*}(x)
            +\beta_i^\top\nabla_\theta f_{\theta_i^*}(x)
            +\frac12
             \operatorname{tr}
             \bigl(C_i\nabla_\theta^2 f_{\theta_i^*}(x)\bigr)
        \right\}
        +\int_\Theta f_\theta(x)\,d\nu(\theta)
    \right],
\end{equation}
where $f_\theta:=f(\cdot\mid\theta)$.
Here, $\alpha_i$ and $\beta_i$ describe changes in mass and location,
the positive-semidefinite matrices $C_i$ describe local splitting,
and the positive atomic measure $\nu$, supported away from the true
atoms, describes nonlocal components with vanishing weights.
These quantities satisfy
\[
    \sum_{i=1}^{p}\alpha_i+\nu(\Theta)=0,
    \qquad
    q+\sum_{i=1}^{p}\operatorname{rank}(C_i)\leq k-p,
    \qquad q:=|\operatorname{supp}(\nu)|.
\]
The latter constraint is necessary and sufficient: it specifies how
the additional fitted components can be allocated between local
splitting and nonlocal perturbations.

This characterization leads to the distributional limit of the
relative-density estimator. For $h_{n,G}:=\sqrt n\,\frac{p_G-p_*}{p_*}$,
a local expansion of the log-likelihood gives
\[
    \ell_n(G)-\ell_n(G_*)
    =
    \frac{1}{\sqrt n}\sum_{r=1}^{n}h_{n,G}(X_r)
    -\frac{1}{2n}\sum_{r=1}^{n}h_{n,G}(X_r)^2
    +R_{n,G}.
\]
Since $P_*h_{n,G}=0$, this becomes
\[
    \ell_n(G)-\ell_n(G_*)
    =
    \mathbb G_n h_{n,G}
    -\frac12\|h_{n,G}\|_*^2
    +\widetilde R_{n,G},
\]
where $\mathbb G_n:=\sqrt n(n^{-1}\sum_{r=1}^n\delta_{X_r}-P_*)$ and
$\|h\|_*^2:=\int h^2\,dP_*$.
We establish uniform control of the remainder on the relevant
localized classes and a functional central limit theorem for the
empirical process. The limiting likelihood criterion is therefore
\[
    Q(h):=\mathbb{W}(h)-\frac12\|h\|_*^2,
    \qquad h\in\mathcal T_k,
\]
where $\mathbb{W}$ is a centered Gaussian process with covariance $\operatorname{Cov}\bigl(\mathbb{W}(h),\mathbb{W}(g)\bigr)
    =\int hg\,dP_*.$ 
Together with convergence of the attainable sets, tightness, and
almost sure uniqueness of the limiting maximizer, a likelihood
argmax theorem yields
\begin{equation}
\label{eqn:likelihood_argmax}
    \widehat h_n\rightsquigarrow H^\dagger,
    \qquad
    H^\dagger
    :=
    \operatorname*{arg\,max}_{h\in\mathcal T_k}
    \left\{\mathbb{W}(h)-\frac12\|h\|_*^2\right\}.
\end{equation}

Our second geometric result connects this density limit to
mixing-measure loss. Consider any sequence satisfying
\[
    \frac{p_{G_m}-p_*}{t_m p_*}
    \rightarrow h\in\mathcal T_k
    \qquad\text{in }L^2(P_*),
    \qquad t_m\downarrow0.
\]
After passing to a subsequence and labeling the fitted atoms, let
$L_i$ contain those converging to $\theta_i^*$, for $i=1,\ldots,p$,
and let $L_0$ contain the remaining, nonlocal atoms. Define the
intermediate mixing measure
\[
    \overline G_m
    :=
    \sum_{i=1}^{p}
    \left(\sum_{j\in L_i}\pi_{j,m}\right)\delta_{\theta_i^*}
    +\sum_{j\in L_0}\pi_{j,m}\delta_{\theta_{j,m}},
\]
which moves each local atom to its corresponding true center while
retaining the nonlocal atoms. We prove the decomposition
\[
    W_2^2(G_m,G_*)
    =
    A_m+W_2^2(\overline G_m,G_*)+o(t_m),
    \qquad
    A_m:=
    \sum_{i=1}^{p}\sum_{j\in L_i}
    \pi_{j,m}\|\theta_{j,m}-\theta_i^*\|^2.
\]
The tangent representation determines both leading terms:
\[
    \frac{A_m}{t_m}
    \rightarrow
    \sum_{i=1}^{p}\operatorname{tr}(C_i),
    \qquad
    \frac{W_2^2(\overline G_m,G_*)}{t_m}
    \rightarrow
   \mathfrak M_*(\alpha,\nu),
\]
where $\mathcal C(\alpha,\nu)$ is the value of a finite
optimal-transport problem for the residual mass. Consequently,
\begin{equation}
\label{eqn:sequential_order_expansion_prelim_form}
    \frac{W_2^2(G_m,G_*)}{t_m}
    \rightarrow
    \Omega(h)
    :=
    \sum_{i=1}^{p}\operatorname{tr}(C_i)
    +\mathfrak M_*(\alpha,\nu).
\end{equation}
    
Crucially, this limit depends only on the relative-density direction
$h$, not on the particular sequence of atoms and weights that
realizes it. We call this as \emph{sequential first-order expansion}. 

Combining the likelihood argmax theorem \eqref{eqn:likelihood_argmax} with this sequential
transportation expansion \eqref{eqn:sequential_order_expansion_prelim_form} gives the joint weak convergence
\[
    \left(
        \widehat h_n,\sqrt n\,W_2^2(\widehat G_n,G_*)
    \right)
    \rightsquigarrow
    \left(
        H^\dagger,\Omega(H^\dagger)
    \right).
\]
In particular,
\[
    n^{1/4}W_2(\widehat G_n,G_*)
    \rightsquigarrow
    \bigl\{\Omega(H^\dagger)\bigr\}^{1/2}.
\]
This result identifies the generally non-Gaussian leading
distribution of the mixing-measure error, extending the analysis
beyond its convergence rate. 

Although this weak limit provides a distributional characterization
of the estimation error, its dependence on an implicit Gaussian
optimization problem and the underlying tangent geometry makes
direct implementation difficult. We therefore develop an
$m$-out-of-$n$ parametric bootstrap and prove its consistency when
$m\to\infty$ and $m/n\to0$, enabling inference without explicitly
constructing the limiting law. We assess its finite-sample performance
in a toy example and illustrate its use with flow-cytometry data.

\subsection{Organization of the article}\label{subsec:organization}

Section~\ref{sec:model} introduces the model, assumptions, and Gaussian location example. Section~\ref{sec:densitylimit} characterizes the attainable likelihood directions and their Gaussian limit. Section~\ref{sec:wassersteinlimit} establishes the sequential first-order Wasserstein expansion and estimator-level weak law. Section~\ref{sec:bootstrap} develops bootstrap inference. Sections~\ref{sec:exp} and present numerical studies and a flow-cytometry illustration, and Section~\ref{sec:discussion} concludes. Detailed proofs appear in the appendices.

\endgroup

\section{Notations and preliminaries}\label{sec:model}

This section sets out the notation and assumptions used throughout the paper. 

For a random element $Z$, let $\mathcal L(Z)$ denote its
probability law. We write $\rightsquigarrow$ and
$\xrightarrow{d}$ interchangeably for convergence in
distribution, and $\xrightarrow{P}$ for convergence
in probability. The notation $O_P(1)$ denotes boundedness in probability,
while $o_P(1)$ denotes convergence in probability to zero.
Let $\mathbb N=\{1,2,\ldots\}$ and
$\mathbb N_0=\mathbb N\cup\{0\}$.
For $a\in\mathbb R$, write $a^+=\max\{a,0\}$ and
$a^-=\max\{-a,0\}$.
The indicator of a set $A$ is denoted by $\mathbf 1_A$. For vectors, $\|\cdot\|$ denotes the Euclidean norm,
and $v^\top$ denotes the transpose of $v$.
For a matrix $A$, $\operatorname{tr}(A)$ and
$\operatorname{rank}(A)$ denote its trace and rank.
We write $\mathbb S_+^d$ for the set of symmetric
positive-semidefinite $d\times d$ matrices. For nonnegative quantities $a$ and $b$, the notation
$a\lesssim b$ means $a\le Cb$ for a constant $C>0$
independent of the relevant asymptotic index.
We write $a\asymp b$ when both $a\lesssim b$ and
$b\lesssim a$ hold.
% The analysis takes place in two spaces. A relative $L^2(P_*)$ norm measures the mixture-density perturbations entering the likelihood expansion. Quadratic Wasserstein distance measures error in the mixing measure. The assumptions below are designed to control these quantities uniformly over the fixed-order model, including sequences of fitted measures that contain vanishing-weight atoms.

\subsection{Mixing measures and likelihood} In this paper, we consider a compact and convex parameter space
$\Theta\subset\RR^d$ and fix an integer $k\geq1$. Let
\[
\cO_k(\Theta)
=
\left\{
G=\sum_{j=1}^m\pi_j\delta_{\theta_j}:
1\leq m\leq k,\ \pi_j\geq0,\
\sum_{j=1}^m\pi_j=1,\ \theta_j\in\Theta
\right\}
\]
denote the space of mixing measures with at most $k$ atoms. For compactness arguments, we may represent measures with fewer than $k$ atoms by adding zero-weight labels without changing the underlying measure. Let $f_\theta=f(\cdot\mid\theta)$ be a density relative to a dominating measure $\mu$ on an observation space
$\mathcal X\subset\RR^{d'}$. For each $G\in\cO_k(\Theta)$, define
\[
p_G(x)=\int_\Theta f_\theta(x)\,dG(\theta),
\qquad
\mathcal P_k(\Theta)=\{p_G:G\in\cO_k(\Theta)\}.
\]
Throughout the paper, unless otherwise stated, we assume that
$X_1,\ldots,X_n\stackrel{\mathrm{iid}}{\sim}P_*$, where $P_*$
has density $p_*=p_{G_*}$ and
\[
    G_*=\sum_{i=1}^p\pi_i^*\delta_{\theta_i^*}\in\cO_p(\Theta),
    \qquad 1\leq p<k.
\]
The true atoms $\theta_i^*$ are pairwise distinct and belong to
$\operatorname{int}(\Theta)$, and their weights $\pi_i^*$ are
strictly positive. We focus on overspecified models, hence the
restriction $k>p$. For measurable $g$, write $P_*g=\int_{\mathcal X}g\,dP_*$
whenever the integral exists. For $g,h\in L^2(P_*)$, define
\[
    \langle g,h\rangle_*=P_*(gh),
    \qquad
    \|g\|_*=(P_*g^2)^{1/2}.
\]
Let $P_n=n^{-1}\sum_{r=1}^n\delta_{X_r}$ be the empirical measure,
so that $P_ng=n^{-1}\sum_{r=1}^n g(X_r)$. For measurable $g$
with $P_*|g|<\infty$, define the empirical process by
\[
    \mathbb G_n g
    :=
    \sqrt n\,(P_n-P_*)g
    =
    \frac1{\sqrt n}\sum_{r=1}^n
    \{g(X_r)-P_*g\}.
\]
Our main object of study is the maximum likelihood estimator,
defined as a measurable global maximizer
\[
\widehat G_n\in
\argmax_{G\in\cO_k(\Theta)}\ell_n(G),
\qquad
\ell_n(G)=\sum_{r=1}^n\log p_G(X_r).
\]
The maximizer need not be unique. Under the conditions of
the measurable selection theorem in \citet{rockafellar},
we may choose a measurable maximizer. To quantify the local density error, define $\rho(G)
=
\left\|
{(p_G-p_*)}/{p_*}
\right\|_*$.

\subsection{Assumptions} We state the assumptions used to establish consistency, expand the
likelihood, and characterize the local density directions.
% We separate the conditions needed for consistency, likelihood expansion and identification of local directions. Assumptions~\ref{ass:geometry} and~\ref{ass:smooth} control parameter geometry and the relative derivatives of the kernel. Assumption~\ref{ass:consistency} gives global likelihood consistency; Assumption~\ref{ass:strongid} prevents second-order derivative cancellations that would change the local inverse scale. Each assumption is imposed at the fixed fitted order $k$.

% \begin{assumption}[Compact parameter geometry and separated truth]\label{ass:geometry}
% $\Theta$ is compact and convex.  The true atoms $\theta_1^*,\ldots,\theta_p^*$ are pairwise distinct and belong to $\operatorname{int}(\Theta)$.  Moreover, a$\pi_{\min}:=\min_{1\le i\le p}\pi_i^*>0$.
% \end{assumption}

% Interior support points permit the local $O(\sqrt t)$ splits used in the recovery-sequence argument. At a boundary point the attainable local moments would depend on the tangent geometry of $\Theta$; the present characterization excludes that separate complication.

\begin{assumption}[Relative $C^3$ smoothness and a common envelope]\label{ass:smooth}
For $\mu$-almost every $x$, the map $\theta\mapsto f_\theta(x)$ is three times continuously differentiable on an open neighborhood of $\Theta$, and $f_\theta(x)>0$ for every $\theta\in\Theta$. Furthermore, there exist $\eta>0$ and a measurable function $F$ such that $P_*F^{2+\eta}<\infty$ and, for $r=0,1,2,3$,
\[
\sup_{\theta\in\Theta}
\frac{\|\nabla_\theta^r f_\theta(x)\|}{p_*(x)}
\le F(x)
\qquad\text{for }P_*\text{-almost every }x, \]
where tensor norms may be any fixed equivalent norms.
\end{assumption}

Assumption~\ref{ass:smooth} controls the kernel and its first three
derivatives relative to $p_*$ through a common
$L^{2+\eta}(P_*)$ envelope. The first two derivative bounds control the linear and quadratic
terms in the Taylor expansion of $f_\theta$ around each true atom,
while the third-derivative bound controls the remainder after
normalization by $p_*$. This integrability
condition allows the relative quantities to be unbounded in $x$,
accommodating the tail behavior of mixture kernels without requiring
a uniformly bounded density ratio.

% The derivative bound is relative to $p_*$ rather than uniform in the observation $x$. It therefore provides common $L^{2+\eta}(P_*)$ control of the relative scores and Hessians; while the third derivative yields the required remainders. This formulation permits the likelihood expansion under tail behavior typical of mixture kernels without assuming a uniformly bounded density ratio.

\begin{assumption}[Likelihood consistency and ordinary identifiability]\label{ass:consistency}
The class $\{\log p_G:G\in\cO_k(\Theta)\}$ is $P_*$-Glivenko--Cantelli with an integrable envelope.  In addition,
\[
p_G=p_*\quad \mu\text{-almost everywhere}
\quad\Longrightarrow\quad
G=G_*
\qquad(G\in\cO_k(\Theta)).
\]
Equivalently, $P_*\log p_G<P_*\log p_*$ for every $G\ne G_*$. 
\end{assumption}

% This assumption yields consistency of $\widehat G_n$ in $W_2$. Identification of the local scale additionally requires excluding derivative cancellations caused by coalescing fitted atoms.
The Glivenko--Cantelli condition ensures uniform convergence of the
normalized log-likelihood $\frac{1}{n}\ell_n(G)$ to its population criterion
$G\mapsto P_*\log p_G$, while identifiability makes $G_*$ its unique
maximizer through the Kullback--Leibler inequality. Together with
compactness of $\cO_k(\Theta)$ and continuity of the population
criterion, these properties yield
$W_2(\widehat G_n,G_*)\to0$ in probability by the standard
consistency theorem for $M$-estimators
\citep[Theorem~5.7]{vdV1998}.

\begin{assumption}[Second-order strong identifiability up to the fitted order]\label{ass:strongid}
For every collection of pairwise distinct points $\xi_1,\ldots,\xi_m\in\Theta$ with $m\le2k$, the identity
\[
\sum_{\ell=1}^{m}
\left[
 a_\ell f_{\xi_\ell}
 +b_\ell^\top\nabla f_{\xi_\ell}
 +\frac12\tr(C_\ell\nabla^2 f_{\xi_\ell})
\right]
=0
\quad \mu\text{-almost everywhere}
\]
with $a_\ell\in\RR$, $b_\ell\in\RR^d$, and symmetric $C_\ell\in\RR^{d\times d}$ implies
\[
a_\ell=0,\qquad b_\ell=0,\qquad C_\ell=0
\quad\text{for every }\ell.
\]
\end{assumption}

% In the compactification of normalized density secants, one may encounter derivatives at the $p$ true support points and at up to $k$ limiting fitted locations. Their number is bounded by $p+k\le2k$. The stated finite-order condition is therefore sufficient for the inverse inequality and uniqueness of tangent coordinates. It is a version of the strong-identifiability conditions used in mixture-rate analysis \citep{HoNguyen2016,HeinrichKahn2018}.

Assumption~\ref{ass:strongid} rules out cancellations among kernels
and their first two derivatives, providing a standard basis for
establishing inverse inequalities in mixture estimation
\citep{HoNguyen2016,HeinrichKahn2018}. It also ensures that each
$h\in\mathcal T_k$ has a unique representation in
\eqref{eqn:dung_unique_presentation_of_h}: subtracting two such
representations yields an identity involving at most $2k$ distinct
locations, so the assumption forces the corresponding coefficients
to agree.

\subsection{A representative example: Gaussian location mixtures with known covariance}\label{subsec:gaussianexample}

% A principal example is the Gaussian location family $f_\theta(x)=\varphi_\Sigma(x-\theta)$,
% where $\Sigma\succ0$ is known,
% with $\Theta\subset\RR^d$ compact and convex. Its verification is explicit and does not require Gaussian assumptions elsewhere in the theory.
The preceding assumptions are standard in the analysis of strongly
identifiable finite mixtures and hold for the Gaussian location
family with known covariance. The following corollary, whose proof is defered to Appendix~\ref{sec:auxiliary_results}, verifies them
for $f_\theta(x)=\varphi_\Sigma(x-\theta)$, where $\varphi_\Sigma$ denotes the density of zero-mean Gaussian vector with covariance matrix $\Sigma$, $\Sigma\succ0$ is
known and $\Theta\subset\RR^d$ is compact and convex.
\begin{corollary}
\label{cor:gaussian}
% Let $\Theta\subset\RR^d$ be compact and convex, and suppose that the true support points are distinct elements of $\operatorname{int}(\Theta)$.  
For the Gaussian location family with known positive-definite covariance matrix $\Sigma$, Assumptions~\ref{ass:smooth}--\ref{ass:strongid} hold. 

% Consequently, all weak-limit and bootstrap results in Sections~\ref{sec:densitylimit}--\ref{sec:bootstrap} apply to the global MLE over $\cO_k(\Theta)$.
\end{corollary}

% The proof of Corollary~\ref{cor:gaussian} illustrates the two ingredients. Gaussian kernel derivatives are polynomial multiples of a Gaussian density, giving the requisite relative envelopes under compactness. For strong identifiability, Fourier transformation converts a finite derivative relation into
% \[
% \sum_j e^{-it^\top\xi_j}P_j(t)=0,
% \]
% where the $P_j$ are polynomials of degree at most two. Linear independence of exponential-polynomial terms at distinct locations then gives the claimed derivative independence. Appendix~\ref{sec:gaussianproof} verifies the full assumptions, including the likelihood envelope.

\begin{remark}[Gaussian mixtures with unknown component scales]
\label{remark:non-scale_location}
The known-covariance restriction applies only to the Gaussian
example, not to the general theory. The multivariate Gaussian
density $f(x\mid\theta,\Sigma)$ satisfies the heat equation $\nabla_\theta^2 f(x\mid\theta,\Sigma)
    =
    2\nabla_\Sigma f(x\mid\theta,\Sigma)$,
where $\nabla_\Sigma$ denotes the covariance gradient with
respect to the Frobenius inner product on symmetric matrices.
Thus, first-order covariance perturbations can cancel
second-order location perturbations, violating second-order
strong identifiability. These derivative interactions lead to
higher-order estimation regimes \citep{HoNguyen2016Weak}
and require a separate weak-limit analysis.
\end{remark}

\section{Likelihood geometry and the density weak limit}\label{sec:densitylimit}

% Maximum likelihood depends on a mixing measure through its induced density. Accordingly, the first limit theorem describes the density perturbation selected by the global MLE. Its state space is not a linear tangent space: the fixed number of fitted atoms restricts the possible second-order displacements, and vanishing-weight atoms can contribute directions whose locations are nonlocal. We characterize the attainable directions first and then derive the Gaussian likelihood limit over this set.

In this section, we characterize the \textit{attainable tangent set} $\cT_k$
and illustrate how the finite component budget constrains local
splitting and vanishing nonlocal components. We then establish
the weak limit of the relative-density maximum likelihood
estimator $\sqrt n\,\frac{p_{\widehat G_n}-p_*}{p_*}$ as the maximizer of a Gaussian quadratic criterion
over $\cT_k$.
\subsection{The attainable tangent set and the finite component budget}\label{subsec:tangentconfig}
% \subsection{The attainable second-order set and the finite component budget}
% \label{subsec:tangentconfig}

We begin by defining the \textit{attainable tangent set} of relative-density directions attainable
by mixing measures with at most $k$ atoms:
\begin{equation}\label{eq:Tkdef}
    \cT_k
    :=
    \left\{
        h\in L^2(P_*):
        \exists\,t_m\downarrow0,\ G_m\in\cO_k(\Theta)
        \text{ such that }
        \left\|
            \frac{p_{G_m}-p_*}{t_m p_*}-h
        \right\|_*\to0
    \right\}.
\end{equation}
Taking $t_m=n^{-1/2}$ gives the scaling relevant to the
relative-density maximum likelihood estimator. Every element of
$\cT_k$ has mean zero under $P_*$, since
$$P_*\left(\frac{p_{G_m}-p_*}{t_m p_*}\right)=0.$$ The following proposition characterizes the elements of $\cT_k$ through
changes in cluster mass, local first and second moments, and
vanishing nonlocal components.

\begin{proposition}[Characteristics of attainable tangent set]
\label{prop:tangent}
Under Assumptions~\ref{ass:smooth}--\ref{ass:strongid}, a function $h\in L^2(P_*)$ belongs to
$\cT_k$ if and only if it admits the representation
\begin{equation}\label{eq:hgamma}
    h(x)
    =
    \frac{1}{p_*(x)}
    \left[
        \sum_{i=1}^{p}
        \left\{
            \alpha_i f_{\theta_i^*}(x)
            +\beta_i^\top\nabla f_{\theta_i^*}(x)
            +\frac12\tr\left(
                C_i\nabla^2f_{\theta_i^*}(x)
            \right)
        \right\}
        +\int_\Theta f_\theta(x)\,d\nu(\theta)
    \right],
\end{equation}
where $\alpha_i\in\RR$, $\beta_i\in\RR^d$,
$C_i\in\mathbb S_+^d$, and
$\nu=\sum_{\ell=1}^{q}c_\ell\delta_{\eta_\ell}$, here $c_\ell>0$, 
    $\eta_\ell\in
    \Theta\setminus\{\theta_1^*,\ldots,\theta_p^*\}$,
with distinct $\eta_\ell$. The case $q=0$ corresponds to $\nu=0$.
These coordinates satisfy mass conservation,
\begin{equation}\label{eq:masscon}
    \sum_{i=1}^{p}\alpha_i+\nu(\Theta)=0,
\end{equation}
and the finite component constraint,
\begin{equation}\label{eq:budget}
    q+\sum_{i=1}^{p}\rank(C_i)\leq k-p.
\end{equation}
The representation is unique up to permutation of the atoms of $\nu$.
\end{proposition}

We provide the proof of this result in Appendix~\ref{sec:prop_of_attainable_tangent_set}. We call $\gamma=(\alpha,\beta,C,\nu)$ the
\emph{tangent configuration} associated with $h$.
The coefficient $\alpha_i$ describes the signed first-order change
in the mass near $\theta_i^*$, while $\beta_i$ and $C_i$ describe
its scaled local first and second moments. The measure $\nu$
records fitted atoms whose weights vanish but whose limiting
locations lie outside the true support. These mechanisms may
occur within the same sequence, subject to \eqref{eq:masscon} and \eqref{eq:budget}.

Equation~\eqref{eq:masscon} expresses conservation of total probability
mass: the signed changes in the masses of the true clusters must
offset the mass assigned to nonlocal components. In particular,
a nonzero $\nu$ requires a corresponding net decrease in the
total mass of the true clusters. Meanwhile, the component constraint~\eqref{eq:budget} reflects the number of support points
needed to realize each configuration. A local second-moment matrix
of rank $r$ requires at least $r+1$ nearby fitted atoms, consuming
$r$ additional atoms beyond the one needed for the true cluster.
Each support point of $\nu$ consumes one further atom. Their total
cannot exceed the $k-p$ additional components available.
Conversely, local simplex constructions and nonlocal atoms with
weights proportional to $t_m$ realize every configuration satisfying
this constraint. Thus, together with mass conservation, \eqref{eq:budget} is
necessary and sufficient for a configuration $\gamma$ such that $h_{\gamma} \in \cT_k$.

\begin{remark}[Vanishing nonlocal components contribute at the root-$n$ scale]
Consistency of the fitted mixing measure does not require every
fitted atom to approach the true support. A component whose weight
is of order $n^{-1/2}$ may remain separated from the true atoms
and still contribute to the root-$n$ density perturbation.
The coordinate $\nu$ in \eqref{eq:hgamma} captures these
vanishing nonlocal components.
\end{remark}

\begin{remark}[The tangent set need not be linear or convex]
The positivity restrictions $C_i\succeq0$ and $\nu\geq0$, together
with the component constraint \eqref{eq:budget}, imply that $\cT_k$
is generally nonconvex and need not be invariant under multiplication
by $-1$. Consequently, standard projection arguments for closed
convex cones do not apply to the limiting likelihood criterion.
We establish existence and almost sure uniqueness of its maximizer
in Lemma~\ref{lem:unique}.
\end{remark}

\subsection{Example: Tangent and in a one-to-two-component}\label{subsec:one_two}

To illustrate how the finite component budget constrains the
attainable tangent directions, we consider a one-to-two-component
Gaussian location mixture. Let $f_\theta(x)=\varphi(x-\theta)$, where $\varphi$ is the standard
normal density, and let $\Theta\subset\RR$ be compact with $0$ in
its interior. Take $G_*=\delta_0$ and fitted upper order $k=2$.
Define
\[
    u(x):=x,\qquad
    v(x):=\frac12(x^2-1),\qquad
    r_\eta(x):=\exp\{\eta x-\eta^2/2\}-1.
\]
These functions satisfy $u=\frac{f'_0}{f_0}$, $v=\frac{f''_0}{2f_0}$, and
    $r_\eta=\frac{f_\eta-f_0}{f_0}$,
where derivatives are taken with respect to the location parameter.
The attainable tangent set is
\begin{equation}\label{eq:one_two_tangent}
    \cT_2
    =
    \bigl\{bu+cv:b\in\RR,\ c\geq0\bigr\}
    \ \cup\
    \bigl\{
        bu+a r_\eta:
        b\in\RR,\ a\geq0,\ \eta\in\Theta\setminus\{0\}
    \bigr\}.
\end{equation}
The term $bu$ represents a first-order location shift and may
occur in either family. The term $cv$ arises from local splitting of the true atom, whereas
$a r_\eta$ represents the relative-density perturbation obtained
by transferring mass of order $t$ from the true atom at $0$ to
a nonlocal atom at $\eta$.
Since only one additional fitted atom is available, the constraint
$q+\rank(C)\leq1$ prevents nontrivial local splitting and a
nonlocal component from occurring simultaneously.

% Consider a one-dimensional Gaussian location mixture with known unit variance, compact $\Theta$ containing $0$ in its interior, true mixing measure $G_*=\delta_0$, and fitted upper order $k=2$. For $X\sim N(0,1)$ put
% \[
% u(x)=x,\qquad v(x)=\tfrac12(x^2-1),\qquad
% r_\eta(x)=\exp\{\eta x-\eta^2/2\}-1.
% \]
% Here $u=f'_0/f_0$, $v=f''_0/(2f_0)$ and $r_\eta=(f_\eta-f_0)/f_0$. The single redundant atom can be used in exactly one of two ways:
% \begin{equation}\label{eq:one_two_tangent}
% \cT_2=\big\{bu+cv:b\in\RR,\ c\ge0\big\}
% \ \cup\
% \big\{bu+a r_\eta:b\in\RR,\ a\ge0,\ \eta\in\Theta\setminus\{0\}\big\}.
% \end{equation}
% The first family describes an extra atom used to split the true component locally; the second uses that atom for a vanishing-weight component at a separated location. The budget $q+\rank(C)\le1$ excludes $a>0$ and $c>0$ simultaneously.

\subsection{The likelihood-selected Gaussian limit}\label{subsec:densitygauss}

To characterize the direction selected by the likelihood, define $\cH:=\overline{\operatorname{span}(\cT_k)}^{\,L^2(P_*)}$, and let $\mathbb W$ be an isonormal Gaussian process indexed by $\cH$, with covariance $\EE\{\mathbb W(h)\mathbb W(g)\} = \langle h,g\rangle_*$, where $h,g\in\cH$. The limiting local log-likelihood criterion takes the form $Q(h):=\mathbb W(h)-\frac12\|h\|_*^2$. Since the fitted model permits only directions in $\cT_k$, we define
\begin{equation}\label{eq:Hdagger}
    H^\dagger
    :=
    \argmax_{h\in\cT_k}
    \left\{
        \mathbb W(h)-\frac12\|h\|_*^2
    \right\}.
\end{equation}
Lemma~\ref{lem:unique} establishes existence and almost sure
uniqueness of this maximizer. The argument uses a separable version
of $\mathbb W$ on bounded sections of $\cT_k$ and does not require
bounded sample paths on the full Hilbert-space unit ball.

\begin{theorem}[Root-$n$ weak limit of the relative-density MLE]
\label{thm:density}
Under Assumptions~\ref{ass:smooth}--\ref{ass:strongid},
\begin{equation}\label{eq:rhorate}
    \rho(\widehat G_n)
    =
    \left\|
        \frac{p_{\widehat G_n}-p_*}{p_*}
    \right\|_*
    =
    O_{P}(n^{-1/2}),
\end{equation}
and
\begin{equation}\label{eq:densitylimit}
    \widehat h_n
    :=
    \sqrt n\,\frac{p_{\widehat G_n}-p_*}{p_*}
    \rightsquigarrow H^\dagger
    \qquad\text{in }L^2(P_*).
\end{equation}
\end{theorem}

The proof of this theorem is discussed in Appendix~\ref{sec:likelihood_expansion_and_convergence}. Despite the root-$n$ normalization, the limiting relative-density
error in \eqref{eq:densitylimit} need not be Gaussian. Its distribution
is determined by maximizing a Gaussian quadratic criterion over
the generally nonconvex set $\cT_k$.

% The proof first decomposes relative-density perturbations into
% local moment contributions and nonlocal residual components.
% A local inverse inequality and empirical-process bounds establish
% the root-$n$ rate in \eqref{eq:rhorate}. The characterization of
% $\cT_k$ identifies the possible subsequential limits, and a Gaussian
% likelihood argmax theorem yields \eqref{eq:densitylimit}.
% The details are given in
% Appendices~\ref{sec:proof1}--\ref{sec:proof3}.

\begin{proof}[Proof sketch of Theorem \ref{thm:density}]
The starting point is a Taylor expansion of the log-likelihood difference. Writing $h_{n,G}:=\sqrt n\,\frac{p_G-p_*}{p_*}$,
we have
\[
    \ell_n(G)-\ell_n(G_*)
    =
    \sum_{r=1}^{n}
    \log\left(1+\frac{h_{n,G}(X_r)}{\sqrt n}\right).
\]
Whenever the $G$ is sufficiently closed to true mixing measure $G_*$, expansion
of the logarithm $\log(1+x) = x - \frac{1}{2}x^2 + O(x^3)$ for $x>-1$ gives
\begin{equation}\label{eq:sketch_likelihood_expansion}
    \ell_n(G)-\ell_n(G_*)
    =
    \mathbb G_n h_{n,G}
    -\frac12 P_nh_{n,G}^2
    +R_{n,G},
\end{equation}
where we use $P_*h_{n,G}=0$.
We first obtain a rate bound from a one-sided likelihood inequality,
then establish a uniform version of
\eqref{eq:sketch_likelihood_expansion} to identify the weak limit.

\medskip
\noindent\textit{Part I: Consistency and root-$n$ localization.}

\smallskip
\noindent\textit{Step I.1.} Using Lemma \ref{lem:consistency}, we have $W_2(\widehat G_n,G_*) \to 0$. Thus when $n$ sufficiently large, there exists a small
neighborhood $\cU$ of $G_*$ such that $\widehat G_n \in \cU$  with probability tending to one.
Consider the normalized class
\[
    \cS
    :=
    \left\{
        s_G=\frac{p_G-p_*}{\rho(G)p_*}:
        G\in\cU,\ G\ne G_*
    \right\},
\]
whose elements satisfy $P_*s_G=0$ and $P_*s_G^2=1$.
As it is possible to construct a common envelop $F_0$ (Lemma~\ref{lem:donsker}), we choose $K<\infty$ such that $P_*\{F_0^2\ind(F_0>K)\}<\frac14$, then  $$\inf_{s\in\cS}P_*\{s^2\ind(F_0\le K)\}\geq\frac34.$$ In addition, the Glivenko--Cantelli property of the truncated squared class gives
\[
    \sup_{s\in\cS}
    \left|
        (P_n-P_*)\{s^2\ind(F_0\le K)\}
    \right|
    \pto0.
\]
Consequently, the event $\mathcal{E}_n:=
    \left\{
        \inf_{s\in\cS}
        P_n\{s^2\ind(F_0\le K)\}\geq\frac12
    \right\}$ has probability tending to one. 
    
Shrink $\cU$ so that $\rho(G)K\leq1/2$ for $G\in\cU$.
Using $\log(1+y)
    \leq y-\frac14y^2\ind(y\leq1/2)$
    for $y>-1$, and $|s_G|\leq F_0$, we obtain
\[
    \ell_n(G)-\ell_n(G_*)
    \leq
    \sqrt n\,\rho(G)\mathbb G_ns_G
    -\frac14 n\rho(G)^2
        P_n\{s_G^2\ind(F_0\le K)\}.
\]
On $\mathcal E_n$, this yields $\ell_n(G)-\ell_n(G_*)
    \leq
    \sqrt n\,\rho(G)\mathbb G_ns_G
    -c n\rho(G)^2$ for $c=\frac18$.
\smallskip

\noindent\textit{Step I.2.}
Define $Z_n:=\sup_{s\in\cS}|\mathbb G_ns|$.
Since $\ell_n(\widehat G_n)\geq\ell_n(G_*)$, the preceding
inequality implies, on an event whose probability tends to one, $\sqrt n\,\rho(\widehat G_n)\leq c^{-1}Z_n$. In addition, the Donsker property of $\cS$ gives $Z_n=O_{P}(1)$, hence $\rho(\widehat G_n)=O_{P}(n^{-1/2})$. 

\medskip
\noindent\textit{Part II: Weak convergence of the relative-density MLE.}

\smallskip
\noindent\textit{Step II.1.}
Fix $M<\infty$ and restrict attention to
\[
    \cH_{n,M}
    :=
    \left\{
        h_{n,G}:
        G\in\cO_k(\Theta),\ \sqrt n\,\rho(G)\leq M
    \right\}.
\]
For sufficiently large $n$, the corresponding mixing measures
belong to $\cU$, and $h_{n,G}=a_Gs_G$ with $0\leq a_G\leq M$ and $s_G \in \mathcal{S}$.
The Glivenko--Cantelli property of the squared class gives
\[
    \sup_{\sqrt n\rho(G)\leq M}
    |(P_n-P_*)h_{n,G}^2|
    \leq
    M^2\sup_{s\in\cS}|(P_n-P_*)s^2|
    =o_{P}(1).
\]
Furthermore, we can control the remainder $R_{n,G}$ in equation \eqref{eq:sketch_likelihood_expansion} as $\sup_{\sqrt n\rho(G)\le M}|R_{n,G}|=o_P(1)$ using the common envelop $F_0$. Together, these bounds control the Taylor remainder uniformly and
give
\[
    \sup_{\sqrt n\rho(G)\leq M}
    \left|
        \ell_n(G)-\ell_n(G_*)
        -\mathbb G_nh_{n,G}
        +\frac12\|h_{n,G}\|_*^2
    \right|
    \pto0.
\]
Thus, for $G$ such that $\sqrt n\,\rho(G)\leq M$, the likelihood discrepancy $\ell_n(G) - \ell_n(G_*)$ can be written as a function of $h \in \mathcal{H}_{n,M}$:
\[
    Q_n(h)
    =
    \mathbb G_nh-\frac12\|h\|_*^2+r_n(h),
    \qquad
    \sup_{h\in\cH_{n,M}}|r_n(h)|=o_{P}(1).
\]

\noindent\textit{Step II.2.}
The compact sets $\cH_{n,M}$ and
$A_M:=\cT_k\cap\{\|h\|_*\leq M\}$ lie in a common compact
set $K_M$. As established in the detailed proof, $d_H(\cH_{n,M},A_M)\to0$.
Define $Z_n(h):=\mathbb G_n h-\frac12\|h\|_*^2$,
for $h\in K_M$. The Donsker lemma gives $Z_n\rightsquigarrow Q$
in $\ell^\infty(K_M)$,
where $Q$ has almost surely continuous paths on $K_M$. The uniform likelihood expansion yields
\[
Q_n(h)=Z_n(h)+r_n(h),
\qquad
\Delta_{n,M}:=
\sup_{h\in\cH_{n,M}}|r_n(h)|=o_P(1).
\]
Consequently, any restricted likelihood maximizer
$\widehat h_{n,M}\in\argmax_{h\in\cH_{n,M}}Q_n(h)$ satisfies
\[
Z_n(\widehat h_{n,M})
\geq \sup_{h\in\cH_{n,M}}Z_n(h)-2\Delta_{n,M}.
\]
By Lemma~\ref{lem:unique}, the restriction of $Q$ to $A_M$ admits an
almost surely unique maximizer $H_M^\dagger$.
The argmax theorem for varying domains
\citep[Theorem~1]{cox2022generalized} therefore gives $\widehat h_{n,M}\rightsquigarrow H_M^\dagger$
in $L^2(P_*)$.

\smallskip
\noindent\textit{Step II.3.}
Equation~\eqref{eq:rhorate} gives
$\|\widehat h_n\|_*=O_P(1)$, and hence $\lim_{M\to\infty}\limsup_{n\to\infty}
\mathbb P\{\|\widehat h_n\|_*>M\}=0$. Choose $\widehat h_{n,M}$ to agree with $\widehat h_n$
whenever $\|\widehat h_n\|_*\leq M$. Moreover,
$H_M^\dagger=H^\dagger$ whenever $\|H^\dagger\|_*\leq M$,
by uniqueness of the global maximizer. Since
$\|H^\dagger\|_*<\infty$ almost surely, a bounded-Lipschitz
comparison, first letting $n\to\infty$ and then $M\to\infty$,
yields
\[
\widehat h_n\rightsquigarrow H^\dagger
\quad\text{in }L^2(P_*).
\]
\end{proof}
The density limit alone does not determine the asymptotic behavior
of the fitted component parameters or their Wasserstein loss.
For example, local splitting at distance $n^{-1/4}$ and mass changes
of order $n^{-1/2}$ both contribute at the root-$n$ density scale,
but affect transportation costs differently. Connecting
Theorem~\ref{thm:density} to mixing-measure loss therefore requires
an additional deterministic expansion.

\section{Sequential expansion and weak convergence of Wasserstein error}\label{sec:wassersteinlimit}
% In the next section, we show that the leading transportation cost
% depends only on the limiting density direction, independently of
% the sequence of fitted mixing measures that realizes it.

% The likelihood limit selects a direction $h\in\cT_k$ in relative-density space. To obtain a weak limit for the mixing-measure error, we must identify the leading quadratic transportation cost of any sequence realizing $h$. This cost has two parts: the local second moments of fitted atoms converging to true support points, and the transport of residual mass between separated locations. The main deterministic result is a sequential expansion, not a derivative along a single parameterized curve.
To connect the density limit established in
Section~\ref{sec:densitylimit} to the Wasserstein estimation error
$W_2(\widehat G_n,G_*)$, we establish a sequential expansion of
the squared Wasserstein distance. Specifically, whenever
$t_m\downarrow0$ and $\frac{p_{G_m}-p_*}{t_m p_*}
    \to h\in\cT_k$ in $L^2(P_*)$,
we show that
\[
    W_2^2(G_m,G_*)=t_m\Omega(h)+o(t_m).
\]
The coefficient $\Omega(h)$ depends only on the limiting
relative-density direction $h$, not on the sequence realizing it.
It combines the contribution of local second moments around
the true atoms with the transportation cost of residual mass
between separated locations. This expansion applies to every
sequence attaining $h$, rather than only to a prescribed
parameterized path, and yields the Wasserstein weak limit
when combined with Theorem~\ref{thm:density}.
\subsection{First-order transportation functional and sequential expansion}\label{subsec:transportfunctional}

We construct the transportation functional by separating local
spreading from residual mass transport. Let
$\gamma=(\alpha,\beta,C,\nu)$ be the unique tangent configuration
associated with $h_\gamma\in\cT_k$, where $\nu=\sum_{\ell=1}^{q}c_\ell\delta_{\eta_\ell}.$  Consider a sequence satisfying
\[
    \frac{p_{G_m}-p_*}{t_mp_*}\longrightarrow h_\gamma
    \qquad\text{in }L^2(P_*),
    \qquad t_m\downarrow0.
\]
Our goal is to determine the limit of the normalized squared
Wasserstein distance $W_2^2(G_m,G_*)/t_m$. We proceed in two steps.

\medskip
\noindent\textit{Step 1: Decomposition of $W_2^2(G_m,G_*)$.}
After passing to a subsequence and labeling the fitted atoms,
let $L_i$ contain those approaching $\theta_i^*$, and let $L_0$
contain the remaining atoms. Define the intermediate measure $\overline G_m$ by moving each local cluster of $G_m$ onto its true support point while retaining the nonlocal atoms: 
\[
    M_{i,m}:=\sum_{j\in L_i}\pi_{j,m},
    \qquad
    \overline G_m
    :=
    \sum_{i=1}^{p}M_{i,m}\delta_{\theta_i^*}
    +\underbrace{\sum_{j\in L_0}\pi_{j,m}\delta_{\theta_{j,m}}}_{:=\nu_m}
\]
The cost of this local collapse is $A_m
    :=
    \sum_{i=1}^{p}\sum_{j\in L_i}
    \pi_{j,m}\|\theta_{j,m}-\theta_i^*\|^2$. 
Using the intermediate measure $\overline G_m$, we can show in Lemma~\ref{lem:collapse} the decomposition of $W_2^2(G_m,G_*)$ into quadratic local collapse $A_m$ and the transportation from $\overline G_m$ to $G_*$: 
\begin{equation}\label{eq:collapse_main}
    W_2^2(G_m,G_*)
    =
    A_m+W_2^2(\overline G_m,G_*)+o(t_m).
\end{equation}
% To obtain this decomposition, we compare transport plans before
% and after collapse. The inverse inequality bounds the optimal
% cost by $O(t_m)$, so separation of the true atoms limits
% cross-cluster transport to $O(t_m)$ mass. Since the local
% displacements tend to zero, the cross terms in the quadratic
% cost are $o(t_m)$.

The term $A_m/t_m$ converges directly to
$\sum_{i=1}^{p}\operatorname{tr}(C_i)$, the sum of the traces of
the local second-moment limits.

\medskip
\noindent\textit{Step 2: Asymptotic behavior of $W_2^2(\overline G_m,G_*)$.}
It remains to determine the limit of
$W_2^2(\overline G_m,G_*)/t_m$. For sufficiently large $m$,
this cost accounts for cluster mass discrepancies and nonlocal
atoms. Put
$d_{ij}^2:=\|\theta_i^*-\theta_j^*\|^2$, let $x_{ij}\geq0$ denote the first-order mass transported from
true support point $i$ to true support point $j$, and let
$y_{\ell i}\geq0$ denote the first-order mass transported from
$\eta_\ell$ to $\theta_i^*$. Define
$\mathfrak M_*(\alpha,\nu)$ as the optimal value of
\begin{equation}\label{eq:LPobj}
    \min_{x_{ij},y_{\ell i}\geq0}
    \left\{
        \sum_{i\ne j}d_{ij}^2x_{ij}
        +
        \sum_{\ell=1}^{q}\sum_{i=1}^{p}
        \|\eta_\ell-\theta_i^*\|^2y_{\ell i}
    \right\},
\end{equation}
subject to
\begin{equation}\label{eq:LPnonlocal}
    \sum_{i=1}^{p}y_{\ell i}=c_\ell,
    \qquad \ell=1,\ldots,q,
\end{equation}
and
\begin{equation}\label{eq:LPbalance}
    \sum_{j\ne i}(x_{ij}-x_{ji})
    =
    \alpha_i+\sum_{\ell=1}^{q}y_{\ell i},
    \qquad i=1,\ldots,p.
\end{equation}
Constraint~\eqref{eq:LPnonlocal} assigns all nonlocal mass to
the true support. Constraint~\eqref{eq:LPbalance} balances mass
at each true atom: its net outgoing mass equals its initial
surplus $\alpha_i$ plus the mass received from nonlocal atoms.

The tangent configuration describes precisely the limiting inputs
to this finite transport problem. In particular, the scaled local mass discrepancies satisfy ${(M_{i,m}-\pi_i^*)}/{t_m}\rightarrow\alpha_i$  for $i=1,\ldots,p.$ The nonlocal mass is represented by $\nu_m:=\sum_{j\in L_0}\pi_{j,m}\delta_{\theta_{j,m}}$,
and its scaled measure satisfies $t_m^{-1}\nu_m\Rightarrow\nu$,
where the convergence is weak convergence of finite measures. Stability of the finite transportation problem
(Lemma~\ref{lem:transportstability}) therefore yields
\[
    \frac{W_2^2(\overline G_m,G_*)}{t_m}
    \longrightarrow\mathfrak M_*(\alpha,\nu).
\]
Combining these two limits $A_m/t_m$ and ${W_2^2(\overline G_m,G_*)}/{t_m}$ into \eqref{eq:collapse_main}, we define
the first-order transportation functional by
\begin{equation}\label{eq:Omega}
    \Omega(h_\gamma)
    :=
    \sum_{i=1}^{p}\tr(C_i)
    +\mathfrak M_*(\alpha,\nu),
\end{equation}
where we denote by convention $\Omega(0) = 0$. Uniqueness of the tangent coordinates makes $\Omega$ a
well-defined functional on $\cT_k$. The first-moment coordinates
$\beta_i$ do not contribute directly: a displacement of order
$t_m$ of an atom with positive limiting mass has squared
transportation cost of order $t_m^2$. They nevertheless enter
the likelihood criterion and may affect the distribution of
$\Omega(H^\dagger)$ through the selected direction $H^\dagger$.

The following proposition states the resulting expansion.
Although the construction uses classified subsequences,
uniqueness of the tangent coordinates ensures that every such
subsequence yields the same limit.

\begin{proposition}[Sequential first-order expansion of the squared Wasserstein distance]
\label{prop:transportmain}
Let $t_m\downarrow0$ and $G_m\in\cO_k(\Theta)$ satisfy
\[
    \frac{p_{G_m}-p_*}{t_mp_*}
    \longrightarrow h\in\cT_k
    \qquad\text{in }L^2(P_*).
\]
Then
\begin{equation}\label{eq:directionalmain}
    \frac{W_2^2(G_m,G_*)}{t_m}
    \longrightarrow\Omega(h).
\end{equation}
\end{proposition}

We show the proof of Proposition~\ref{prop:transportmain} in Appendix~\ref{app:sequential_expansion}. Proposition~\ref{prop:transportmain} shows that every sequence
attaining the same relative-density direction $h$ satisfies
\[
    W_2^2(G_m,G_*)=t_m\Omega(h)+o(t_m).
\]
We call this a \emph{sequential first-order expansion}:
it is first-order because its leading term is linear in
the density-perturbation scale $t_m$, and sequential because
it holds for every sequence realizing $h$, rather than only
along a prescribed parameterized path.
The coefficient $\Omega(h)$ combines local spreading,
described by $C_i$, with residual mass transport, determined
by $\alpha$ and $\nu$. The latter cost is the optimal value
of a finite linear program in transported masses, with
fixed squared distances as coefficients.

\begin{remark}[The expansion applies beyond maximum likelihood]
The sequential expansion holds for any deterministic sequence
$G_m\in\cO_k(\Theta)$ with a convergent normalized density
perturbation, regardless of how $G_m$ is obtained. Every sequence
attaining the same direction $h$ has the same leading
transportation cost $\Omega(h)$. This property allows us to
combine the expansion with the density weak limit of the MLE
in the next theorem.
\end{remark}

\subsection{Joint weak convergence and the fourth-root law}\label{subsec:w2joint}
Combining the density limit in Theorem~\ref{thm:density} with the
sequential expansion in Proposition~\ref{prop:transportmain}, we
obtain joint weak convergence of the relative-density estimator
and its squared Wasserstein error.

\begin{theorem}[Joint relative-density and Wasserstein weak limit]
\label{thm:wasserstein}
Under Assumptions~\ref{ass:smooth}--\ref{ass:strongid},
\begin{equation}\label{eq:jointlimit}
    \left(
        \widehat h_n,
        \sqrt n\,W_2^2(\widehat G_n,G_*)
    \right)
    \rightsquigarrow
    \left(
        H^\dagger,\Omega(H^\dagger)
    \right)
    \qquad\text{in }L^2(P_*)\times\RR.
\end{equation}
Consequently,
\begin{equation}\label{eq:w2limit}
    n^{1/4}W_2(\widehat G_n,G_*)
    \rightsquigarrow
    \{\Omega(H^\dagger)\}^{1/2},
\end{equation}
and
\begin{equation}\label{eq:w2rate}
    W_2(\widehat G_n,G_*)=O_{P}(n^{-1/4}).
\end{equation}
The limiting distribution in \eqref{eq:w2limit} is generally
non-Gaussian.
\end{theorem}

The proof this theorem is given in Appendix~\ref{sec:maincompletion}. Theorem~\ref{thm:wasserstein} identifies how the density and mixing-measure errors fluctuate jointly. The likelihood criterion selects the random density direction $H^\dagger$, and the functional $\Omega$ translates this direction into the leading squared transportation cost. Thus, the Wasserstein limit is determined by the same random direction that governs the density limit.

A few remarks clarify the relation to existing rate results
and the nature of the limiting distribution. First, inverse inequalities provide a route from density
estimation to mixing-measure convergence. In particular,
the bound $W_2^2(G,G_*)\lesssim\rho(G)$ in Lemma~\ref{lem:w2bound}, combined with
$\rho(\widehat G_n)=O_{P}(n^{-1/2})$, yields the
fourth-root upper rate in \eqref{eq:w2rate}. This bound alone,
however, does not identify a limiting distribution.
Theorem~\ref{thm:wasserstein} follows the same density-to-mixing-measure
route but replaces the inequality with a sequential expansion,
thereby identifying the random leading term
$\Omega(H^\dagger)$.

Second, the Wasserstein limit in \eqref{eq:w2limit} is generally
non-Gaussian, as it is obtained by applying the nonlinear map
$h\mapsto\{\Omega(h)\}^{1/2}$ to the likelihood-selected direction
$H^\dagger$. Non-Gaussian distance limits also arise in regular
models: for a one-dimensional Gaussian location model with
$p=k=1$ and known variance $\sigma^2$,
\[
    \sqrt n\,W_1(\widehat G_n,G_*)
    \rightsquigarrow |Z|,
    \qquad Z\sim N(0,\sigma^2).
\]
Here the root-$n$ distance limit is non-Gaussian, while the
fourth-root limit in \eqref{eq:w2limit} degenerates to zero.

\subsection{Example continued: Wasserstein expansion in a one-to-two-component mixture}
\label{subsec:one_two_transport}

We return to the Gaussian location mixture in
Section~\ref{subsec:one_two} and compute the transportation
functional for its two tangent families. Using the notation
introduced there, consider
\begin{align*}
    G_t^{\mathrm{split}}
    =
    \frac12\delta_{bt+\sqrt{ct}}
    +\frac12\delta_{bt-\sqrt{ct}}, \quad 
    G_t^{\mathrm{remote}}
    =
    (1-at)\delta_{bt}+at\delta_\eta,
\end{align*}
where $b\in\RR$, $c,a\geq0$, and
$\eta\in\Theta\setminus\{0\}$.
Both paths belong to $\cO_2(\Theta)$ for sufficiently small
$t>0$ and satisfy
\[
    \frac{p_{G_t^{\mathrm{split}}}-p_*}{tp_*}
    \longrightarrow bu+cv,
    \qquad
    \frac{p_{G_t^{\mathrm{remote}}}-p_*}{tp_*}
    \longrightarrow bu+a r_\eta
    \qquad\text{in }L^2(P_*).
\]
Since $G_*=\delta_0$, every fitted atom must transport its mass
to $0$. The squared Wasserstein costs are therefore
\begin{align}
    W_2^2(G_t^{\mathrm{split}},G_*)
    =tc+b^2t^2,\quad 
    W_2^2(G_t^{\mathrm{remote}},G_*)
    =at\eta^2+(1-at)b^2t^2.
\end{align}
Dividing by $t$ and taking limits gives
\begin{equation}\label{eq:toy_omega}
    \Omega(bu+cv)=c,
    \qquad
    \Omega(bu+a r_\eta)=a\eta^2.
\end{equation}
The local splitting cost is thus determined by the second-moment
coefficient $c$, while the nonlocal cost is the transferred mass
coefficient $a$ multiplied by the squared distance $\eta^2$.
The location-shift coefficient $b$ contributes only at order
$t^2$ and disappears from the first-order expansion.
Both formulas give zero on the family
$\{bu:b\in\RR\}$.

Although these calculations use explicit paths,
Proposition~\ref{prop:transportmain} ensures that
\eqref{eq:toy_omega} holds for every sequence attaining the
corresponding density direction. The formulas therefore determine
the Wasserstein limit after the likelihood selects $H^\dagger$:
\[
    n^{1/4}W_2(\widehat G_n,G_*)
    \rightsquigarrow
    \begin{cases}
        \sqrt c,
        & \text{if } H^\dagger=bu+cv,\\[2pt]
        \sqrt a\,|\eta|,
        & \text{if } H^\dagger=bu+a r_\eta.
    \end{cases}
\]
The two expressions agree and exactly equal to 0 wherever the families overlap.

The likelihood optimization takes place over the union of
the two families in \eqref{eq:one_two_tangent}. It compares their
optimized Gaussian quadratic criteria, including optimization
over the nonlocal location $\eta$. Consequently, the selected
direction depends on the Gaussian process indexed by these
locations, and its transportation cost cannot be recovered
from the optimized likelihood value alone. This example makes
the transportation functional explicit, while the unrestricted
Gaussian argmax and its induced Wasserstein distribution remain
implicitly characterized.

\section{Bootstrap inference for Wasserstein uncertainty}\label{sec:bootstrap}

% A plug-in simulation of $\{\Omega(H^\dagger)\}^{1/2}$ would require estimating the unknown attainable set and the covariance of the associated Gaussian process. The following bootstrap avoids these quantities. Its validity follows from a direct comparison between the fitted and true experiments at a smaller sample size.

% Theorem~\ref{thm:density} gives $\rho^2(\widehat G_n)=O_P(n^{-1})$. Consequently, an $m$-sample from the fitted density becomes asymptotically indistinguishable, in total variation, from an $m$-sample from the true density when $m=o(n)$. A separate centering argument uses $m^{1/4}W_2(\widehat G_n,G_*)=o_P(1)$. Together these comparisons identify the conditional bootstrap limit.
Although Theorem~\ref{thm:wasserstein} characterizes the limiting distribution of the Wasserstein estimation error, direct simulation is difficult because the tangent set $\cT_k$ has a complex, generally nonconvex structure, and the optimization
criterion involves a Gaussian process indexed by this set. In this section, we develop an $m$-out-of-$n$ parametric bootstrap that approximates this distribution without explicitly constructing the limiting Gaussian process or transportation functional. We prove its consistency when $m\to\infty$ and $m/n\to0$, and establish convergence of bootstrap quantiles and asymptotically correct coverage of Wasserstein confidence regions under a quantile regularity condition. We also show that these regions can be implemented using Monte Carlo bootstrap replicates while retaining their asymptotic coverage as the number of replicates increases.

\subsection{Consistency of the $m$-out-of-$n$ parametric bootstrap}\label{subsec:bootstrapmain}
Let $\mathcal F_n=\sigma(X_1,\ldots,X_n)$ be the $\sigma$-field
generated by the original data, and let $Q_n=P_{\widehat G_n}$
be the fitted sampling distribution. For each sample size $r$,
fix a measurable map $\mathcal A_r$ selecting a global maximum
likelihood estimator over $\cO_k(\Theta)$. Let $m=m_n$ be a
deterministic bootstrap sample size satisfying
\begin{equation}\label{eq:mconditions}
    m\to\infty,
    \qquad m/n\to0.
\end{equation}
Conditional on $\mathcal F_n$, draw
$X_1^*,\ldots,X_m^*\stackrel{\mathrm{iid}}{\sim}Q_n$ and define
$\widehat G_m^*=\mathcal A_m(X_1^*,\ldots,X_m^*)$ as the bootstrap
maximum likelihood estimator. The statistic
\begin{equation}\label{eq:Tstar}
    T_{n,m}^*
    :=m^{1/4}W_2(\widehat G_m^*,\widehat G_n)
\end{equation}
measures the scaled bootstrap estimation error around the
original estimate. For comparison, let
\begin{equation}\label{eq:TnTinf}
    T_n:=n^{1/4}W_2(\widehat G_n,G_*),
    \qquad
    T_\infty:=\{\Omega(H^\dagger)\}^{1/2}
\end{equation}
be the original estimation-error statistic and its weak limit.
Write $\mathcal L^*$ for conditional law given $\mathcal F_n$.
For probability laws $P$ and $Q$ on $\RR$, define the
bounded-Lipschitz distance by
\[
    d_{\mathrm{BL}}(P,Q)
    :=
    \sup_{f\in\mathrm{BL}_1}
    \left|\int f\,dP-\int f\,dQ\right|,
    \qquad
    \]
    where $\mathrm{BL}_1 = \left\{
        f:\RR\to\RR:
        \|f\|_\infty\leq1,\ 
        |f(x)-f(y)|\leq|x-y|,\ \forall x,y\in\RR
    \right\}.$ Although $\Theta$ is compact, the scaled Wasserstein statistic
$T_n=n^{1/4}W_2(\widehat G_n,G_*)$ need not be uniformly bounded
in $n$. We therefore use $d_{\mathrm{BL}}$ to express convergence
of its distribution without imposing additional moment conditions.

\begin{theorem}[Consistency of the $m$-out-of-$n$ parametric bootstrap]
\label{thm:bootstrap}
Under Assumptions~\ref{ass:smooth}--\ref{ass:strongid} and
\eqref{eq:mconditions},
\begin{equation}\label{eq:bootlimit}
    d_{\mathrm{BL}}
    \left(
        \mathcal L^*(T_{n,m}^*),
        \mathcal L(T_\infty)
    \right)
    \pto0.
\end{equation}
Consequently,
\begin{equation}\label{eq:bootfiniteNtarget}
    d_{\mathrm{BL}}
    \left(
        \mathcal L^*(T_{n,m}^*),
        \mathcal L(T_n)
    \right)
    \pto0.
\end{equation}
\end{theorem}

The proof of Theorem~\ref{thm:bootstrap} is provided in
Appendix~\ref{app:bootstrapproofs}. It establishes bootstrap consistency
by comparing the fitted and true sampling distributions at
sample size $m$, then controlling the change of center from
$G_*$ to $\widehat G_n$.

Specifically, let $Y_1,\ldots,Y_m\stackrel{\mathrm{iid}}{\sim}P_*$
be independent of the original data, and define
$\widetilde G_m:=\mathcal A_m(Y_1,\ldots,Y_m)$ and
$T_m^0:=m^{1/4}W_2(\widetilde G_m,G_*)$.
Using the same measurable MLE selection for both samples gives
\begin{equation}\label{eq:explicitBLbound}
\begin{split}
    d_{\mathrm{BL}}
    \left(
        \mathcal L^*(T_{n,m}^*),\mathcal L(T_m^0)
    \right)
    \leq
    m^{1/4}W_2(\widehat G_n,G_*)+
   2\|Q_n^{\otimes m}-P_*^{\otimes m}\|_{\mathrm{TV}}
\end{split}
\end{equation}
The first term controls the centering error through the triangle
inequality and is $(m/n)^{1/4}T_n=o_{ P}(1)$.
while the second controls the total variation distance between
$Q_n^{\otimes m}$ and $P_*^{\otimes m}$. It also vanishes because
$\rho(\widehat G_n)^2=O_{P}(n^{-1})$, $m/n\to0$, and the inequality 
\begin{equation*}
    2\|Q_n^{\otimes m}-P_*^{\otimes m}\|_{\mathrm{TV}}
\le\left[\{1+\rho(\widehat G_n)^2\}^m-1\right]^{1/2}.
\end{equation*}
Finally, $T_m^0\rightsquigarrow T_\infty$ as $m\to\infty$,
which yields \eqref{eq:bootlimit}; comparison with
$T_n\rightsquigarrow T_\infty$ gives \eqref{eq:bootfiniteNtarget}.

The condition $m/n\to0$ provides a sufficient asymptotic rule,
rather than an optimal finite-sample choice of $m$. Smaller
bootstrap samples reduce the two comparison terms but may
increase the finite-$m$ error in approximating the limiting
distribution. Section~\ref{sec:exp} examines this tradeoff.

\subsection{Critical values, confidence regions, and Monte Carlo implementation}
\label{subsec:bootstrapinference}

The bootstrap approximation yields inference for both the
Wasserstein error and its square. Since
$\sqrt m\,W_2^2(\widehat G_m^*,\widehat G_n)=(T_{n,m}^*)^2$,
the continuous mapping theorem gives the following consequence.

\begin{corollary}[Bootstrap validity for the squared Wasserstein statistic]
\label{cor:bootsquared}
Under the conditions of Theorem~\ref{thm:bootstrap},
\begin{equation}\label{eq:bootsquared}
    \mathcal L^*
    \left\{
        \sqrt m\,W_2^2(\widehat G_m^*,\widehat G_n)
    \right\}
    \weakto_P
    \mathcal L\{\Omega(H^\dagger)\}.
\end{equation}
\end{corollary}

To construct confidence regions, let
$F_\infty(t):=P(T_\infty\leq t)$ be the limiting distribution
function and $F_{n,m}^*(t):=P^*(T_{n,m}^*\leq t)$ its conditional
bootstrap counterpart. For $p_0\in(0,1)$, define
\[
    q_{p_0}:=\inf\{t:F_\infty(t)\geq p_0\},
    \qquad
    q_{n,m}^*(p_0):=\inf\{t:F_{n,m}^*(t)\geq p_0\}
\]
as their generalized inverse quantiles. Convergence of the
bootstrap quantile and exact asymptotic coverage require
regularity of the limiting distribution at the chosen quantile.

\begin{corollary}[Bootstrap quantiles and Wasserstein confidence regions]
\label{cor:bootquantile}
Under the conditions of Theorem~\ref{thm:bootstrap}, fix
$\alpha\in(0,1)$ and put $p_0=1-\alpha$. Suppose that
$q_{p_0}:=\inf\{t:F_\infty(t)\geq p_0\}$ is a regular quantile,
meaning that
\begin{equation}\label{eq:quantreg}
    P(T_\infty=q_{p_0})=0,
    \qquad F_\infty(q_{p_0})=p_0,
    \qquad F_\infty(q_{p_0}+\varepsilon)>p_0
    \quad\text{for every }\varepsilon>0.
\end{equation}
Then the bootstrap quantile is consistent,
$q_{n,m}^*(p_0)\pto q_{p_0}$, and its associated critical value
has asymptotically exact coverage:
\begin{equation}\label{eq:coverage}
    P\left\{
        n^{1/4}W_2(\widehat G_n,G_*)
        \leq q_{n,m}^*(1-\alpha)
    \right\}\longrightarrow1-\alpha.
\end{equation}
Consequently, the Wasserstein confidence region
\begin{equation}\label{eq:confball}
    \mathcal C_{n,1-\alpha}
    :=
    \left\{
        G:W_2(G,\widehat G_n)
        \leq n^{-1/4}q_{n,m}^*(1-\alpha)
    \right\}
\end{equation}
satisfies $P\{G_*\in\mathcal C_{n,1-\alpha}\}\to1-\alpha$.
\end{corollary}

Generalized inverse quantiles and their empirical approximations
are standard in statistical practice and can be computed without
assuming that the underlying distribution is continuous
\citep{EmbrechtsHofert2013,HyndmanFan1996}. Condition~\eqref{eq:quantreg}
is needed for the quantile-consistency and exact-coverage
conclusions, rather than for defining or computing the bootstrap
quantile.

For implementation, generate $B=B_n$ conditionally independent
bootstrap replicates $T_{n,m}^{*(1)},\ldots,T_{n,m}^{*(B)}$ and define
\[
    \widehat F_{n,m,B}^*(t)
    :=
    \frac1B\sum_{b=1}^{B}\ind\{T_{n,m}^{*(b)}\leq t\}
\]
as their empirical distribution function. Its generalized inverse,
\[
    \widehat q_{n,m,B}^*(p_0)
    :=
    \inf\{t:\widehat F_{n,m,B}^*(t)\geq p_0\},
\]
provides a Monte Carlo approximation to the conditional bootstrap
quantile.

\begin{proposition}[Monte Carlo approximation of the conditional bootstrap]
\label{prop:MCbootstrap}
Under the conditions of Theorem~\ref{thm:bootstrap}, if
$B=B_n\to\infty$, then
\begin{equation}\label{eq:DKWconv}
    \sup_t
    \left|
        \widehat F_{n,m,B}^*(t)-F_{n,m}^*(t)
    \right|
    \pto0.
\end{equation}
Under the quantile conditions of Corollary~\ref{cor:bootquantile},
\begin{equation}\label{eq:MCquantile}
    \widehat q_{n,m,B}^*(1-\alpha)
    \pto q_{1-\alpha},
\end{equation}
and the implementable confidence region
\begin{equation}\label{eq:MCball}
    \widehat{\mathcal C}_{n,1-\alpha}
    :=
    \left\{
        G:
        W_2(G,\widehat G_n)
        \leq n^{-1/4}\widehat q_{n,m,B}^*(1-\alpha)
    \right\}
\end{equation}
has asymptotic coverage $1-\alpha$.
\end{proposition}

Thus, replacing the conditional bootstrap quantile by its
Monte Carlo estimate preserves asymptotic coverage as
$B\to\infty$. The proofs of both corollaries and
Proposition~\ref{prop:MCbootstrap} are given in
Appendix~\ref{app:bootstrapproofs}.

For nominal coverage $1-\alpha$, Algorithm~\ref{alg:WassersteinCI}
constructs a Wasserstein confidence region without explicitly
constructing the tangent set or the limiting Gaussian process. The coverage guarantee holds when $m\to\infty$, $m/n\to0$,
and $B\to\infty$, at a regular target quantile. Numerical
approximation of the global maximizer and the finite-sample
choice of $m$ remain implementation choices, which we assess
in the following experiments.

\begin{algorithm}[t]
\caption{Constructing a Wasserstein confidence region}
\label{alg:WassersteinCI}
\begin{algorithmic}[1]
\Require Observations $X_1,\ldots,X_n$; fitted upper order $k$;
bootstrap sample size $m<n$; number of replicates $B$;
significance level $\alpha\in(0,1)$.
\Ensure Wasserstein confidence region
$\widehat{\mathcal C}_{n,1-\alpha}$.

\State Fit the global mixture likelihood over $\cO_k(\Theta)$
to obtain $\widehat G_n$.

\For{$b=1,\ldots,B$}
    \State Generate $X_1^{*(b)},\ldots,X_m^{*(b)}
    \stackrel{\mathrm{iid}}{\sim}P_{\widehat G_n}$,
    independently across replicates.
    \State Refit the global likelihood over $\cO_k(\Theta)$
    to obtain $\widehat G_m^{*(b)}$.
    \State Compute
    $T_{n,m}^{*(b)}
    =m^{1/4}W_2(\widehat G_m^{*(b)},\widehat G_n)$.
\EndFor

\State Compute the empirical generalized inverse quantile
\[
    \widehat q_{n,m,B}^*(1-\alpha)
    :=
    \inf\left\{
        t:\frac1B\sum_{b=1}^{B}
        \ind\{T_{n,m}^{*(b)}\leq t\}\geq1-\alpha
    \right\}.
\]
\State Set $r_{n,\alpha}
=n^{-1/4}\widehat q_{n,m,B}^*(1-\alpha)$.
\State \Return
$\widehat{\mathcal C}_{n,1-\alpha}
=\{G\in\cO_k(\Theta):W_2(G,\widehat G_n)\leq r_{n,\alpha}\}$.
\end{algorithmic}
\end{algorithm} 
% \paragraph{Implementing a Wasserstein confidence region.}
% For a nominal coverage probability $1-\alpha$, the procedure can be implemented without constructing the tangent set or the limiting Gaussian process:
% \begin{enumerate}[label=\textbf{\arabic*.},leftmargin=2.6em,itemsep=0.25em]
% \item Fit the global mixture likelihood over $\cO_k(\Theta)$ to obtain $\widehat G_n$.
% \item Choose $m<n$ and, for each of $B$ replicates, generate $m$ observations from $P_{\widehat G_n}$ and refit the likelihood with the same upper order $k$.
% \item Compute $T_{n,m}^{*(b)}=m^{1/4}W_2(\widehat G_m^{*(b)},\widehat G_n)$ and its empirical $(1-\alpha)$ quantile $\widehat q_{n,m,B}^*(1-\alpha)$.
% \item Report the Wasserstein ball around $\widehat G_n$ with radius $n^{-1/4}\widehat q_{n,m,B}^*(1-\alpha)$, as in \eqref{eq:MCball}.
% \end{enumerate}

% The coverage result concerns the joint asymptotic regime $m\to\infty$, $m/n\to0$ and $B\to\infty$, at a regular target quantile. Numerical approximation to a global maximizer and finite-sample selection of $m$ remain implementation choices. We assess both in the following experiments.

\section{Experimental Results}\label{sec:exp}
In this section, we empirically assess the finite-sample implications of our theory through Gaussian location-mixture simulations with known covariance and a flow-cytometry application. Our evaluation focuses on two aspects: (i) the fourth-root scale of mixing-measure error and the influence of transportation geometry and the fitted component budget on its distribution; and (ii) the use of these regions to obtain model-based uncertainty bounds for population-level summaries that are invariant to component labeling.

\begin{figure}[t]
\centering
\includegraphics[width=0.95\textwidth]{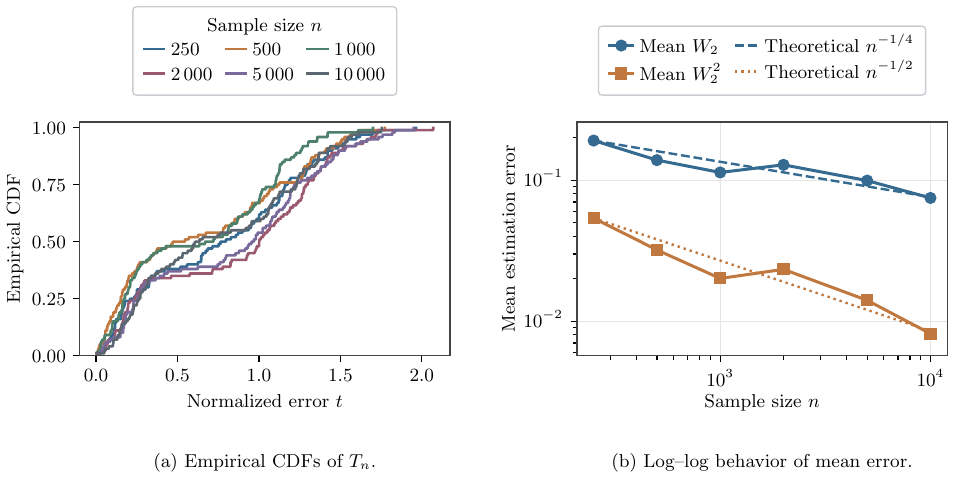}
\caption{Wasserstein estimation error for a one-component Gaussian mixture fitted with at most two components. (a) Empirical CDFs of \(n^{1/4}W_2(\widehat G_n,G^\star)\) across sample sizes. (b) Mean \(W_2\) and \(W_2^2\) errors on logarithmic axes, with reference lines proportional to \(n^{-1/4}\) and \(n^{-1/2}\), respectively.}
\label{fig:simAoverview}
\end{figure}

\subsection{Simulation Studies}
We conduct numerical studies in Gaussian location mixtures with known covariance, beginning by investigating whether the fourth-root normalization stabilizes error distributions across sample sizes. To interpret the errors on this scale, we separate local displacement from mass reallocation and quantify their contributions to Wasserstein loss. Because the component budget constrains the attainable perturbations, we further study how changing this budget affects the error distribution and multivariate splitting geometry. This dependence of the error distribution on the fitted component budget also matters for inference, since upper quantiles of the normalized error determine the radii of Wasserstein confidence regions.

\noindent Throughout the simulated experiments, we keep the covariance fixed at $I_d$ so that fitting varies only the component locations and mixing weights. For each study, we generate observations from a specified true mixing measure $G^\star$ and fit a model allowing at most $k$ components. All centers are restricted to $\Theta=[-5,5]^d$, unless stated otherwise. For each configuration, we repeat the experiment on independently generated datasets. The number of repetitions is reported in the corresponding subsection.

\noindent To approximate the global MLE considered in the theory, we use the expectation--maximization (EM) algorithm \citep{dempster1977maximum} with multiple initializations based on data quantiles, nearby component centers, small component weights, and random starts. At each EM iteration, we restrict the component centers to $\Theta$ and keep the mixing weights positive and normalized to sum to one. Finally, we evaluate the error of the fitted mixing measures using Wasserstein distances, computed by monotone matching in one dimension and by solving a finite transport linear program in two dimensions.

\subsubsection{Scaling of Mixing-Measure Error}\label{subsec:mixing-error-scaling}
We first examine how mixing-measure estimation error varies with sample size in an over-specified Gaussian mixture.
\paragraph{Setup.} We consider a one-dimensional Gaussian location mixture with $G^\star=\delta_0$ and fitted upper order $k=2$. Although the observations come from a single Gaussian, the additional fitted component permits local splitting around the true center. This provides a simple over-specified setting in which to study the estimation error while keeping the sampling distribution fixed. We generate observations from $N(0,1)$ and repeat the experiment on $100$ independent datasets for each $n\in\{250,500,1000,2000,5000,10000\}$. For each fit, we record
\[
T_n=n^{1/4}W_2(\widehat G_n,G^\star),
\qquad
V_n=\sqrt n\,W_2^2(\widehat G_n,G^\star).
\]
The normalized errors allow us to compare distributions across sample sizes, while the corresponding unnormalized errors reveal how the magnitude of the estimation error decreases with $n$.

\paragraph{Results.} In Figure~\ref{fig:simAoverview}(a), the empirical distributions of $T_n$ remain broadly similar across sample sizes after fourth-root normalization, suggesting that this normalization captures a common scale for the estimation error. By comparing the error distributions, we also consider their upper quantiles, which are not captured by the average error alone.

\noindent To evaluate the rate of decay, in Figure~\ref{fig:simAoverview}(b) we plot the mean unnormalized errors against sample size on logarithmic axes and compare these estimates with $-1/4$ and $-1/2$, which is the exponents associated with the theoretical scalings $n^{-1/4}$ for $W_2$ and $n^{-1/2}$ for $W_2^2$, respectively.  The mean errors decrease and generally follow the reference trends. Both theoretical exponents lie within the corresponding $95\%$ Monte Carlo intervals in Table~\ref{tab:simAslopes}, which indicates that the observed mean errors are compatible with these scalings over the simulated sample sizes.

\begin{table}[t]
\centering
\caption{Estimated slopes of log mean Wasserstein error against log sample size. Each mean is computed from $100$ independent datasets. The $95\%$ Monte Carlo intervals are obtained by resampling the replicate errors within each sample size and refitting the slope, using $2{,}000$ bootstrap resamples.}
\vspace{10pt}
\label{tab:simAslopes}
\small
\begin{tabular}{lccc}
\toprule
Error measure & Theoretical exponent & Estimated slope & $95\%$ Monte Carlo interval \\
\midrule
$W_2(\widehat G_n,G^\star)$
& $-1/4$ & $-0.214$ & $[-0.257,-0.168]$ \\
$W_2^2(\widehat G_n,G^\star)$
& $-1/2$ & $-0.451$ & $[-0.515,-0.385]$ \\
\bottomrule
\end{tabular}
\end{table}

\noindent For comparison, fitting a single Gaussian to the same data gives a center estimate equal to the sample mean with probability tending to one, and its Wasserstein error has scale $n^{-1/2}$. This error decreases faster than the $n^{-1/4}$ scale observed when two components are allowed. Consequently, multiplying the one-component error by $n^{1/4}$ would make it converge to zero, whereas the normalized errors from the two-component fits remain on a comparable scale in Figure~\ref{fig:simAoverview}(a). This comparison illustrates how allowing a redundant component can slow mixing-measure recovery, even though the data still come from the same single Gaussian.

\begin{table}[t]
\centering
\caption{Median Wasserstein error and transport costs at $n=5000$ for a balanced two-component truth, fitted with at most three components. Each separation uses $60$ independent datasets. Here $T_n=n^{1/4}W_2(\widehat G_n,G^\star)$, while the transport costs and absolute residual are scaled by $\sqrt n$. Medians are computed separately for each column and do not provide an additive decomposition.}
\vspace{12pt}
\label{tab:simBsep}
\small
\begin{tabular}{ccccc}
\toprule
& & \multicolumn{3}{c}{\textit{Transport costs and residual}} \\
\cmidrule[0.6pt](lr){3-5}
\shortstack{Separation\\$\Delta$}
& \shortstack{Wasserstein error\\$T_n$}
& \shortstack{Local displacement\\$\sqrt n L_n$}
& \shortstack{Mass reallocation\\$\sqrt n M_n$}
& \shortstack{Absolute residual\\$\sqrt n|R_n^{\mathrm{OT}}|$} \\
\midrule
$2$ & $2.183$ & $0.344$ & $4.363$ & $0.159$ \\
$4$ & $2.258$ & $0.875$ & $4.549$ & $0.054$ \\
\bottomrule
\end{tabular}
\end{table}
\subsubsection{Local Displacement and Mass Reallocation}\label{subsec:transport-decomposition}

In preceding experiments, we investigate how quickly mixing-measure error decreases with sample size. We continue by studying the contributions to this error when the true mixture contains more than one component. In this setting, the fitted measure can differ from the truth both in the locations of its atoms and in the mass assigned to each true component. Motivated by the transportation expansion in Proposition~\ref{prop:transportmain}, we separate these contributions to assess how much of the Wasserstein error comes from local displacement and how much comes from mass reallocation between true components.

\paragraph{Setup.} We consider two-point true mixing measures with separations $\Delta\in\{1,2,4\}$, using both balanced and unequal mixing weights, and fit mixtures with at most three components. This allows the fitted measure to split a true component while also assigning different total masses to the two true components. For each configuration, we generate $60$ independent datasets at sample sizes ranging from $500$ to $5000$.

\paragraph{Error decomposition.} To measure local displacement, we assign each fitted center $\widehat\theta_j$ to a nearest true center $\theta^\star_{\iota(j)}$ and compute $L_n=\sum_j\widehat\pi_j\bigl(\widehat\theta_j-\theta^\star_{\iota(j)}\bigr)^2$. This sum gives the total transportation cost of moving the fitted atoms onto their assigned true centers. After this movement, the resulting measure is
\[
\widehat G_n^c=\sum_i\widehat M_i\delta_{\theta_i^\star},
\qquad
\widehat M_i=\sum_{j:\,\iota(j)=i}\widehat\pi_j.
\]
The centers of $\widehat G_n^c$ now agree with those of $G^\star$, but their weights may still differ. To quantify this difference in weights, we compute $M_n=W_2^2(\widehat G_n^c,G^\star)$. For a two-point truth, this cost equals the mass that must be transferred between the two true centers multiplied by their squared separation. From this perspective, $M_n$ measures the cost of reallocating excess mass from one true component to the other.

\noindent The costs $L_n$ and $M_n$ describe two steps, namely moving fitted atoms onto their nearest true centers and then reallocating mass between the true centers. The optimal coupling can instead transport mass directly from the fitted centers to the true centers, so its cost need not equal $L_n+M_n$. We measure this discrepancy by the signed transport residual
\[
R_n^{\mathrm{OT}}
=
\underbrace{W_2^2(\widehat G_n,G^\star)}_{\text{optimal transport cost}}
-
\underbrace{L_n}_{\text{local displacement cost}}
-
\underbrace{M_n}_{\text{mass-reallocation cost}}
\]
The residual measures how closely the sum of the two contributions approximates the actual Wasserstein loss. Since it can be positive or negative, we use its absolute value to quantify the size of this discrepancy. For a direct comparison with the scaled Wasserstein loss $\sqrt n W_2^2(\widehat G_n,G^\star)$, we report the local displacement cost, mass-reallocation cost, and absolute residual as $\sqrt n L_n$, $\sqrt n M_n$, and $\sqrt n|R_n^{\mathrm{OT}}|$, respectively.

\begin{figure}[t]
\centering
\includegraphics[width=1\textwidth]{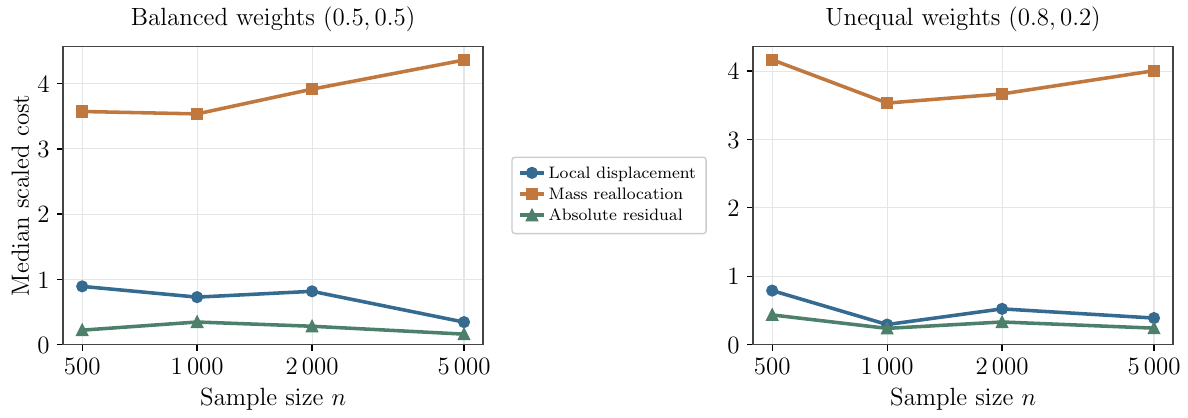}
\caption{Local displacement, mass reallocation, and absolute transport residual for a two-component Gaussian mixture fitted with at most three components. The true centers have separation \(\Delta=2\), with mixing weights \((0.5,0.5)\) on the left and \((0.8,0.2)\) on the right. Curves show the median of each quantity after scaling by \(\sqrt n\).}
\label{fig:simBdecomp}
\end{figure}

\paragraph{Results.} Figure~\ref{fig:simBdecomp} shows how these quantities vary with sample size at separation $\Delta=2$, for both balanced and unequal mixing weights. Table~\ref{tab:simBsep} complements this comparison by reporting the median quantities at $n=5000$ for balanced truths with separations $\Delta=2$ and $\Delta=4$. In both configurations, the median mass-reallocation cost is substantially larger than the median local displacement cost. For the balanced truth at $\Delta = 4$, the median scaled mass-reallocation cost is $4.549$, compared with $0.875$ for local displacement and $0.054$ for the absolute residual. These results show that discrepancies in component masses can make a substantial contribution to Wasserstein error, even after the fitted atoms are moved onto the true centers.

The decomposition also reveals differences that are less apparent from the total error. Although the median normalized Wasserstein errors $T_n$ are similar at the two separations, the median local contribution increases from $0.344$ at $\Delta=2$ to $0.875$ at $\Delta=4$. Meanwhile, the median mass-reallocation contributions remain close, at $4.363$ and $4.549$, respectively. Despite their similar total Wasserstein errors, the two configurations differ in the contributions from local displacement and mass reallocation. By separating these contributions, we can distinguish errors in the fitted locations from errors in the mass assigned to each true component.

\subsubsection{Effect of Component Budget on Wasserstein Error}\label{subsec:simC}

The fourth-root scale presented in Section \ref{subsec:mixing-error-scaling} describes the magnitude of mixing-measure error, but it does not by itself determine the error distribution. Our theory further characterizes this distribution through the local perturbations available to the fitted model, which depend on the permitted number of components. More specifically, increasing the component budget expands the set of attainable perturbations and may change the direction selected by likelihood maximization. To clarify the resulting effect on estimation error, we compare different fitted component budgets while keeping the true mixing measure fixed.

\paragraph{Setup.} We generate observations from the balanced two-point mixing measure $G^\star=\frac12\delta_{-1}+\frac12\delta_1$ To extend the analysis beyond the lower-fluorescence group, we compare the complete fitted mixing measures using Wasserstein distance. The case $k=2$ matches the true number of components, whereas $k=3$ and $k=4$ allow one and two additional components, respectively. These are upper limits on the fitted component count, so the fit is not required to use every available component. For each configuration, we use $100$ independently generated datasets at sample sizes ranging from $500$ to $5000$. We compare the resulting distributions of $T_n=n^{1/4}W_2(\widehat G_n,G^\star)$ through their medians and $90^{\mathrm{th}}$ percentiles. The medians summarize typical estimation error, whereas the upper percentiles describe larger errors and are relevant to confidence-region calibration.

\begin{figure}[t]
\centering
\includegraphics[width=0.75\textwidth]{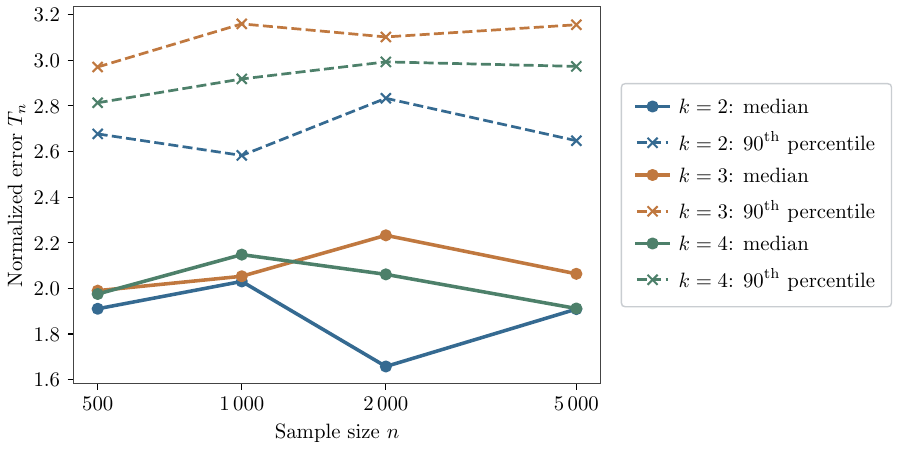}
\caption{Wasserstein estimation error $n^{1/4}W_2(\widehat G_n,G^\star)$ across sample sizes for fitted component budgets $k=2,3,4$, with $G^\star=\tfrac12\delta_{-1}+\tfrac12\delta_1$. Solid and dashed curves show the medians and $90^{\mathrm{th}}$ percentiles, respectively.}
\label{fig:simCbudget}
\end{figure}

\paragraph{Results.} For the three fitted component budgets, Figure~\ref{fig:simCbudget} compares the medians and $90^{\mathrm{th}}$ percentiles of $T_n$ across sample sizes. In particular, at $n=5000$, the $90^{\mathrm{th}}$ percentile increases when the component budget changes from $k=2$ to $k=3$, but decreases when the budget is further increased to $k=4$. These differences indicate that the estimated upper quantiles depend on the fitted component budget, even when the true mixing measure and sample size are held fixed.

\noindent The dependence on the fitted component budget is consistent with its role in our weak-limit theory. The attainable tangent sets satisfy $\mathcal T_2\subset\mathcal T_3\subset\mathcal T_4$, so increasing $k$ gives likelihood maximization additional local perturbations to choose from. Because the selected perturbation maximizes the likelihood objective rather than minimizes transportation cost, this inclusion does not imply that Wasserstein error or its upper quantiles must increase with $k$. This distinction explains why the nonmonotone ordering observed at $n=5000$ is compatible with the nested tangent sets, although the theory does not predict the relative sizes of these quantiles.

\noindent This experiment implies that the fourth-root scale alone is insufficient to determine a confidence-region radius. The radius also depends on an upper quantile of the error distribution, whose theoretical characterization involves the fitted component budget. To account for this dependence, bootstrap refits should use the same component budget as the original fit, so that the resulting critical values reflect the estimation procedure being assessed.
\subsubsection{Component Budget and Multivariate Splitting Geometry}\label{subsec:simE}

In Section~\ref{subsec:simC}, we show that changing the fitted component budget can affect the distribution of Wasserstein error. In multiple dimensions, additional components also allow different spatial arrangements around a true center. For example, fitted centers may spread mainly along a line or extend along a second independent direction. These arrangements need not be distinguishable from the total Wasserstein error alone, which motivates us to examine how the fitted component budget affects splitting geometry and how the observed configurations relate to the rank constraint in \eqref{eq:budget}.

\paragraph{Setup.} We generate observations from $\mathcal N(0,I_2)$, corresponding to the true mixing measure $G^\star=\delta_{(0,0)}$, and fit Gaussian location mixtures with at most $k=2$ or $k=3$ components. By considering a single true component, we exclude mass reallocation between distinct true centers and focus on how the fitted atoms are arranged around the origin. Both budgets permit over-specification, and comparing them isolates the effect of allowing one additional fitted atom while leaving the sampling distribution unchanged. For each budget, we use $100$ independently generated datasets at each sample size $n\in\{500,1000,2000,5000\}$.

\paragraph{Measuring splitting geometry.} To describe the arrangement of the fitted atoms, we compute their rescaled second-moment matrix about the true center, $\widehat C_n
=
\sqrt n\sum_j\widehat\pi_j
\widehat\theta_j\widehat\theta_j^\top.$
For any unit vector $u$, the quantity $u^\top\widehat C_n u
=
\sqrt n\sum_j\widehat\pi_j
\bigl(u^\top\widehat\theta_j\bigr)^2$
measures the weighted squared displacement of the fitted centers along direction $u$. Through these directional second moments, $\widehat C_n$ describes both how far the fitted atoms lie from the origin and how their displacements are distributed across directions. The factor $\sqrt n$ places these squared displacements on the scale used for second-order local perturbations in our theory.

\noindent Let $\lambda_1(\widehat C_n)\geq\lambda_2(\widehat C_n)\geq0$ denote its eigenvalues. The larger eigenvalue measures the second moment along the dominant direction, while the smaller eigenvalue measures the second moment along the perpendicular direction. Whenever $\lambda_1(\widehat C_n)>0$, we summarize their relative contributions by $\frac{\lambda_2(\widehat C_n)}{\lambda_1(\widehat C_n)}$.
A ratio close to zero indicates that the fitted second moment is concentrated mainly along one direction. A larger ratio indicates a more substantial contribution from the second direction, with a ratio of one corresponding to equal contributions along the two principal directions. Because the common scaling factor cancels, the ratio describes the relative directional contributions independently of their overall magnitude.

\noindent This distinction matters because Wasserstein error summarizes the total displacement without specifying how it is distributed across directions. To be precise, for the single-atom truth considered here,
\[
\sqrt n\,W_2^2(\widehat G_n,G^\star)
=
\operatorname{tr}(\widehat C_n)
=
\lambda_1(\widehat C_n)+\lambda_2(\widehat C_n).
\]
Consequently, two fits can have the same scaled Wasserstein error even if one spreads mainly along a line and the other has a substantial contribution from a second direction. The eigenvalue ratio distinguishes these configurations. We report its median and $90^{\mathrm{th}}$ percentile across datasets to describe both typical fitted arrangements and those with a larger relative contribution from the second direction.

\noindent Our theory explains why these arrangements can depend on the component budget. Under \eqref{eq:budget}, a local second-moment matrix of rank $r$ requires at least $r+1$ nearby fitted atoms. With $k=2$, a nonzero limiting local splitting matrix therefore has rank at most one. Increasing the budget to $k=3$ makes rank-two local splitting attainable. We investigate whether this additional splitting direction appears in finite-sample fits by comparing the eigenvalue ratios under $k=2$ and $k=3$.

\begin{figure}[t]
\centering
\includegraphics[width=0.91\textwidth]{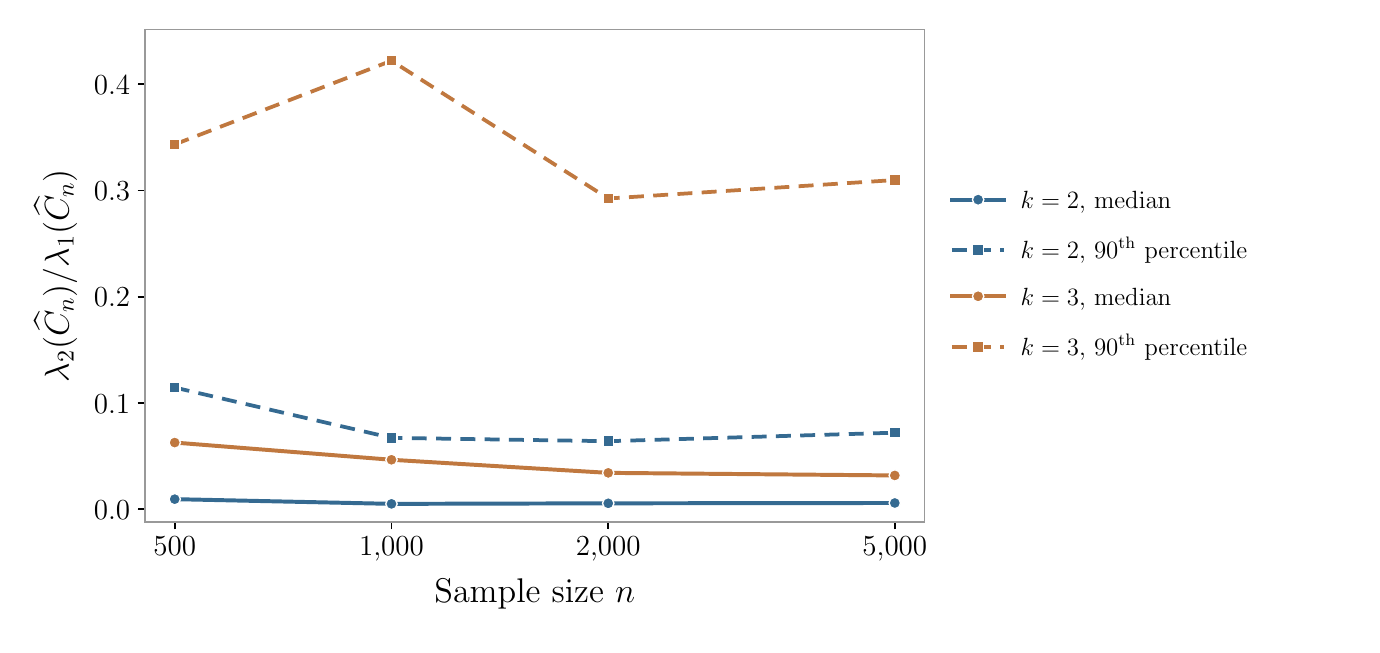}
\caption{Medians and $90^{\mathrm{th}}$ percentiles of $\lambda_2(\widehat C_n)/\lambda_1(\widehat C_n)$ across sample sizes for observations from $\mathcal N(0,I_2)$, fitted with at most $k=2$ or $k=3$ components. The quantities $\lambda_1(\widehat C_n)$ and $\lambda_2(\widehat C_n)$ are the largest and second-largest eigenvalues of the fitted second-moment matrix about the origin. Solid and dashed curves represent the medians and $90^{\mathrm{th}}$ percentiles, respectively.}
\label{fig:simErank}
\end{figure}

\paragraph{Results.} Figure~\ref{fig:simErank} reports the median and $90^{\mathrm{th}}$ percentile of the eigenvalue ratio across sample sizes. At $n=5000$, the two component budgets give the following results.
\begin{enumerate}
    \item[(i)] With $k=2$, the median ratio is close to zero, and the $90^{\mathrm{th}}$ percentile is below $0.1$. The second principal direction contributes little relative to the first, both at the median and at the upper percentile.

    \item[(ii)] With $k=3$, the $90^{\mathrm{th}}$ percentile is around $0.3$, compared with less than $0.1$ under $k=2$. This higher upper percentile suggests that the additional component can accommodate more substantial spreading along the second direction.
\end{enumerate}

\noindent The observed difference is consistent with the additional splitting directions permitted by \eqref{eq:budget}. This constraint concerns limiting local configurations, so it does not require the second eigenvalue of every finite-sample matrix under $k=2$ to be exactly zero. In particular, $\widehat C_n$ is measured about the true center, and finite-sample shifts of the fitted centers can contribute to its second eigenvalue. Similarly, making rank-two splitting attainable under $k=3$ does not require every dataset to produce a fit with substantial spreading in both directions. The experiment illustrates how the component budget restricts the geometry available to likelihood maximization, while the eigenvalue ratios reveal which directional contributions appear in the resulting fits. The eigenvalue ratios complement Wasserstein error by describing how the fitted atoms spread across directions, beyond the total magnitude of their displacements.
\subsection{Flow-Cytometry Analysis with Wasserstein Uncertainty}\label{sec:flowcytometry}

Our theory provides uncertainty regions for an entire mixing measure, from which we derive bounds for population-level summaries without matching individual component labels across fits. We apply this approach to flow-cytometry measurements, where a broad fluorescence group may be represented by different numbers of mixture components. First, we examine how the fitted representation changes with the component budget. Furthermore, we construct bootstrap Wasserstein regions within each specified Gaussian model and use them to bound summaries of the mixing mass near a selected fluorescence profile.

\subsubsection{Data and Gaussian Mixture Model}\label{subsec:flowdata}

\paragraph{Data and preprocessing.} We use the public \emph{Rituximab} dataset from the Bioconductor \texttt{flowClust} package \citep{Gasparetto2004,LoHahne2009}, containing $n=1{,}545$ cell events and eight recorded channels. We analyze \texttt{FL1.H}, which measures anti-BrdU FITC fluorescence related to DNA synthesis, and \texttt{FL3.H}, which measures 7-AAD fluorescence related to DNA content. These two channels describe complementary aspects of the fluorescence profile of each event. We retain all events and transform their recorded intensities by $Z_i=
\left(
\operatorname{asinh}(F_{i1}/150),
\operatorname{asinh}(F_{i3}/150)
\right)^\top$
, where $F_{i1}$ and $F_{i3}$ denote the two channel intensities. The transformation accommodates zero measurements and compresses large values, and all subsequent fitting and transport calculations use these transformed coordinates.

\paragraph{Model and fitting procedure.} We fit Gaussian location mixtures of the form
\begin{equation}\label{eq:flowmodel}
p_G(z)=\sum_{j=1}^{k}\pi_j\varphi_2(z;\mu_j,\sigma^2I_2),
\qquad
G=\sum_{j=1}^{k}\pi_j\delta_{\mu_j},
\qquad
\mu_j\in[-1,5]^2,
\end{equation}
using fixed $\sigma=0.30$ and component budgets $k\in\{2,3,4\}$. Each atom of $G$ represents a latent fluorescence-profile center, with its weight giving the mixing mass assigned to that profile. Fixing the covariance leaves the centers and weights as the quantities to be estimated and gives the Gaussian location family considered in Corollary~\ref{cor:gaussian}. We treat this specification as a working model: its fitted centers describe the fluorescence distribution under the chosen kernel, and their biological interpretation requires additional validation.

\noindent For each budget, we fit the full dataset by multistart EM with $35$ initializations and retain the fit with the highest likelihood. We also assess density prediction using a common fixed 80--20 training/test split. For this comparison, each model is fitted on the training observations and evaluated on the held-out observations. We use these held-out likelihoods to assess predictive fit and the full-data estimates to obtain the population summaries and bootstrap regions below.

\subsubsection{Component Splitting and Mixing-Mass Stability}\label{subsec:flowmodel}

We report the fitted weights and likelihood summaries for the three component budgets in Table~\ref{tab:flowfits}. For component budgets $k=2,3,4$, the held-out log likelihood per event is $-1.382$, $-1.263$, and $-1.187$, respectively. The increase across these budgets indicates better density prediction on the common test set. To understand what changes in the latent representation accompany this improvement, we examine the fitted centers and the mass assigned to the main fluorescence groups.

\begin{table}[t]
\centering
\caption{Rituximab mixture fits with fixed $\sigma=0.30$. Log likelihood, BIC, and component weights use the full-data fits; held-out log likelihood uses a common $80$--$20$ training/test split. Within each fit, components are ordered by the first coordinate of their centers. Component indices do not imply matching across budgets, and a dash denotes an absent component.}
\vspace{12pt}
\label{tab:flowfits}
\small
\setlength{\tabcolsep}{5pt}
\begin{tabular}{cccccccc}
\toprule
& \multicolumn{2}{c}{\emph{Full-data fit}}
& & \multicolumn{4}{c}{\emph{Component weights}} \\
\cmidrule[0.6pt](lr){2-3}
\cmidrule[0.6pt](lr){5-8}
Budget $k$
& \shortstack{Log\\likelihood}
& BIC
& \shortstack{Held-out log likelihood\\(per event)}
& $\widehat\pi_1$
& $\widehat\pi_2$
& $\widehat\pi_3$
& $\widehat\pi_4$ \\
\midrule
$2$ & $-1916.4$ & $3869.5$ & $-1.382$
& $0.742$ & $0.258$ & \textemdash & \textemdash \\
$3$ & $-1757.5$ & $3573.8$ & $-1.263$
& $0.726$ & $0.201$ & $0.072$ & \textemdash \\
$4$ & $-1655.4$ & $3391.6$ & $-1.187$
& $0.098$ & $0.643$ & $0.186$ & $0.074$ \\
\bottomrule
\end{tabular}
\end{table}

\begin{figure}[t]
\centering
\includegraphics[width=0.98\linewidth]{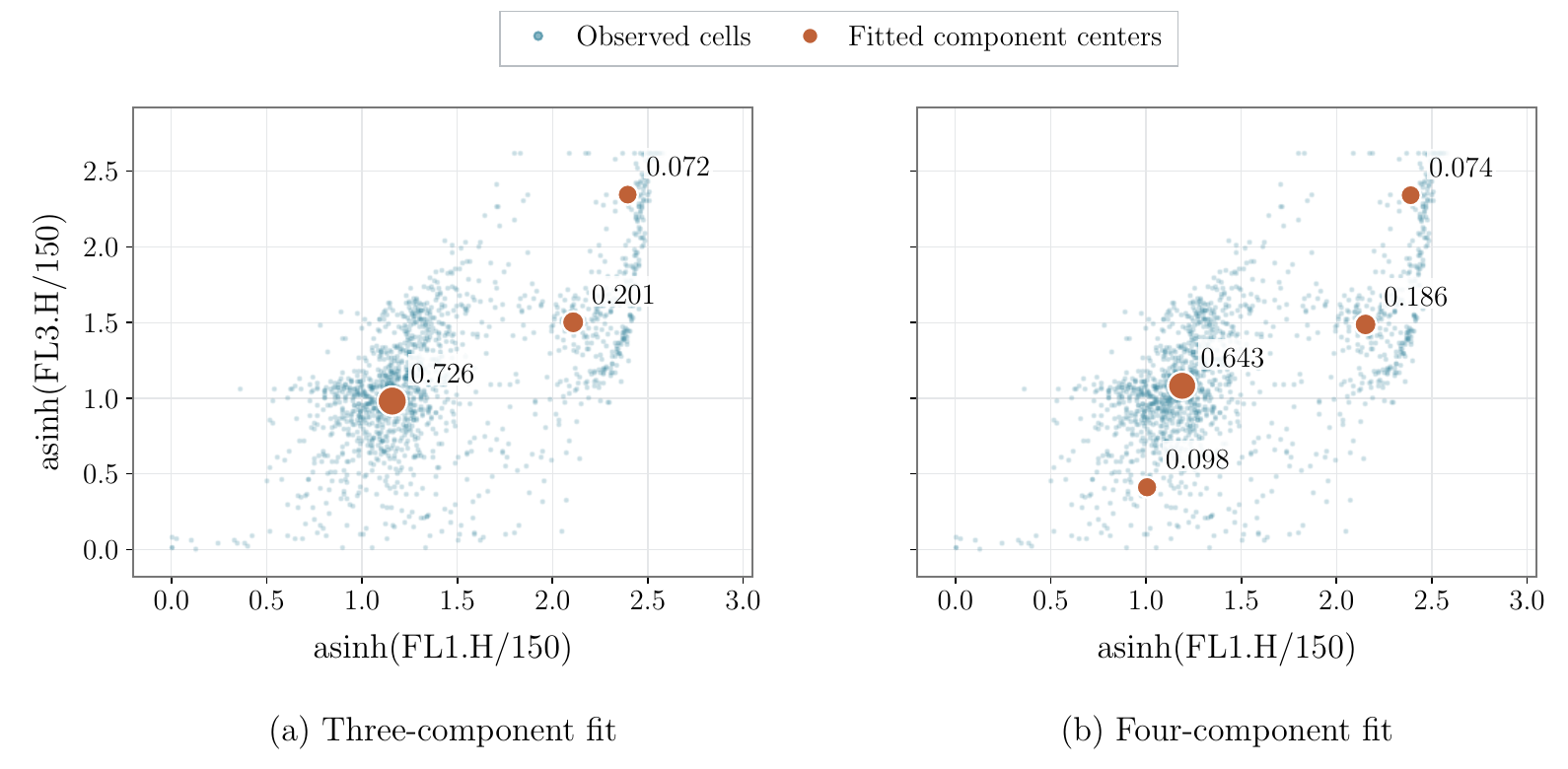}
\caption{Transformed Rituximab observations and fitted mixing atoms, with their weights, for (a) three and (b) four components. The principal lower-fluorescence region receives about $73\%$ and $74\%$ of the total mixing mass, respectively, while being represented by different numbers of atoms.}
\label{fig:flowpopulations}
\end{figure}

\noindent Figure~\ref{fig:flowpopulations} shows the three- and four-component fits over the transformed observations. With three components, the principal lower-fluorescence group is represented by an atom near $(1.2,1.0)$ with weight $0.726$. As for the four-component fit, two nearby atoms describe this region, carrying weights 0.098 and 0.643. Their combined weight is about $0.74$, close to that of the single atom in the three-component fit. Although the representation changes from one atom to two, the total mass assigned to the broad fluorescence region remains similar. This change makes a componentwise comparison ambiguous. The weight $0.726$ in the three-component fit has no single counterpart in the four-component fit because two atoms contribute to the same broad region. Combining their contributions gives a more direct description of the mass in that region. To compare this region across different component representations, we use summaries based on the locations and weights of all contributing atoms.

\noindent We also compare the complete fitted mixing measures using Wasserstein distance. Exact discrete transport gives $W_2(\widehat G_2,\widehat G_3)=0.245 \text{ and } W_2(\widehat G_3,\widehat G_4)=0.237.$ These distances account jointly for changes in locations and weights without requiring a one-to-one correspondence between atoms. They describe sensitivity to the component budget. To assess sampling uncertainty for a given budget, we next keep that budget fixed and use the bootstrap to construct a region around its fitted measure.

\subsubsection{Bootstrap Wasserstein Regions}\label{subsec:flowbootstrap}

\paragraph{Bootstrap construction.} For each $k\in\{3,4\}$, we generate $B=250$ parametric bootstrap datasets from the corresponding fitted observation distribution. Each bootstrap dataset is refitted using eight EM initializations, with the component budget, covariance, and center constraints held fixed. Through these refits, we quantify sampling uncertainty separately for each component budget.

\noindent We use two bootstrap sample sizes, $m\in\{118,355\}=\{\lfloor n^{0.65}\rfloor,\lfloor n^{0.80}\rfloor\}$, and record $T_{n,m}^{*(b)}=m^{1/4}W_2\bigl(\widehat G_m^{*(b)},\widehat G_n\bigr)$ for $b=1,\ldots,B$.
We use these normalized bootstrap errors to estimate a critical value on the fourth-root scale. Following Corollary~\ref{cor:bootquantile}, we convert this critical value to a radius for the original sample size by
\begin{equation}\label{eq:flowradius}
\widehat r_{.95}(m)
=
n^{-1/4}\operatorname{Quantile}_{.95}
\left\{
T_{n,m}^{*(1)},\ldots,T_{n,m}^{*(B)}
\right\}.
\end{equation}
The resulting ball consists of mixing measures within distance $\widehat r_{.95}(m)$ of $\widehat G_n$ in the specified model class. Its nominal $95\%$ coverage interpretation follows from the assumptions of the corollary, with the likelihood maximizers approximated numerically by multistart EM.

\paragraph{Results.} The estimated critical values and radii are given in Table~\ref{tab:flowbootstrap}. Figure~6(a) shows that all four radii are close to $0.2$. At $m=355$, the radii under $k=3$ and $k=4$ are nearly equal, although the corresponding uncertainty regions are centered at different fitted mixing measures.

\begin{table}[t]
\centering
\caption{Bootstrap critical values and radii of nominal $95\%$ Wasserstein confidence regions for the Rituximab mixing measure ($n=1545$, fixed $\sigma=0.30$). Each $(k,m)$ configuration uses $B=250$ parametric bootstrap samples of size $m$, with at most $k$ components in both the original fit and bootstrap refits. Here $\widehat q_{.95}^{*}$ is the empirical $95$th percentile of $m^{1/4}W_2(\widehat G_m^{*},\widehat G_n)$, and \mbox{$\widehat r_{.95}(m)=n^{-1/4}\widehat q_{.95}^{*}$} is the radius of the ball centered at $\widehat G_n$.}
\vspace{12pt}
\label{tab:flowbootstrap}
\small
\setlength{\tabcolsep}{5pt}
\renewcommand{\arraystretch}{1.1}
\begin{tabular}{cccccc}
\toprule
& & \multicolumn{2}{c}{\textit{Component budget $k=3$}}
  & \multicolumn{2}{c}{\textit{Component budget $k=4$}} \\
\cmidrule[0.6pt](lr){3-4}
\cmidrule[0.6pt](lr){5-6}
\shortstack{Bootstrap\\sample size $m$}
& \shortstack{Sample-size\\ratio $m/n$}
& \shortstack{Critical value\\$\widehat q_{.95}^{*}$}
& \shortstack{Wasserstein radius\\$\widehat r_{.95}(m)$}
& \shortstack{Critical value\\$\widehat q_{.95}^{*}$}
& \shortstack{Wasserstein radius\\$\widehat r_{.95}(m)$} \\
\midrule
$118$ & $0.076$ & $1.183$ & $0.189$ & $1.259$ & $0.201$ \\
$355$ & $0.230$ & $1.255$ & $0.200$ & $1.238$ & $0.198$ \\
\bottomrule
\end{tabular}
\end{table}

\noindent The between-budget distance and the bootstrap radii address separate questions. The distance $W_2(\widehat G_3,\widehat G_4)=0.237$ measures the change produced by allowing an additional component. Each bootstrap radius instead measures uncertainty within one fixed model specification. We retain both comparisons because a population summary can be sensitive to the permitted representation as well as to sampling variability.

\subsubsection{Bounds for Population-Level Summaries}\label{subsec:flowinterpretation}

We now use the Wasserstein regions to bound the mixing mass near a selected fluorescence profile. A direct sum of nearby component weights requires deciding which fitted atoms belong to the region. Such a sum can change abruptly when an atom crosses the chosen boundary. To obtain a summary that varies gradually as fitted centers move, we use the smooth gate
\[
\psi_{c,h}(\mu)
=
\left(1-\frac{\|\mu-c\|}{h}\right)_+,
\qquad
M_{c,h}(G)
=
\int\psi_{c,h}(\mu)\,dG(\mu),
\]
where $c\in\mathbb R^2$ specifies the profile of interest and $h>0$ specifies the neighborhood width. The gate assigns value one at $c$, decreases linearly with distance, and assigns zero beyond distance $h$.

\noindent For a fitted discrete measure, this gives $M_{c,h}(\widehat G_n)=\sum_j\widehat\pi_j\psi_{c,h}(\widehat\mu_j)$. Every fitted atom contributes according to both its weight and its distance from $c$. Whether the region is represented by one atom or several, this definition combines their contributions without depending on component labels. The quantity describes the latent fluorescence-profile centers under the mixture model; it differs from the fraction of observed events falling inside a hard gate.
\paragraph{Uncertainty bounds.} To quantify uncertainty in this summary, we bound its deviation from the fitted value in terms of the Wasserstein distance. Since $\psi_{c,h}$ has Lipschitz constant $1/h$, the Kantorovich--Rubinstein bound and $W_1\leq W_2$ imply
\begin{equation}\label{eq:flowgatebound}
\bigl|M_{c,h}(G)-M_{c,h}(\widehat G_n)\bigr|
\leq
\frac{1}{h}W_1(G,\widehat G_n)
\leq
\frac{1}{h}W_2(G,\widehat G_n).
\end{equation}
Hence, every mixing measure in a Wasserstein ball of radius $r$ has smooth-gate mass within $r/h$ of the fitted value. Because the mass lies in $[0,1]$, the induced interval is
\[
\left[
\max\{0,M_{c,h}(\widehat G_n)-r/h\},
\;
\min\{1,M_{c,h}(\widehat G_n)+r/h\}
\right].
\]
This inequality holds for every $c$ and $h>0$. On the event that the target mixing measure belongs to the confidence ball, all these bounds hold simultaneously, so the same region can be used to study several fluorescence profiles without constructing a separate bootstrap region for each summary.

\noindent For a fixed Wasserstein radius $r$, the bandwidth $h$ determines both the size of the neighborhood around $c$ and the uncertainty bound $r/h$. In particular, a smaller $h$ focuses the summary on a narrower neighborhood, but increases the uncertainty allowance $r/h$ for a fixed Wasserstein radius. By contrast, a larger $h$ averages contributions over a broader neighborhood and gives a smaller allowance. This dependence follows directly from how quickly the gate changes with the location of an atom.

\begin{table}[t]
\centering
\caption{Estimates and conservative confidence intervals for the lower-fluorescence smooth-gate mass in the Rituximab data. The mass estimate is $M_{c,h}(\widehat G_n)=\sum_j\widehat\pi_j(1-\|\widehat\mu_j-c\|/h)_+$, with $c=(1.16,1.05)$ and $h=1.20$. Each interval is the estimate plus or minus the error bound $\widehat r_{.95}(m)/h$, truncated to $[0,1]$. The radii are from the nominal $95\%$ Wasserstein regions in Table~\ref{tab:flowbootstrap}, using $m=355$, $B=250$, and fixed $\sigma=0.30$. Endpoints use unrounded values.}
\vspace{12pt}
\label{tab:flowgate}
\small
\setlength{\tabcolsep}{8pt}
\renewcommand{\arraystretch}{1.1}
\begin{tabular}{cccc}
\toprule
\shortstack{Component budget\\$k$}
& \shortstack{Estimated smooth-gate mass\\$M_{c,h}(\widehat G_n)$}
& \shortstack{Error bound\\$\widehat r_{.95}(m)/h$}
& \shortstack{Confidence interval\\for smooth-gate mass} \\
\midrule
$3$ & $0.710$ & $0.167$ & $[0.543,\,0.876]$ \\
$4$ & $0.681$ & $0.165$ & $[0.517,\,0.846]$ \\
\bottomrule
\end{tabular}
\end{table}

\begin{figure}[t]
\centering
\includegraphics[width=0.98\linewidth]{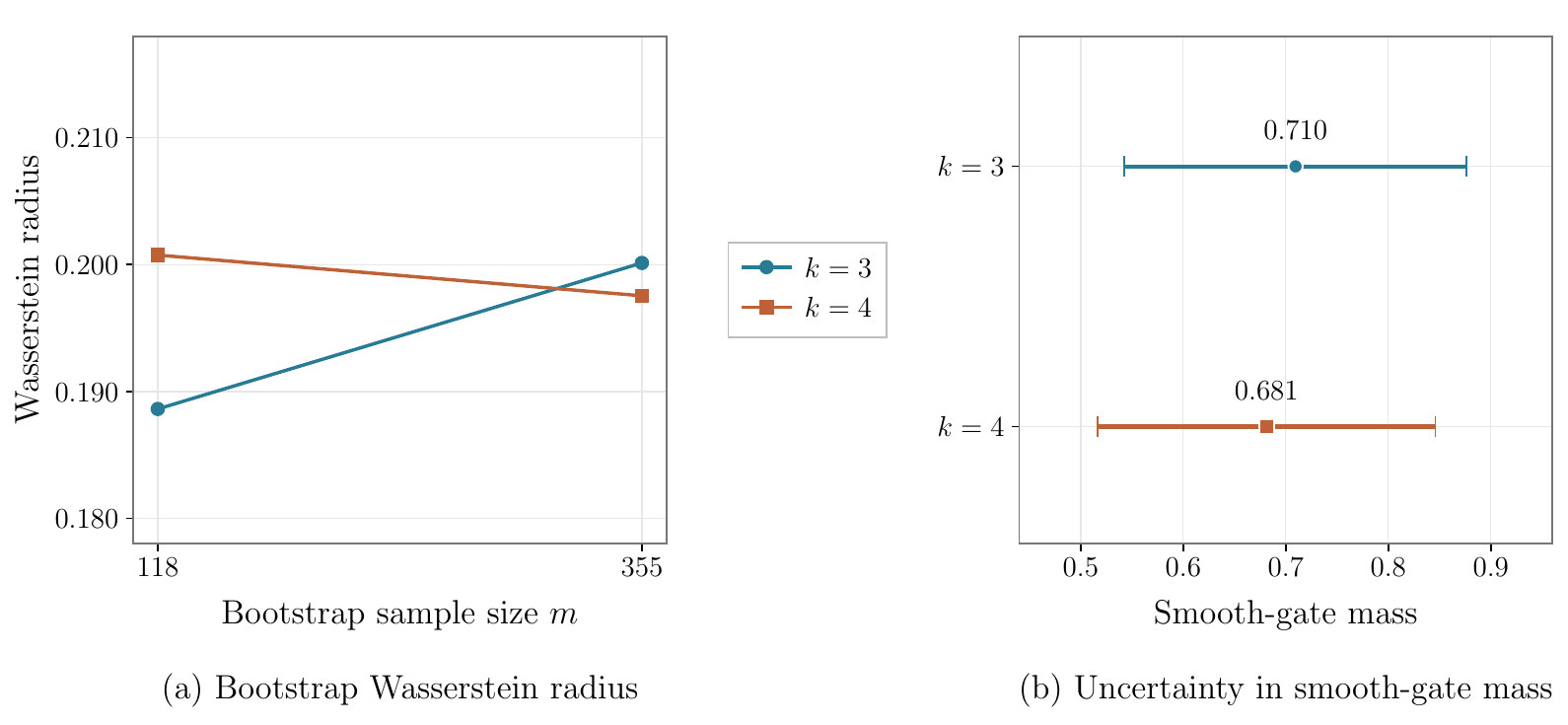}
\caption{Bootstrap uncertainty estimates for the Rituximab mixing measures with fitted component budgets $k=3,4$. (a) Radii of nominal $95\%$ Wasserstein regions for bootstrap sample sizes $m=118,355$. (b) Smooth-gate mass estimates and intervals using $m=355$, $c=(1.16,1.05)$, and $h=1.20$. Points denote the estimates, and horizontal bars show the intervals obtained by adding and subtracting the corresponding Wasserstein radius divided by $h$.}
\label{fig:flowinference}
\end{figure}

\paragraph{Results.} For the lower-fluorescence region, we take $c=(1.16,1.05)$ and $h=1.20$, using the bootstrap radii obtained with $m=355$. The estimated smooth-gate mass is $0.710$ under $k=3$ and $0.681$ under $k=4$. As shown in Table~\ref{tab:flowgate} and Figure~\ref{fig:flowinference}(b), the corresponding confidence intervals have similar widths and overlap substantially. To interpret these estimates, recall that the smooth gate weights each atom's contribution by its distance from $c$, whereas the combined weights reported in Section~\ref{subsec:flowmodel} simply sum the weights of the atoms assigned to the group. Although the lower-fluorescence group is represented by different numbers of atoms, the same summary remains well defined under both budgets, and each fitted model supplies an associated uncertainty interval through its Wasserstein region. This suggests a practical use of inference for the complete mixing measure, as a bound on its estimation error yields uncertainty bounds for population summaries that remain meaningful as the atomic representation changes.

\section{Discussion}\label{sec:discussion}

We studied the weak convergence of the Wasserstein estimation
error for global maximum likelihood estimation in finite mixtures
with a fixed upper bound on the number of components. Our analysis
follows the density-to-mixing-measure approach used in
convergence-rate theory, but develops a distributional counterpart
to the inverse inequalities underlying that approach. We first
derive a uniform quadratic expansion of the log-likelihood
difference and characterize the attainable relative-density
directions under the finite component budget. The density weak
limit is then determined by maximizing a Gaussian quadratic
functional over a generally nonconvex set. A sequential expansion
of the squared Wasserstein distance connects this selected
direction to the leading transportation cost, yielding joint
weak convergence of the relative-density estimator and its
Wasserstein error. We also establish an $m$-out-of-$n$ parametric
bootstrap that approximates the resulting distribution and
provides Wasserstein confidence regions at regular target
quantiles.

Our approach opens several directions for further investigation.
On the theoretical side, a natural extension is to heteroscedastic
location-scale mixtures, where scale parameters vary across
components. As discussed in Remark~\ref{remark:non-scale_location},
the derivative identities of Gaussian location-scale kernels
lead to higher-order cancellations and require a different
characterization of attainable directions. Related progress for
infinite mixtures with an unknown shared scale is provided by
\citep{BarilettoLeRinaldoHo2026}, who establish inverse inequalities
connecting density discrepancies to errors in the location
mixing measure and scale parameter. Extending the present
weak-limit theory to component-specific scales would require
identifying the appropriate higher-order density directions
and the corresponding transportation expansion. The interaction
between derivative cancellations and the number of redundant
components would be central to this analysis.

On the methodological and experimental side, a natural direction
is to connect our likelihood geometry with mixture-order selection and penalized estimation \citep{Keribin2000,ManoleKhalili2021}, accounting for the singular structure that also motivates information criteria beyond regular likelihood theory \citep{DrtonPlummer2017}. Increasing the fitted
order enlarges the model and can improve the maximized
log-likelihood, even when the additional components are redundant. A joint weak-limit analysis of these improvements across candidate orders could clarify how penalties affect the probabilities of selecting different component counts. Numerical studies could then examine penalty calibration and bootstrap approximation of selection probabilities under varying component separation and optimization accuracy.

\clearpage
\appendix
\startcontents[appendices]

\begin{center}
\textbf{\Large Appendices for\\[0.5em]
``Weak Limits of Parameter Estimation\\
in Strongly Identifiable Mixture Models''}
\end{center}

\printcontents[appendices]{}{1}{\setcounter{tocdepth}{2}}

% Appendix content
% Appendix content

\section{Background and Preliminary Results}
\label{app:preliminaries}
Before presenting the detailed proofs, we recall several facts
used repeatedly below. First, Cauchy--Schwarz inequality gives
\[
    \int|p_G-p_*|\,d\mu
    =
    P_*\left|\frac{p_G-p_*}{p_*}\right|
    \leq
    \left\|\frac{p_G-p_*}{p_*}\right\|_*
    =\rho(G),
\]
so convergence in relative $L^2(P_*)$ norm implies
$L^1(\mu)$ convergence of the mixture densities. Throughout the proofs, unless otherwise stated, $\nabla f_\theta$,
$\nabla^2f_\theta$, and $\nabla^r f_\theta$ denote derivatives
with respect to the parameter $\theta$.
Second, Assumption~\ref{ass:smooth} implies, uniformly in $\theta\in\Theta$,
\[
\left|\frac{f_\theta}{p_*}\right|\le F,
\qquad
\left\|\frac{\nabla^r f_\theta}{p_*}\right\|\le F,
\quad r=1,2,3,
\]
with $F\in L^{2+\eta}(P_*)$. In particular, all displayed relative derivatives belong to $L^2(P_*)$. The dominated convergence theorem also gives $L^2(P_*)$-continuity of $\theta\mapsto f_\theta/p_*$ and of its derivatives through order two. Whenever we pass a kernel, gradient, Hessian, or exact secant to a parameter limit below, this is the domination being used.
We also use the subsequence characterization of convergence
in probability: $Z_n\to0$ in probability if and only if every
subsequence contains a further subsequence along which $Z_n\to0$
almost surely. For deterministic sequences, compactness and
the bounds established below allow us to extract subsequences
along which the relevant coefficients and atom locations
converge.
% We collect a few facts used throughout the appendices. Most are standard consequences of Wasserstein topology, or elementary norm inequalities; see \citet{vdV1998,Villani2009,Santambrogio2015} for details.  They are recorded here to keep the later proofs focused on the mixture-specific arguments.

% \paragraph{Relative-density norm controls total variation of densities.}
% For every $G\in\cO_k(\Theta)$,
% \begin{equation}
% \label{eqn:dung_appendix_A_relative_density_inequality}
%     \int |p_G-p_*|\,d\mu
% =P_*\left|\frac{p_G-p_*}{p_*}\right|
% \le \left\|\frac{p_G-p_*}{p_*}\right\|_{L^2(P_*)}.
% \end{equation}
% Thus convergence in the relative $L^2(P_*)$ norm implies $L^1(\mu)$ convergence of the mixture densities. This observation is used repeatedly to combine compactness with identifiability.

% \paragraph{Consequences of the derivative envelope.}

% \paragraph{Subsequence principle.}
% Many assertions below are proved by the following standard device. To show that a random quantity $Z_n$ converges in probability to zero, it is enough to show that every subsequence contains a further subsequence along which $Z_n\to0$ almost surely. Likewise, to identify a deterministic sequential limit, we may pass to subsequences until every one of the finitely many atom labels, coefficients, and support locations converges. Because $k$ is fixed, no diagonal extraction over an increasing number of labels is ever required.

\subsection{Consistency and local geometry}

We establish consistency, an exact density decomposition,
and local comparisons between relative-density $\rho(G)$ and Wasserstein
discrepancies $W_2(G,G_*)$. These results provide the geometric ingredients
for the convergence-rate and weak-limit arguments developed
below, after establishing the required empirical-process
properties.

\begin{lemma}[Consistency]\label{lem:consistency}
Under Assumptions~\ref{ass:smooth}-\ref{ass:strongid}, we have 
\[ 
W_1(\widehat G_n,G_*)\pto0,
\qquad
\rho(\widehat G_n):=\left\|\frac{p_{\widehat G_n}-p_*}{p_*}\right\|_*\pto0.
\]
\end{lemma}

\begin{proof}
This is a standard compact $M$-estimation argument once the parameter is viewed as a mixing measure, see \citet{KieferWolfowitz1956,Leroux1992,vdV1998} for comparison. For completeness, we provide a proof below. 

\medskip\noindent\emph{Consistency of $\widehat G_n$.} We prove that $W_1(\widehat G_n,G_*)\pto0$. As $\Theta$ is a compact set, we can verify that $\mathcal{O}_k(\Theta)$ is a weakly compact set. This weak compactness can similarly be understood that $\mathcal{O}_k(\Theta)$ is compact under any Wasserstein metric, as these distances can metricize the weak convergence in compact support. 

Consider the functional $M:\mathcal{O}_k(\Theta) \to \mathbb{R}$, $G\mapsto P_*\log p_G$, then $M$ is continuous. Indeed, under Assumption~\ref{ass:consistency} about integrable envelope for the log-likelihood class, for $G_m \Rightarrow G$, using dominated convergence theorem, we have $M(G_m) \to M(G)$. Furthermore, for any $G \in \mathcal{O}_k(\Theta)$, 
\begin{equation*}
    M(G)-M(G_*)
=P_*\log\frac{p_G}{p_*}
=-D_{\mathrm{KL}}(P_*\|P_G)\le0,
\end{equation*}
where the equality holds only if $p_G = p_*$, $\mu$-a.e., which implies that $G = G_*$ by Assumption~\ref{ass:consistency}. Hence $G_*$ is the unique maximizer of $M$.  

Fix $\eps>0$ and let $K_\eps=\{G\in\cO_k(\Theta):W_1(G,G_*)\ge\eps\}$, which is a compact set. Since $M$ is continuous and $G_*$ is its unique maximizer, there exists a positive gap $\delta_\eps$ between $M(G_*)$ and $M(G)$ for $G \in K_{\varepsilon}$, in other words, 
\begin{equation*}
\delta_\eps
:=M(G_*)-\sup_{G\in K_\eps}M(G)>0.
\end{equation*}
Assumption~\ref{ass:consistency} that the class $\{\log p_{G}:G\in \mathcal{O}_k(\Theta)\}$ is $P_*$-Glivenko-Cantelli with integrable envelop implies that
\begin{equation*}
\sup_{G\in\cO_k(\Theta)}|P_n\log p_G-P_*\log p_G|\pto0.
\end{equation*}
On the event that the preceding supremum is smaller than $\delta_\eps/3$, every $G\in K_\eps$ satisfies
\[
P_n\log p_G
\le M(G)+\frac{\delta_\eps}{3}
\le M(G_*)-\frac{2\delta_\eps}{3}
\le P_n\log p_* -\frac{\delta_\eps}{3}.
\]
Such a $G$ cannot maximize the empirical likelihood. Therefore $\PP\{W_1(\widehat G_n,G_*)\ge\eps\}\to0$.

\medskip\noindent\emph{Continuity of the $\rho(\widehat G_n)$.}\ We prove that $\rho(\widehat G_n) \pto 0$. Consider a sequence $\{G_m\}_{m\geq 0}\subset \mathcal{O}_k(\Theta)$ such that $G_m\Rightarrow G_*$. For $P_*$-a.e. $x$, we have 
${p_{G_m}(x)}/{p_*(x)}\longrightarrow1$.
Moreover, Assumption~\ref{ass:smooth} with $r=0$ gives
\[
0\le \frac{p_{G_m}(x)}{p_*(x)}
=\int \frac{f_\theta(x)}{p_*(x)}\,dG_m(\theta)
\le F(x),
\]
which implies
\[
\left|\frac{p_{G_m}(x)-p_*(x)}{p_*(x)}\right|^2
\le (F(x)+1)^2.
\]
Because $F\in L^{2+\eta}(P_*)$, the dominating function $(F+1)^2$ is integrable. Applying Dominated Convergence Theorem, we have 
\begin{equation}
\label{eqn:dung_weak_convergence_implies_a_s_convergence_of_density}
    \left\|\frac{p_{G_m}-p_*}{p_*}\right\|_*\to0.
\end{equation}
We now pass from $W_2$ consistency of $\widehat G_n$ to the limit $\rho(\widehat G_n)\to 0$. Suppose that $\rho(\widehat G_n)$ does not converge in probability to zero, then there would exist a subsequence and $\eps>0$ on which $\PP\{\rho(\widehat G_n)>\eps\}$ stays bounded away from zero. A further subsequence, denote as $\{\widehat G_{n_k}\}_{k\geq 1}$, can be chosen along which $\widehat G_{n_k}\to G_*$ almost surely; then the result in equation \eqref{eqn:dung_weak_convergence_implies_a_s_convergence_of_density} above would then imply $\rho(\widehat G_{n_k})\to0$ almost surely, a contradiction. Therefore
\[
\rho(\widehat G_n)=\left\|\frac{p_{\widehat G_n}-p_*}{p_*}\right\|_*\pto0.
\]
\end{proof}

% We next state Lemma \ref{lem:secant} that expresses the density discrepancy \(p_G-p_*\) through second-order Taylor expansions around the true center nearest to each fitted atom. Such expansions are standard in the analysis of overfitted mixtures; see, for example, \citet{Chen1995,HoNguyen2016,HeinrichKahn2018}.
The following lemma expresses the density discrepancy $p_G-p_*$
through second-order Taylor expansions around the true atom
nearest to each fitted atom, called \textit{exact second-order secant decomposition}. Such decompositions are standard
in the analysis of overfitted mixtures
\citep{Chen1995,HoNguyen2016,HeinrichKahn2018}. We state this
representation as an intermediate step toward establishing
the convergence rate of $W_2^2(\widehat G_n,G_*)$ and, more
importantly, to explain the local moment structure underlying
the sequential first-order transportation expansion in
Section~\ref{subsec:transportfunctional}.
\begin{lemma}[Exact second-order secant decomposition]\label{lem:secant}
Fix a measurable nearest-true-atom map
$\iota(\theta)\in\argmin_{1\le i\le p}\norm{\theta-\theta_i^*}$,
with deterministic tie breaking. For $G=\sum_{j=1}^{m}\pi_j\delta_{\theta_j}\in\cO_k(\Theta)$, define
$i_j=\iota(\theta_j)$, $\Delta_j=\theta_j-\theta_{i_j}^*$, and $r_j=\norm{\Delta_j}$ be the distance from $\theta_j$ to its nearest true center.
For $r>0$, and unit vector $u\in\mathbb S^{d-1}$, set
\[
Q_{i,r,u}=
\frac{f_{\theta_i^*+ru}-f_{\theta_i^*}-r u^\top\nabla f_{\theta_i^*}}{r^2},
\]
with continuous extension in $L^2(P_*)$ when $r = 0$ as  $Q_{i,0,u}=\frac12u^\top\nabla^2f_{\theta_i^*}u$.
Let
\[
a_i(G)=\left(\sum_{j:i_j=i}\pi_j\right)-\pi_i^*,
\qquad
b_i(G)=\sum_{j:i_j=i}\pi_j\Delta_j,
\qquad
\lambda_j(G)=\pi_jr_j^2.
\]
The following representation is exact: 
\begin{equation}\label{eq:exactsecant}
p_G-p_*
=\sum_{i=1}^{p}a_i(G)f_{\theta_i^*}
+\sum_{i=1}^{p}b_i(G)^\top\nabla f_{\theta_i^*}
+\sum_{j:r_j>0}\lambda_j(G)Q_{i_j,r_j,u_j},
\end{equation}
where $u_j=\Delta_j/r_j$. Moreover, the zero-order parameters $a_i(G)$ have vanishing sum $\sum_{i=1}^{p}a_i(G)=0$.
\end{lemma}
\begin{proof}
\medskip\noindent\emph{Taylor expansion for each true atom.}\ We develop the Taylor expansion for density function $f_{\theta}$ at each true atom $\theta_i^*$. Fix $i,r>0$, and $u\in\mathbb S^{d-1}$ such that $\theta_i^*+ru\in\Theta$. Define
$g(s)=f_{\theta_i^*+sru}$ for $0\le s\le1$. By Assumption~\ref{ass:smooth}, $g$ is twice continuously differentiable with respect to $s$ and
\[
g'(s)=r u^\top\nabla f_{\theta_i^*+sru},
\qquad
g''(s)=r^2 u^\top\nabla^2 f_{\theta_i^*+sru}u.
\]
Taylor's formula with integral remainder yields
\[
g(1)=g(0)+g'(0)+\int_0^1(1-s)g''(s)\,ds.
\]
Consequently, we have 
\begin{equation}\label{eq:QintegralDetailed}
Q_{i,r,u}
=\int_0^1(1-s)u^\top\nabla^2 f_{\theta_i^*+sru}u\,ds.
\end{equation}
By continuity of $\nabla^2 f_{\theta}$ in $\theta$, we can naturally extend $Q_{i,r,u}(x)$ to the case when $r = 0$ as 
\[
Q_{i,r,u}(x)\to \frac12u^\top\nabla^2f_{\theta_i^*}(x)u
:=Q_{i,0,u}(x)
\]
for $P_*$-a.e. $x$. Moreover, from \eqref{eq:QintegralDetailed} and Assumption~\ref{ass:smooth},
$\left|{Q_{i,r,u}(x)}/{p_*(x)}\right|
\le \frac12 F(x)$.
Hence the extension is also continuous in $L^2(P_*)$ by dominated convergence.

\medskip\noindent\emph{Aggregation over fitted atoms.}\ Now we apply this expansion to the fitted atom $f_{\theta_j}$ at the true atom $\theta^*_{i_{j}}$ before summing up. For each fitted atom $\theta_j$ with $r_j>0$, the above expansion gives
\[
f_{\theta_j}
=f_{\theta_{i_j}^*}
+\Delta_j^\top\nabla f_{\theta_{i_j}^*}
+r_j^2Q_{i_j,r_j,u_j}.
\]
If $r_j=0$, the last two terms are zero and the same formula is valid after interpreting $u_j$ arbitrarily. Multiplying this identity by $\pi_j$ and summing over $j$, we have 
\begin{equation*}
    p_G
=\sum_{i=1}^{p}\left(\sum_{j:i_j=i}\pi_j\right)f_{\theta_i^*}
+\sum_{i=1}^{p}\left(\sum_{j:\iota(j)=i}\pi_j\Delta_j\right)\nabla f_{\theta_i^*}
+\sum_{j:r_j>0}\pi_jr_j^2Q_{i_j,r_j,u_j},
\end{equation*}
Subtracting the equation above to $p_*=\sum_i\pi_i^*f_{\theta_i^*}$ gives \eqref{eq:exactsecant}. This completes our proof.

\end{proof}

Next, we define the following geometric local size, which incorporates the coefficients of secant decomposition in Lemma \ref{lem:secant}. 
\begin{equation}\label{eq:Ddef}
D(G)=\sum_{i=1}^{p}\abs{a_i(G)}+\sum_{i=1}^{p}\norm{b_i(G)}+\sum_{j=1}^{m}\pi_jr_j^2.
\end{equation}

The next result establishes local equivalence between the
relative-density distance $\rho(G)$ and $D(G)$ under second-order
strong identifiability. The lower bound is a local inverse
inequality, while the upper bound follows from smoothness of
the kernel. Related inverse inequalities use Hellinger or
$L^1$ discrepancies between mixture densities to control
mixing-measure errors
\citep{Nguyen2013,HoNguyen2016,HeinrichKahn2018,ManoleHo2022}.

\begin{lemma}[Local equivalence of relative-density and mixing-measure discrepancies]
\label{lem:inverse}
There exist a $W_2$ neighborhood $\cU$ of $G_*$ in
$\cO_k(\Theta)$ and constants $0<c<C<\infty$ such that
\begin{equation}\label{eq:inverse}
    cD(G)\leq\rho(G)\leq CD(G),
    \qquad G\in\cU.
\end{equation}
\end{lemma}
\begin{proof}
The upper bound is immediate from \eqref{eq:exactsecant} and the relative derivative envelope. Indeed, with $Q_{i,r,u}$ as in Lemma~\ref{lem:secant}, the integral Taylor formula gives $\|Q_{i,r,u}/p_*\|_*\le\|F\|_2/2$, uniformly in its feasible arguments. Thus $\rho(G)\le C D(G)$.

Suppose the lower bound does not hold locally. Then for each $m$ there exists $G_m\ne G_*$ with $W_2(G_m,G_*)\le m^{-1}$, $D_m:=D(G_m)>0$, and $\rho(G_m)/D_m\to0$. Pad every measure to $k$ labeled atoms (allowing zero weights) and use the exact secant decomposition. Set
\[
A_{i,m}=a_i(G_m)/D_m,\quad B_{i,m}=b_i(G_m)/D_m,\quad
L_{j,m}=\pi_{j,m}r_{j,m}^2/D_m.
\]
By definition,
\begin{equation}\label{eq:normalizedCoefficientMassDetailed}
\sum_i|A_{i,m}|+\sum_i\|B_{i,m}\|+\sum_j L_{j,m}=1,\qquad L_{j,m}\ge0.
\end{equation}
After passing to a subsequence, we may assume that the nearest-center
indices $i_j$ are fixed and that all bounded coefficients converge
to limits denoted by $\bar A_i$, $\bar B_i$, and $\bar L_j$. Write $r_{j,m}=\|\theta_{j,m}-\theta_{i_j}^*\|$ and choose unit directions $u_{j,m}$ (arbitrarily when $r_{j,m}=0$) such that $\theta_{j,m}-\theta_{i_j}^*=r_{j,m}u_{j,m}$. Using the compactness of $[0,\operatorname{diam}(\Theta)]\times\mathbb S^{d-1}$, we can suppose that $r_{j,m}\to\bar r_j$ and $u_{j,m}\to\bar u_j$. The $L^2(P_*)$ continuity of $Q_{i,r,u}/p_*$ follows from its integral Hessian representation and the common $L^2(P_*)$ envelope. Therefore $\rho(G_m)/D_m\to0$ yields, in $L^2(P_*)$ and hence $\mu$-almost everywhere,
\begin{equation}\label{eq:normalizedsecantzero}
\sum_i\bar A_i f_{\theta_i^*}+\sum_i\bar B_i^\top\nabla f_{\theta_i^*}
+\sum_j\bar L_j Q_{i_j,\bar r_j,\bar u_j}=0.
\end{equation}
If $\bar r_j>0$, set $\eta_j=\theta_{i_j}^*+\bar r_j\bar u_j$. Nearest-center assignment ensures that $\eta_j$ cannot equal a different true center. Expanding $Q_{i_j,\bar r_j,\bar u_j}$, the coefficient of $f_{\eta}$ at each distinct nontrue location $\eta$ is the nonnegative sum $\sum_{j:\,\eta_j=\eta,\,\bar r_j>0}\frac{\bar L_j}{\bar r_j^2}.$ At $\bar r_j=0$, we have the contribution of the Hessian can be expressed as 
\[
\sum_{j:\,i_j=i,\,\bar r_j=0}\bar L_j \bar u_j^\top\nabla^2f_{\theta^*_i}\,\bar u_j = \sum_{j:\,i_j=i,\,\bar r_j=0}\bar L_j \tr(\bar u_j\bar u_j^\top\nabla^2f_{\theta^*_i}\,) = \tr(\bar C_i\nabla^2f_{\theta^*_i}),
\]
where $C_i = \sum_{j:\,i_j=i,\,\bar r_j=0}L_j\bar u_j\bar u_j^\top$. Therefore, by rewriting the identity \eqref{eq:normalizedsecantzero} based on above observation, we receive a linear combination of $f_{\theta_i^*}$, $\nabla f_{\theta_i^*}$, $\nabla^2f_{\theta_i^*}$, and the $f_{\eta_j}$ for each $\bar r_j>0$ such that this combination contains at most $p+k\le2k$ distinct locations of $\theta_i^*$ and $\eta_j$. Assumption~\ref{ass:strongid} implies that every kernel, gradient and Hessian coefficient vanishes. By considering the coefficient of $f_{\eta_j}$, as $\bar L_j \geq 0$, we have $\bar L_j=0$ for every $\bar r_j>0$. By plugging in back to \eqref{eq:normalizedsecantzero} and considering the contribution of Hessian matrix $\nabla^2 f_{\theta^*_i}$, we have $\operatorname{tr}(\bar C_i)=\sum_{j:i_j=i,\bar r_j=0}\bar L_j=0$ forces $\bar L_j=0$. The true-location coefficients then imply $\bar A_i=0$ and $\bar B_i=0$. This contradicts the limit of \eqref{eq:normalizedCoefficientMassDetailed}, whose left-hand side equals one. The local lower bound follows.
\end{proof}
The preceding lemma compares $\rho(G)$ with the coefficient
discrepancy $D(G)$. We next relate these quantities to the
squared Wasserstein distance. The first bound controls
$\rho(G)^2$ by $W_2^2(G,G_*)$, while the second bounds the
transportation cost through cluster mass discrepancies and
local second moments. Combined with the local equivalence
in Lemma~\ref{lem:inverse}, the latter yields
$W_2^2(G,G_*)\lesssim\rho(G)$.
\begin{lemma}\label{lem:w2bound}
\begin{enumerate}
    \item [(i)] There exists $c< \infty$ such that $c\rho^2(G) \leq W_2^2(G,G_*)$. 
    \item [(ii)] There exists $c,C<\infty$ such that for all $G$ sufficiently close to $G_*$ in Wasserstein metric, \begin{equation}\label{eq:w2D}
W_2^2(G,G_*) \le C\left\{\sum_i|a_i(G)|+\sum_j\pi_jr_j^2\right\}\le CD(G)\le C'\rho(G).
\end{equation}
\end{enumerate}

\end{lemma}
\begin{proof}
For the part $(i)$, let $\pi$ be a coupling between $G$ and $G_*$. We observe that for any couple of parameter $\theta$ and $\theta'$, from Assumption~\ref{ass:smooth},  
\begin{equation*}
    \left|\frac{f_\theta(x) - f_{\theta'}(x)}{p_*(x)}\right| \leq  F(x)\|\theta-\theta'\|
\end{equation*}
Therefore, we have 
\begin{equation*}
    \left|\dfrac{p_G(x) - p_{G'}(x)}{p_*(x)}\right| = \left|\int \dfrac{f_\theta(x)-f_{\theta'}(x)}{p_*(x)}\,d\pi(\theta,\theta')\right|\leq  F(x)\int\|\theta-\theta'\|\,d\pi(\theta,\theta').
\end{equation*}
Thus, by plugging in the optimal coupling with respect to 1-Wasserstein distance and using the integrability of $F$, we have
\begin{equation*}
    \rho(G) \lesssim W_1(G,G_*)\leq W_2(G,G_*). 
\end{equation*}
For part $(ii)$, from Lemma~\ref{lem:inverse}, it is sufficient to prove that there exists a constant $C$ satisfying
\begin{equation}
\label{eqn:dung_W_2_decomposition}
W_2^2(G,G_*) \leq  C\left\{\sum_i|a_i(G)|+\sum_j\pi_jr_j^2\right\}
\end{equation}
To prove this, we consider the following intermediate measure $G^c=\sum_{i=1}^p(\pi_i^*+a_i)\delta_{\theta_i^*}$. We can interprete the mass of each center $\theta_i^*$ to be $\sum_{j:i_j = i}\pi_i$. Using triangle inequality and Cauchy-Schwartz inequality, we have 
\begin{equation*}
    W_2^2(G,G_*)
\le 2(W_2^2(G,G^{c}) + W_2^2(G^{c},G_*)).
\end{equation*}

\textit{Estimate $W_2^2(G,G^{c})$.} Consider the transportation plan in which each fitted centers $\theta_j$ are transferred to their closet true center $\theta^*_{i_j}$. Then, we can bound $W_2^2(G,G^c)$ by the squared cost of this transportation 
\begin{equation}
\label{eqn:dung_W_2_G_G^c}
    W_2^2(G,G^c) \leq \sum_{j=1}^m\pi_j\|\theta_j - \theta^*_{i_j}\|^2. 
\end{equation}

\textit{Estimate $W_2^2(G^{c},G_*)$.} We take into account the transportation plan from $G^c$ to $G_*$ which keep a mass of $\min\{\pi_i^*+a_i,\pi_i^*\}$ for the center $\theta_i$, while the remaining mass can be transferred to any other centers. Then, we can show that the total remaining mass is equal to $\frac{1}{2}\sum_{i}|a_i|$, and as the cost to move the mass from one center to another one is no more than $\mathrm{diam}(\Theta)$, we have the following estimation 
\begin{equation}
\label{eqn:dung_W_2_G_star_G^c}
    W_2^2(G^{c},G_*) \leq \frac{1}{2}\mathrm{diam}(\Theta)^2\sum_{i}|a_i|. 
\end{equation}

\textit{Putting everything together.} Combining the estimations in \eqref{eqn:dung_W_2_decomposition}, \eqref{eqn:dung_W_2_G_G^c}, and \eqref{eqn:dung_W_2_G_star_G^c}, we achieve the estimation \eqref{eq:w2D}, completing the proof. 
\end{proof}
Once the inverse inequality in Lemma~\ref{lem:w2bound}(ii)
is established, obtaining a convergence rate for
$W_2^2(\widehat G_n,G_*)$ reduces to controlling the
relative-density discrepancy $\rho(\widehat G_n)$.
The following lemma establishes its root-$n$ rate, which
immediately yields the same rate for the squared Wasserstein error.

\subsection{Empirical-process properties and convergence rates}

The local comparisons established above reduce control of the
Wasserstein error to that of the relative-density discrepancy
$\rho(G)$. We now establish envelope and bracketing bounds for
the normalized relative-density class, yielding the Donsker
and Glivenko--Cantelli properties needed for likelihood
localization \citep{vdVW1996,vanDeGeer2000}. These properties
give $\rho(\widehat G_n)=O_P(n^{-1/2})$, and
Lemma~\ref{lem:w2bound}(ii) then yields
$W_2(\widehat G_n,G_*)=O_P(n^{-1/4})$. Although this Wasserstein
rate is known under suitable conditions
\citep{Nguyen2013,HoNguyen2016}, we derive it within our framework
to support the subsequent likelihood expansion and weak-limit
analysis.

\begin{lemma}[Compact normalized secants and Donsker property]\label{lem:donsker}
Let
\[
\cS=\left\{s_G=\frac{(p_G-p_*)/p_*}{\rho(G)}:G\in\cU,\;G\ne G_*\right\},
\]
where $\cU$ is the neighborhood in Lemma~\ref{lem:inverse}. Then:
\begin{enumerate}[label=(\roman*)]
\item $\cS$ admits an envelope $F_0\in L^{2+\eta}(P_*)$;
\item $\cS$ is totally bounded in $L^2(P_*)$ and for sufficiently small $\eps$,
\[
\log N_{[]}\!\left(\eps,\cS,L^2(P_*)\right)\le C\log(1/\eps),
\]
where $N_{[]}$ denotes the bracketing number. Consequently, $\cS$ is a relatively compact subset in $L^2(P_*)$;
\item $\cS$ is $P_*$-Donsker and $\cS^2=\{s^2:s\in\cS\}$ is $P_*$-Glivenko--Cantelli. In addition, the truncated class of $\cS^2$, defined as $\{s^2\ind(F_0\le K):s\in\mathcal{S}\}$, is also $P_*$-Glivenko--Cantelli for any $K \in \mathbb
R$
\end{enumerate}
\end{lemma}

\begin{proof}
Firstly, we prove part $(i)$. For each $G \in \mathcal{U}$, by dividing secant identity \eqref{eq:exactsecant} by $\rho(G)p_*$, we have 
\begin{equation*}
s_G = \frac{p_G-p_*}{\rho(G)p_*}
=\sum_{i=1}^{p}\frac{a_i(G)}{\rho(G)}\frac{f_{\theta_i^*}}{p_*}
+\sum_{i=1}^{p}\frac{b_i(G)^\top}{\rho(G)} \frac{\nabla f_{\theta_i^*}}{p_*}
+\sum_{j:r_j>0}\frac{\lambda_j(G)}{\rho(G)}\frac{Q_{i_j,r_j,u_j}}{p_*},
\end{equation*}
Noting that Lemma~\ref{lem:inverse} gives ${D(G)}/{\rho(G)}\le c^{-1}$, we have each element $s(G)$ can be expressed as a linear combination of ${f_{\theta_i^*}}/{p_*}$, ${\nabla f_{\theta_i^*}}/{p_*}$, and ${Q_{i_j,r_j,u_j}}/{p_*}$ with coefficients bounded by $c^{-1}$. The fixed kernel and gradient terms satisfy $|f_{\theta_i^*}/p_*|\le F$ and $\|\nabla f_{\theta_i^*}/p_*\|\le F$ thanks to Assumption~\ref{ass:smooth}. To study the boundedness of quadratic remainder $Q_{i,r,u}$, we employ the integral representation 
\begin{equation}
\label{eqn:quadratic_remainder_integral_representation}
    Q_{i,r,u}=\int_0^1(1-s)u^\top\nabla^2f_{\theta_i^*+sru}u\,ds.
\end{equation}
Assumption~\ref{ass:smooth} similarly gives us $|\nabla^2f_{\theta}/p_*| \leq F$, hence the integral representation implies $|Q_{i,r,u}/p_*| \leq F/2$. Together with the fact that all coefficients in $s_G$ are bounded, we have 
\begin{equation}\label{eq:SenvDetailed}
|s_G(x)|\le F_0(x):=C_0F(x),\qquad s_G\in\cS,
\end{equation}
for a fixed $C_0<\infty$, hence $F_0\in L^{2+\eta}(P_*)$. This means that $\mathcal{S}$ admits an envelope $F_0 \in L^{2+\eta}(P_*)$.

Next, we prove part $(ii)$. To establish the bound for the bracketing number, we first prove the Lipschitz property of quadratic remainder $Q_{i,r,u}$ with respect to $r$ and $u$. We observe that for any $r,u,r',u'$ and $s\in(0,1)$, we have 
\begin{align*}
    |u^{\top}\nabla^2f_{\theta_i^*+sru}u- (u')^{\top}\nabla^2f_{\theta_i^*+sr'u'}u'| &\leq |u^{\top}(\nabla^2f_{\theta_i^*+sru}-\nabla^2f_{\theta_i^*+sr'u'})u|\\
    &\hspace{1cm} + |u^{\top}\nabla^2f_{\theta_i^*+sr'u'}u-(u')^{\top}\nabla^2f_{\theta_i^*+sr'u'}u'|.
\end{align*}
For the first term, the mean-value theorem and the third-derivative envelope give
\[
\bigl|u^{\top}(\nabla^2f_{\theta_i^*+sru}-\nabla^2f_{\theta_i^*+sr'u'})u\bigr|
\le p_*F\,s(|r-r'|+R\|u-u'\|).
\]
For the second term, using the fact that $\|\nabla^2f_{\theta}/p_*\|\leq F$ and $\|u\|=\|u'\|=1$ give
\[
\bigl|u^{\top}\nabla^2f_{\theta_i^*+sr'u'}u
 -(u')^{\top}\nabla^2f_{\theta_i^*+sr'u'}u'\bigr|
\le2p_*F\|u-u'\|.
\]
% The two contributions are controlled by different derivative orders; keeping their roles separate gives the stated Lipschitz estimate.
Consequently, using the integral representation in equation \eqref{eqn:quadratic_remainder_integral_representation} and Assumption~\ref{ass:smooth}, we have 
\begin{align*}
    \left|\frac{Q_{i,r,u} - Q_{i,r',u'}}{p_*}\right| &\leq \frac{1}{p_*}\int_{0}^1(1-s)|u^{\top}\nabla^2f_{\theta_i^*+sru}u- (u')^{\top}\nabla^2f_{\theta_i^*+sr'u'}u'|\,ds \\
    &\leq F\left[\|u-u'\|\int_{0}^12(1-s)\,ds + (|r-r'| + R\|u-u'\|)\int_{0}^1s\,ds\right],
\end{align*}
therefore, there is a constant $C$ such that the following inequality holds point-wise
\begin{equation}\label{eq:QLipschitzDetailed}
\left|\frac{Q_{i,r,u}-Q_{i,r',u'}}{p_*}\right|
\le C F\{|r-r'|+\|u-u'\|\}.
\end{equation}

We continue by proving the total boundedness and estimation for bracket number of $\mathcal{S}$. Consider the following compact set 
\[
\mathcal R
=\overline{\{(r,u):r>0,\ u\in\mathbb S^{d-1},\exists i:\ \theta_i^*+ru\in\Theta\}}
\subset[0,R]\times\mathbb S^{d-1}.
\]
After padding to $k$ labels, for each second-order secant decomposition (cf. equation~\eqref{eq:exactsecant}) of an element $s_G \in \mathcal{S}$, the scalar coefficients are in $[-c^{-1},c^{-1}]^p$, the gradient coefficients are in $[-c^{-1},c^{-1}]^{pd}$, and the $k$ nonnegative secant coefficients are bounded by $c^{-1}$. Therefore, the resulting parameter set is contained in a bounded subset $\widetilde\Theta$ of a fixed Euclidean space of dimension $D=D(p,k,d)$, and the associated functions may be written $s_\vartheta$, for $\vartheta \in \widetilde\Theta$. Then, using linearity in coefficient, and the pointwise Lipschitz property of $f_{\theta_i^*}/p_*$, $\nabla f_{\theta_i^*}/p_*$ and quadratic remainder $Q_{i,r,u}/p_*$ in \eqref{eq:QLipschitzDetailed}, there is an envelope $L\in L^{2+\eta}(P_*)$ such that
\begin{equation}\label{eq:parametricPointwiseLip}
|s_\vartheta(x)-s_{\vartheta'}(x)|
\le L(x)\|\vartheta-\vartheta'\|_{\mathrm E}.
\end{equation}

For each assignment pattern, take a Euclidean $\delta$-net $\mathcal C$ of the compact parameter set $\widetilde\Theta$, with at most $C\delta^{-D}$ centers. Then, the union of brackets $\{[s_{\vartheta}-\delta L, s_{\vartheta}+\delta L],\vartheta \in \mathcal{C}\}$, such that the $L^2(P_*)$ sizes of each bracket is exactly $2\delta\|L\|_2$, contains all elements of $\mathcal{S}$. By choosing $\delta = \varepsilon/(2\|L\|_2)$, we achieve that 
\begin{equation}\label{eq:bracketEntropyDetailed}
N_{[]}\!\left(\eps,\cS,L^2(P_*)\right)\le C\eps^{-D},
\qquad
\log N_{[]}\!\left(\eps,\cS,L^2(P_*)\right)\le C_1+C_2\log(1/\eps).
\end{equation}
This proves total boundedness of $\cS$ in $L^2(P_*)$ as well as upper bound for bracket number. As $L^2(P_*)$ is a complete space, $\cS$ is relatively compact. Using similar argument, we also have all the result above also holds for $L^1(P_*)$ space with $L^1$ norm. 

Lastly, we prove part $(iii)$. To prove that $\cS$ is $P_*$-Donsker, we start with the observation that thanks to equation \eqref{eq:bracketEntropyDetailed}, we have 
\begin{equation*}
    \int_0^1
\sqrt{\log N_{[]}\!\left(\eps\|F_0\|_2,\cS,L^2(P_*)\right)}\,d\eps \lesssim \int_{0}^1\sqrt{\log(1/\varepsilon)}\,d\varepsilon < \infty.
\end{equation*}
Together with the square-integrable envelope $F_0$, the bracketing central limit theorem yields that $\cS$ is $P_*$-Donsker, see Chapter~2.5 of \citet{vdVW1996}. In particular,
\[
\sup_{s\in\cS}|\mathbb G_ns|=O_P(1),
\]
and the empirical process is asymptotically uniformly continuous in the $L^2(P_*)$ semimetric on the closure of $\cS$.

We show that $\cS^2$ and its truncated class are $P_*$-Glivenko-Cantelli. For two parameter $\vartheta,\vartheta'\in \widetilde\Theta$, we have
\[
|s_\vartheta^2-s_{\vartheta'}^2|
\le (|s_\vartheta|+|s_{\vartheta'}|)
|s_\vartheta-s_{\vartheta'}|
\le 2F_0L\|\vartheta-\vartheta'\|_{\mathrm E}.
\]
For every center $s_{\vartheta_0}$ in the finite Euclidean $\delta$-net $
\mathcal{C}$, the preceding inequality gives valid pointwise brackets
$[s_{\vartheta_0}^2-2\delta F_0L,\,s_{\vartheta_0}^2+2\delta F_0L]$
for the squares of functions indexed by its $\delta$-ball. Their $L^1(P_*)$ width is at most $4\delta P_*(F_0L)$, which is finite due to the fact that $F_0,L\in L^2(P_*)$. Hence $\cS^2$ has finite $L^1(P_*)$ bracketing numbers at every positive accuracy. In addition, its envelope $F_0^2$ is integrable, thus, using the bracketing Glivenko--Cantelli theorem, we have 
\[
\sup_{s\in\cS}|(P_n-P_*)s^2|\pto0.
\]

Multiplying the functions by a fixed indicator such as $\ind(F_0\le K)$ only decreases the bracket sizes, so the same conclusion holds for the truncated squared class used in Lemma~\ref{lem:rootn}. This establishes the stated envelope, bracketing, Donsker, and Glivenko--Cantelli properties.
\end{proof}

% We verify the empirical-process assertions through an explicit bracketing construction. This is slightly stronger than an $L^2(P_*)$ covering-number calculation and avoids using a covering bound where a bracketing or uniform-entropy condition is required; see Chapters~2.4--2.5 of \citet{vdVW1996}.

% This is a localized likelihood-curvature argument of the type used in empirical-process analyses of $M$-estimators; see \citet{vdV1998,vdVW1996}.  We extract the quadratic term only on a truncated set where the relative perturbation is uniformly small.  In particular, the proof never requires uniform convergence of $p_G/p_*$ over the observation space.

% \section{Exact secants, local inversion, and root-$n$ density localization}\label{sec:proof1}
% We next turn to the empirical-process arguments.  Compact $M$-estimation, bracketing central limit theorems, stochastic equicontinuity, and local likelihood expansion are standard tools; for systematic treatments, see \citet{vdV1998,vdVW1996,Dudley1999}. We must handle the indexing class carefully because normalized mixture secants approach singular configurations. In the proofs below, we therefore construct the relevant envelopes and brackets directly instead of relying on a generic finite-dimensional parameterization.

\begin{lemma}[Root-$n$ density localization]\label{lem:rootn}
Under the standing assumptions, the likelihood-selected relative density satisfies $\rho(\widehat G_n)=O_P(n^{-1/2})$.
Consequently, the squared Wasserstein loss satisfies
\[
W_2^2(\widehat G_n,G_*)=O_P(n^{-1/2}).
\]
\end{lemma}

\begin{proof}
\textit{Step 1: A local quadratic upper bound for the log-likelihood.} In this step, we establish an upper bound for $\ell_n(G)-\ell_n(G_*)$ based on $\rho(G)$ for $G \in \mathcal{U}$, where $\mathcal{U}$ is the neighbor in Lemma~\ref{lem:inverse}, shrinking it as needed. By definition, $P_*s_G^2=1$. Additionally, since both $p_G$ and $p_*$ integrate to one,
\begin{equation}\label{eq:meanZeroSecantDetailed}
P_*s_G=\rho(G)^{-1}\int(p_G-p_*)\,d\mu=0,
\end{equation}
Consider $F_0$ be the envelope from Lemma~\ref{lem:donsker}. As $F_0^2$ is integrable, we can choose a positive $K<\infty$ such that $P_*\{F_0^2\ind(F_0>K)\}<\frac14$. Then, for every $s\in\cS$, $s^2\le F_0^2$, we have 
\[
P_*\{s^2\ind(F_0\le K)\}
=1-P_*\{s^2\ind(F_0>K)\}
\ge\frac34.
\]
From part $(iii)$ of Lemma \ref{lem:donsker} about the Glivenko--Cantelli property of the truncated squared class, 
\[
\sup_{s\in\cS}
\left|(P_n-P_*)\{s^2\ind(F_0\le K)\}\right|\pto0.
\]
Therefore the event
$\mathcal E_n:=\left\{
\inf_{s\in\cS}P_n\{s^2\ind(F_0\le K)\}\ge\frac12
\right\}$ has probability tending to one. 

Using the fact that $c\rho^2(G) \leq W_2^2(G,G_*)$ proved in Lemma~\ref{lem:w2bound}, it is possible to shrink the neighbor $\mathcal{U}$ such that in $\mathcal{U}$, as $W_2(G,G_*)\leq m$, then $\rho(G)\lesssim \sqrt{m}$. By shrinking neighbor $\mathcal{U}$ if necessary, we can suppose that $\rho(G)K \leq 1/2$ for all $G \in \mathcal{U}$. On $\{F_0\le K\}$, we have 
\[
\left|\frac{p_G - p_*}{p_*}\right|=\rho(G)|s_G|\le\rho(G)K\le\frac12. 
\]
Applying the inequality $\log(1+y) \leq y - \frac{1}{4}y^2\ind(y\leq \frac{1}{2})$, we have 
\begin{equation*}
    \log\left(1+\frac{p_G-p_*}{p_*}\right) = \log\left(1+\rho(G)s(G)\right) \leq \rho(G)s(G) -\frac{1}{4}\rho(G)^2s(G)^2\ind(F_0\leq K). 
\end{equation*}
By taking the expectation with respect to empirical measure $P_n$, we have 
\[
\ell_n(G)-\ell_n(G_*)
\le \sqrt n\,\rho(G)\mathbb G_ns_G
-\frac{1}{4}n\rho(G)^2P_n\{s_G^2\ind(F_0\le K)\}.
\]
On the event $\mathcal{E}_n$, by letting $c = 1/8$, this inequality becomes 
\begin{equation}\label{eq:rootnMasterDetailed}
\ell_n(G)-\ell_n(G_*)
\le \sqrt n\rho(G)\mathbb G_ns_G-c n\rho(G)^2
\end{equation}
\textit{Step 2:} We prove that $\sqrt n\rho(\widehat G_n) = O_P(n^{-1/2})$. Lemma~\ref{lem:consistency} gives $\rho(\widehat G_n)\pto0$. Since $K$ is fixed, $\PP\{\rho(\widehat G_n)K\le1/2\}\to1$. Moreover, $\widehat G_n$ belongs to the neighborhood $\cU$ with probability tending to one. Hence using the estimation~\eqref{eq:rootnMasterDetailed}, on an event whose probability tends to one, we have 
\begin{equation*}
\ell_n(\widehat G_n)-\ell_n(G_*)
\le \sqrt n\rho(\widehat G_n)\mathbb G_ns_{\widehat G_n}-c n\rho(\widehat G_n)^2
\end{equation*}
Moreover, by the definition of the maximum likelihood estimator, we have $0\le \ell_n(\widehat G_n)-\ell_n(G_*)$. Therefore, by defining $Z_n:=\sup_{s\in\cS}|\mathbb G_ns|$, on an event whose probability tends to one, 
\begin{equation*}
    0\leq  \sqrt n\rho(\widehat G_n)Z_n-cn\rho(\widehat G_n)^2,
\end{equation*}
which implies $\sqrt n\rho(\widehat G_n) \leq c^{-1}Z_n$. Furthermore, Lemma~\ref{lem:donsker} implies $Z_n=O_P(1)$, which leads to $\rho(\widehat G_n)=O_P(n^{-1/2})$. 

For the rate of $W_2^2(\widehat G_n,G_*)$, Lemma~\ref{lem:w2bound} gives 
\[
W_2^2(\widehat G_n,G_*)\le C\rho(\widehat G_n)=O_P(n^{-1/2}).
\]
This proves both assertions.
\end{proof}

\section{Proofs for Section~\ref{sec:densitylimit}}
\label{app:densityproofs}

\subsection{Properties of attainable tangent set $\mathcal{T}_k$}
\label{sec:prop_of_attainable_tangent_set}
%\section{Characterization of the attainable tangent set}\label{sec:proof2}
In this section, we prove the geometric properties of the
attainable tangent set $\cT_k$. Proposition~\ref{prop:tangent}
characterizes its elements and establishes uniqueness of their
canonical representations, while Lemma~\ref{lem:tangentcompact}
shows that each bounded section $\cT_{k,M}$ is compact in
$L^2(P_*)$. These results describe the structure of $\cT_k$
and provide key ingredients for proving the density weak limit
in Theorem~\ref{thm:density}.

% \begin{proposition}[Attainable second-order tangent set]\label{prop:tangent}
% A function $h\in L^2(P_*)$ belongs to $\cT_k$ if and only if there exist $t_m\downarrow0$ and $G_m\in\cO_k(\Theta)$ such that
% \begin{equation}\label{eq:tangentseq}
% \left\|\frac{p_{G_m}-p_*}{t_mp_*}-h\right\|_*\to0.
% \end{equation}
% Moreover, under Assumption~\ref{ass:strongid}, the canonical representation $h=h_\gamma$ is unique up to permutation of the atoms of $\nu$.
% \end{proposition}

\begin{proof}[Proof of Proposition~\ref{prop:tangent}]
The proof has three parts: necessity, sufficiency, and uniqueness of the canonical representation. The geometric approach is related to tangent-cone methods for constrained estimation \citep{Geyer1994} and locally conic
parametrizations of mixture models \citep{DacunhaGassiat1999}. Here, we explicitly track the number of fitted atoms to derive the component-budget constraint \eqref{eq:budget}.

\medskip\noindent\textit{Part 1: Necessity.} Consider a sequence of measures $\{G_m\}_{m} \subset \mathcal{O}_k(\Theta)$ along with $t_m\downarrow 0$ satisfying $\|(p_{G_m}-p_*)/(t_mp_*) -h\|_{*}\to 0$ for some $h$. We construct a quadruple $\gamma = (\alpha,\beta,C,\nu)$ such that 
\begin{equation*}
h(x)
=h_{\gamma}(x) = 
\frac{1}{p_*(x)}
\left[
\sum_{i=1}^{p}
\left\{
\alpha_i f_{\theta_i^*}(x)
+\beta_i^\top\nabla f_{\theta_i^*}(x)
+\frac12\tr\!\left(C_i\nabla^2f_{\theta_i^*}(x)\right)
\right\}
+\int f_\theta(x)\,d\nu(\theta)
\right]
\end{equation*}
and show that the quadruple $\gamma$ satisfies the conditions in \eqref{eq:masscon} and \eqref{eq:budget}. 

\textit{Step 1: Construction of quadruple $\gamma$.} 
From the assumption about $G_m$, we have 
\[
\rho(G_m)
=\left\|\frac{p_{G_m}-p_*}{p_*}\right\|_*
=t_m\left\|h+o_{L^2}(1)\right\|_*
=O(t_m).
\]
Therefore, $\rho(G_m)\to0$, which leads to $G_m\Rightarrow G_*$. For all sufficiently large $m$, $G_m$ lies in the local neighborhood $\mathcal{U}$ in Lemma~\ref{lem:inverse}, so that  $D(G_m)=O(t_m)$. In particular,
\begin{equation}
\label{eqn:dung_a_b_C_consistency}
    \sum_i|a_i(G_m)|=O(t_m),\qquad
\sum_i\|b_i(G_m)\|=O(t_m),\qquad
\sum_j\pi_{j,m}r_{j,m}^2=O(t_m).
\end{equation}
By adding zero-weight atoms, we can pad to exactly $k$ labels . Since $\Theta^k$ is compact, every subsequence has a further subsequence for which $\theta_{j,m}\to\bar\theta_j\in\Theta$ for $j=1,\ldots,k$. We call a label $j$ to be \emph{local to $i$} if $\bar\theta_j=\theta_i^*$, and call it \emph{non-local} if $\bar\theta_j$ is not a true atom. Because the true support is finite and separated, an atom local to $i$ is eventually assigned by the nearest-center map to $i$. Conversely, an atom whose limit is another true center cannot remain assigned to $\theta^*_i$ for large $m$.

Now we construct the measure $\nu$, which stands for the non-local component of $G_m$. Fix a nonlocal label $j$, define 
$\delta_j:=\min_i\|\bar\theta_j-\theta_i^*\|>0$ as its limiting distance from the true support.  For large $m$, $r_{j,m}\ge\delta_j/2$. Since every term in the nonnegative sum $\sum_j\pi_{j,m}r_{j,m}^2$ is bounded by $D(G_m)$, it follows that 
\[
\pi_{j,m}\frac{\delta_j^2}{4}
\le \pi_{j,m}r_{j,m}^2
=O(t_m),
\]
so $\pi_{j,m}=O(t_m)$. After another subsequence, ${\pi_{j,m}}/{t_m}\to c_j\ge0$. By merging labels having the same nontrue limiting location and retaining only the positive limiting coefficients, we obtain a finite positive measure $\nu=\sum_{\ell=1}^q c_\ell\delta_{\eta_\ell}$, where $c_\ell>0$, with pairwise distinct nontrue $\eta_\ell$. 

The next step is to construct the component  $\alpha$, $\beta$, $C$ of the quadruple $\gamma$. Let $L_i$ be the set of labels local to $\theta_i^*$ and define $M_{i,m}=\sum_{j\in L_i}\pi_{j,m}$. Because $G_m\Rightarrow G_*$ and the true atoms are separated, choose disjoint balls $B_i$ centered at $\theta_i^*$ containing no other true atom. Eventually all local labels in $L_i$ lie in $B_i$, all labels local to another true atom lie outside $B_i$, and all genuinely nonlocal labels have total mass $O(t_m)=o(1)$. Weak convergence therefore gives
\[
M_{i,m}\to G_*(\{\theta_i^*\})=\pi_i^*>0.
\]
In particular, every $L_i$ is nonempty for $m$ sufficiently large. Define
\[
\alpha_{i,m}=\frac{1}{t_m}(M_{i,m}-\pi_i^*),
\]
\[
\beta_{i,m}=\frac1{t_m}
\sum_{j\in L_i}\pi_{j,m}(\theta_{j,m}-\theta_i^*),
\]
\[
C_{i,m}=\frac1{t_m}
\sum_{j\in L_i}\pi_{j,m}(\theta_{j,m}-\theta_i^*)(\theta_{j,m}-\theta_i^*)^\top.
\]
We show that $\{a_{i,m}\}_m$, $\{\beta_{i,m}\}_m$, and $\{C_{i,m}\}_m$ are bounded. For $\{a_{i,m}\}_m$, we have 
\begin{equation*}
   a_{i}(G_m) = \sum_{j:i_j =i}\pi_{j,m} - \pi_i^* = \underbrace{t_m \alpha_{i,m}}_{\text{local component}} + \underbrace{\sum_{\substack{j:i_j = i\\ i \notin L_i}}\pi_{j,m}}_{\text{non-local component}}. 
\end{equation*}
From the construction of $\nu$, we have $\pi_{j,m} = \mathcal{O}(t_m)$ for non-local index $j$. Therefore, as $a_{i}(G_m) = \mathcal{O}(t_m)$, we have $\alpha_{i,m} = \mathcal{O}(1)$, which is bounded. Using similar argument, we also have that $\beta_{i,m}$ is bounded. For the second-moment components $C_{i,m}$ are positive semidefinite matrices whose traces are bounded because
\[
\sum_i\operatorname{tr}(C_{i,m})
\le \frac{1}{t_m}\sum_j\pi_{j,m}r_{j,m}^2=O(1).
\]
As a consequence, $C_{i,m}$ are bounded. 
By passing to a subsequence, there exist a triple $(\alpha,\beta,C)$ such that 
\[
\alpha_{i,m}\to\alpha_i,
\qquad
\beta_{i,m}\to\beta_i,
\qquad
C_{i,m}\to C_i\succeq0.
\]

\textit{Step 2: Prove that $h=h_{\gamma}$}. For $\gamma$ established in previous step, we show that $h_{\gamma}$ is exactly the $h$ such that $\|(p_{G_m}-p_*)/(t_mp_*)-h\|_{*} \to 0$. For $j\in L_i$, write $\Delta_{j,m}=\theta_{j,m}-\theta_i^*$, then $\max_{j\in L_i}\|\Delta_{j,m}\|\to0$. Second-order Taylor's formula along the segment lying $\theta_i^*$ and $\theta_{j,m}$ gives
\begin{equation}
\label{eqn:dung_second_order_Taylor_expansion_for_f_theta_star}
    f_{\theta_i^*+\Delta_{j,m}}
=f_{\theta_i^*}+\Delta_{j,m}^\top\nabla f_{\theta_i^*}
+\frac12\operatorname{tr}(\Delta_{j,m}\Delta_{j,m}^\top\nabla^2f_{\theta_i^*})
+R_{j,m}.
\end{equation}
Here, by using Assumption~\ref{ass:smooth}, the pointwise relative bound for the remainder $R_{j,m}$ can be expressed as 
\[
\left|\frac{R_{j,m}(x)}{p_*(x)}\right|
\le C F(x)\|\Delta_{j,m}\|^3.
\]
Therefore, after weighting and summing, and combining with estimation \eqref{eqn:dung_a_b_C_consistency}, we have 
\[
\begin{aligned}
\left\|\frac{\sum_{j\in L_i}\pi_{j,m}R_{j,m}}{t_mp_*}\right\|_*
&\le C\|F\|_2
\left(\max_{j\in L_i}\|\Delta_{j,m}\|\right)
\frac{\sum_{j\in L_i}\pi_{j,m}\|\Delta_{j,m}\|^2}{t_m}
\to 0,
\end{aligned}
\]
Consequently, the local component sum converges in $L^2(P_*)$: 
\[
\frac1{t_mp_*}\left\{
\sum_{j\in L_i}\pi_{j,m}f_{\theta_{j,m}}-\pi_i^*f_{\theta_i^*}
\right\}
\to
\frac{\alpha_if_i+\beta_i^\top\nabla f_{\theta_i^*}+\frac12\operatorname{tr}(C_i\nabla^2f_{\theta_i^*})}{p_*}
\]
For nonlocal labels, consider the center $\theta_{j,m}\to\eta$ and $\pi_{j,m}/t_m\to c$, then we have 
\[
\left\|
\frac{\pi_{j,m}}{t_m}\frac{f_{\theta_{j,m}}}{p_*}
-c\frac{f_\eta}{p_*}
\right\|_*
\to0.
\]
Indeed, the coefficient converges and $\theta\mapsto f_\theta/p_*$ is $L^2(P_*)$ continuous by the derivative envelope. Summing the finitely many nonlocal labels gives the term $\int f_\theta\,d\nu(\theta)/p_*$. Combining the local and nonlocal limits proves $h=h_\gamma$ with the stated canonical form.

\textit{Step 3: Verify the conditions \eqref{eq:masscon} and \eqref{eq:budget} for $\gamma$.} For the mass conservation condition in \eqref{eq:masscon}, noting that for $m$ sufficiently large,
\[
\sum_iM_{i,m}+\sum_{j\in L_0}\pi_{j,m}=1,
\]
where $L_0$ is the set of nonlocal labels. Subtracting $\sum_i\pi_i^*=1$ from the preceding identity, dividing by $t_m$, and passing to the limit, noting that the non-local labels produce the sum of mass of $\nu$ and local label $(M_{i,m}-\pi_i^*)/t_m\to \alpha_i$, we obtain mass conservation in \eqref{eq:masscon}: 
\[
\sum_i\alpha_i+\sum_{\ell=1}^q c_\ell=0,
\]
For rank condition \eqref{eq:budget}, the definition 
\begin{equation*}
    C_{i,m}=\frac1{t_m}
\sum_{j\in L_i}\pi_{j,m}(\theta_{j,m}-\theta_i^*)(\theta_{j,m}-\theta_i^*)^\top.
\end{equation*}
directly yields $\rank(C_{i,m})\leq r_i$ for  $r_i=|L_i|$, which leads to $\rank(C_i) \leq r_i$. We establish the sharper rank bound that $\rank(C_i) \leq r_i-1$. Define the weighted local barycenter
\[
\bar\Delta_{i,m}
=\frac{\sum_{j\in L_i}\pi_{j,m}\Delta_{j,m}}{M_{i,m}}
=\frac{t_m\beta_{i,m}}{M_{i,m}}
\]
and the weighted centered covariance
\[
V_{i,m}
=\sum_{j\in L_i}\pi_{j,m}
(\Delta_{j,m}-\bar\Delta_{i,m})
(\Delta_{j,m}-\bar\Delta_{i,m})^\top.
\]
Let $r_{i,m}^+$ denote the number of labels in $L_i$ having positive weight at sample size $m$ such that they are only the labels contributing to $V_{i,m}$. Their centered displacement vectors satisfy the weighted relation
\[
\sum_{j\in L_i:\,\pi_{j,m}>0}\pi_{j,m}(\Delta_{j,m}-\bar\Delta_{i,m})=0.
\]
Hence their span has dimension at most $r_{i,m}^+-1$, and therefore $\rank(V_{i,m})\le r_{i,m}^+-1\le r_i-1$. Direct expansion gives
\[
C_{i,m} = \frac{V_{i,m}}{t_m} + \frac{t_m\beta_{i,m}\beta_{i,m}^\top}{M_{i,m}}.
\]
Since $M_{i,m}\to\pi_i^*>0$ and $\beta_{i,m}$ is bounded, the second term tends to zero, therefore $V_{i,m}/t_m\to C_i$. In addition, the set $\{A:\rank(A)\le r_i-1\}$ is closed, thus, we obtain the rank bound $\rank(C_i)\le r_i-1.$

The sets $L_1,\ldots,L_p$ and the nonlocal labels are disjoint subsets of the $k$ labels. Each distinct positive-mass nonlocal limit $\eta_\ell$ requires at least one nonlocal label. Hence, $q+\sum_{i=1}^p r_i\le k$.
Using the rank bound above, we obtain the budget restriction \eqref{eq:budget}: 
\[
q+\sum_i\rank(C_i)
\le q+\sum_i(r_i-1)
\le k-p.
\]
This completes the necessity argument.

\medskip\noindent\textit{Part 2: Sufficiency.}
Now fix a canonical configuration $\gamma=(\alpha,\beta,C,\nu)$ satisfying the conditions \eqref{eq:masscon} and \eqref{eq:budget}, and denote $q$ the number of support point of $\nu$. We construct a family of discrete measures $(G_t)_{t\downarrow 0}\subset \mathcal{O}_k(\Theta)$ such that $\|(p_{G_t}-p_{G_*})/(tp_{G_*}) - h_{\gamma}\|_* \to 0$. 

Let $r_i=\rank(C_i)$. If $r_i=0$, take one point $v_{i,1}=0$ with weight $1$. Suppose $r_i\ge1$. Choose a factorization ${C_i}/{\pi_i^*}=A_iA_i^\top$,
for $A_i\in\RR^{d\times r_i}$ with full column rank. Choose \(z_{i,1},\ldots,z_{i,r_i+1}\) as the vertices of a regular simplex in \(\mathbb R^{r_i}\), centered at the origin and scaled to have squared norm \(r_i\). With equal weights \(\omega_{i,h}=1/(r_i+1)\), these vectors satisfy
\[
\sum_{h=1}^{r_i+1}\omega_{i,h}z_{i,h}=0,
\qquad
\sum_{h=1}^{r_i+1}\omega_{i,h}z_{i,h}z_{i,h}^\top=I_{r_i}.
\]
By setting $v_{i,h}=A_iz_{i,h}$, we have 
\begin{equation}\label{eq:simplexDetailed}
\sum_{h=1}^{r_i+1}\omega_{i,h}v_{i,h}=0,
\qquad
\sum_{h=1}^{r_i+1}\omega_{i,h}v_{i,h}v_{i,h}^\top=C_i/\pi_i^*.
\end{equation}
For sufficiently small $t>0$, define the local cluster mass $M_i(t)=\pi_i^*+t\alpha_i>0$
and locations
\begin{equation}
\label{eqn:expansion_of_fitted_center}
\theta_{i,h}(t)
=\theta_i^*+\sqrt t\,v_{i,h}+t\frac{\beta_i}{\pi_i^*}.
\end{equation}
When $t$ is small enough, as $\theta_i^*$ lies in the interior of $\Theta$ and there are finitely many $v_{i,h}$, all these points $\theta_{i,h}(t)$ belong to $\Theta$. Thus, the measures $G_t$ can be constructed by giving $\theta_{i,h}(t)$ mass $M_i(t)\omega_{i,h}$ and for each nonlocal atom of $\nu$, putting mass $tc_\ell$ at $\eta_\ell$: 
\begin{equation*}
    G_t = \sum_{i=1}^{p}\left(\sum_{h=1}^{r_i+1}M_i(t)\omega_{i,h}\delta_{\theta_{i,h}(t)}\right) + t\nu. 
\end{equation*}
By \eqref{eq:masscon}, the total mass of $G_t$ is exactly equal to 
\[
\sum_iM_i(t)+t\sum_\ell c_\ell
=1+t\left(\sum_i\alpha_i+\sum_\ell c_\ell\right)=1.
\]
In addition, the number of local support points is $\sum_i(r_i+1)=p+\sum_i\rank(C_i)$, so the total number of support points is at most
\[
p+\sum_i\rank(C_i)+q\le k.
\]
Therefore, we have $G_t\in\cO_k(\Theta)$.

Now we show that $\|(p_{G_t}-p_{G_*})/(tp_{G_*}) - h_{\gamma}\|_* \to 0$. Let $\Delta_{i,h}(t)=\theta_{i,h}(t)-\theta_i^*$. We apply the Taylor expansion for the local component as in the equation \eqref{eqn:dung_second_order_Taylor_expansion_for_f_theta_star} and evaluate the contribution of weighted zero-, first-, second-order components and remainders. Firstly, we observe that the local displacements satisfy $\max_h\|\Delta_{i,h}(t)\|=O(\sqrt t)$. For the first-order displacement, it follows from the equation \eqref{eq:simplexDetailed} and the definition of $\theta_{i,h}$ that as $t\downarrow 0$, 
\[
\begin{aligned}
\sum_hM_i(t)\omega_{i,h}\Delta_{i,h}(t)
=M_i(t)\left\{\sqrt t\sum_h\omega_{i,h}v_{i,h}
+t\frac{\beta_i}{\pi_i^*}\right\}
=t\beta_i+t^2\frac{\alpha_i}{\pi_i^*}\beta_i
=t\beta_i+o(t).
\end{aligned}
\]
Again using the centering in \eqref{eq:simplexDetailed}, noting that all $t^{3/2}$ cross terms vanish after summation, we have 
\[
\begin{aligned}
\sum_{h=1}^{r_i+1}M_i(t)\omega_{i,h}\Delta_{i,h}(t)\Delta_{i,h}(t)^\top =M_i(t)t\sum_h\omega_{i,h}v_{i,h}v_{i,h}^\top+O(t^2)=t\frac{M_i(t)}{\pi_i^*}C_i+O(t^2)
=tC_i+o(t).
\end{aligned}
\]
For the weighted third-order remainder, we have 
\[
\left\|
\frac{1}{tp_*}
\sum_hM_i(t)\omega_{i,h}R_{i,h}(t)
\right\|_*
\le C\|F\|_2\frac{O(t^{3/2})}{t}
=O(\sqrt t)\to0.
\]
The nonlocal atoms contribute exactly $t\int f_\theta\,d\nu(\theta)$. Combining the zeroth-, first-, and second-order calculations gives
\[
\left\|\frac{p_{G_t}-p_*}{tp_*}-h_\gamma\right\|_*\to0.
\]
\medskip\noindent\textit{Part 3: Uniqueness of the canonical representation.}

\par

Suppose that there exist two tangent configurations $\gamma = (\alpha,\beta,C,\nu)$ and $\widetilde\gamma = (\widetilde \alpha,\widetilde\beta,\widetilde C,\widetilde\nu)$ such that $h_\gamma=h_{\widetilde\gamma}$. By multiplying by $p_*$ and regrouping the term, we have 
\begin{equation}
\label{eqn:dung_sum_of_independent_terms}
\sum_{i=1}^{p}\left[(\alpha_i - \widetilde\alpha_i)f_{\theta_i^*} + (\beta_i - \widetilde\beta_i)^{\top}\nabla f_{\theta_i^*}  + \frac{1}{2}\tr\left((C_i-\widetilde C_i)\nabla^2f_{\theta_i^*} \right)\right]+  \int f_{\theta}(x)d(\nu-\widetilde\nu)(\theta) = 0. 
\end{equation}
Let $q$ and $\widetilde q$ denote the numbers of support points of $\nu$ and $\widetilde\nu$, respectively. The budget constraint implies $q,\widetilde q\le k-p$. Therefore, the number of distinct locations in the identity \eqref{eqn:dung_sum_of_independent_terms} is at most
\[
p+q+\widetilde q
\le p+2(k-p)
=2k-p\le2k.
\]
Assumption~\ref{ass:strongid} forces every coefficient in the difference identity to vanish. Hence
\[
\alpha_i=\widetilde\alpha_i,
\qquad
\beta_i=\widetilde\beta_i,
\qquad
C_i=\widetilde C_i
\]
for every true atom, and the signed difference of the two nonlocal finite measures $\nu-\widetilde\nu$ is zero. Therefore $\nu=\widetilde\nu$ as measures, which means equality up to permutation of their listed atoms. This proves uniqueness.
\end{proof}

\begin{lemma}[Compactness of bounded attainable tangent set]\label{lem:tangentcompact}
For every $M<\infty$, the set
\[
\cT_{k,M} := \cT_k\cap\{h:\norm{h}_*\le M\}
\]
is compact in $L^2(P_*)$.
\end{lemma}

\begin{proof}
We first show that the unit attainable tangent section $\mathcal{S}_{k} = \{h\in \cT_k:\|h\|_* = 1\}$ is relatively compact. Let $h\in\cT_k$ with $\|h\|_*=1$. Then there exist $t_m\downarrow0$ and $G_m\in\cO_k(\Theta)$ with
\[
\frac{p_{G_m}-p_*}{t_mp_*}\to h\quad\text{in }L^2(P_*).
\]
Then we have 
\[
\frac{\rho(G_m)}{t_m}
=\left\|\frac{p_{G_m}-p_*}{t_mp_*}\right\|_*\to1.
\]
Consequently the normalized secants satisfy
\[
s_{G_m}=
\frac{t_m}{\rho(G_m)}\,
\frac{p_{G_m}-p_*}{t_mp_*}\to h
\quad\text{in }L^2(P_*).
\]
Thus the unit attainable tangent section $\mathcal{S}_{k}$ lies in the closure of the normalized secant class $\cS$, and therefore this closure is compact in $L^2(P_*)$ due to Lemma~\ref{lem:donsker}. 

We prove that $\mathcal{T}_k$ is closed. Indeed, if $h_m\in\cT_k$ and $h_m\to h$ in $L^2(P_*)$, then from the definition of $\cT_k$m, we can choose positive numbers $t_m<1/m$ and measures $G_m$ such that
\[
\left\|\frac{p_{G_m}-p_*}{t_mp_*}-h_m\right\|_*\le m^{-1}.
\]
The triangle inequality yields convergence of $\frac{p_{G_m}-p_*}{t_mp_*}$ to $h$ in $L^2(P_*)$, consequently $h\in\cT_k$. An implication of this result is that $\mathcal{S}_k$ is compact in $L^2(P_*)$.
% This argument uses $L^2$ throughout; $L^1$ convergence would not suffice for the likelihood geometry.

Finally, we show that $\cT_{k,M}$ is compact in $L^2(P_*)$. For every $a>0$, if $h=h_{(\alpha,\beta,C,\nu)}\in\cT_k$ then $ah=h_{(a\alpha,a\beta,aC,a\nu)}\in\cT_k$: positive scaling preserves mass conservation, matrix ranks and support counts. Hence the bounded section $\cT_{k,M}$ is the continuous image of $[0,M] \times \mathcal{S}_{k} \to \cT_{k}$, $(a,s)\mapsto as$, and as $[0,M] \times \mathcal{S}_{k}$ is a compact set, $\cT_{k,M}$ is a compact set in $L^2(P_*)$ too. 
\end{proof}

% \section{Gaussian quadratic argmax and convergence of the likelihood maximizer}\label{sec:proof3}
% After localization, the likelihood has the familiar Gaussian-quadratic form of nonregular likelihood theory, but the feasible set is nonconvex.  The arguments below therefore combine compactness of the attainable set with a moving-domain argmax theorem rather than an ordinary projection formula.  This follows the general perspective of \citet{Chernoff1954,Shapiro1985,SelfLiang1987,Geyer1994,DacunhaGassiat1997,DacunhaGassiat1999,Drton2009}, while the exact feasible set is the one established in Appendix~\ref{sec:proof2}.

% The structure is a local asymptotic quadratic expansion indexed by a non-Euclidean feasible set; compare the empirical-process likelihood expansions in \citet{vdV1998,vdVW1996}.

\subsection{Likelihood expansion and convergence of the maximizer}
\label{sec:likelihood_expansion_and_convergence}
We now use the tangent-set geometry established above to prove
Theorem~\ref{thm:density}. Lemma~\ref{lem:LAN} provides a uniform
quadratic expansion of the local log-likelihood, and
Lemma~\ref{lem:unique} establishes existence and almost-sure
uniqueness of the limiting Gaussian maximizer. Together with
root-$n$ localization and compactness of bounded tangent
sections, these results yield weak convergence of the
relative-density MLE through an argmax argument.

\begin{lemma}[Uniform local quadratic likelihood expansion]\label{lem:LAN}
For every fixed $M<\infty$,
\begin{equation}\label{eq:uniformLAN}
\sup_{\sqrt n\rho(G)\le M}
\left|\ell_n(G)-\ell_n(G_*)-\mathbb G_n h_{n,G}+\frac12\norm{h_{n,G}}_*^2\right|\pto0,
\end{equation}
where
\[
h_{n,G}=\sqrt n\frac{p_G-p_*}{p_*}.
\]
\end{lemma}

\begin{proof}

\textit{Step 1: }
We precise the neighbors of $G_*$ in $\mathcal{O}_k(\Theta)$ containing measure $G$ that we work through out the proof. Consider the neighborhood $\cU$ in Lemma~\ref{lem:inverse}. By identifiability, the unique zero point of $\rho$ is $G_*$. By compactness of $\mathcal O_k(\Theta)$ and continuity of $G\mapsto\rho(G)$, there is $\delta>0$ such that $\rho(G)<\delta$ implies $G\in\cU$. For $M$ fixed, consider $\mathcal H_{n,M}^{\rm att}
=\left\{h_{n,G}:\sqrt n\rho(G)\le M\right\}$, then for $n\geq (M/\delta)^2$, we have $\mathcal H_{n,M}^{\rm att} \subset \mathcal{U}$. 
In addition, every function in the $H_{n,M}^{\rm att}$ can be written as
\[
h_{n,G}=a_Gs_G,\qquad 0\le a_G\le M,\qquad s_G\in\cS.
\]
where the case $\rho(G)=0$ gives $G=G_*$ by identifiability and hence $h_{n,G}=0$.

\textit{Step 2: } We give an estimation for the process $\max_{r\leq n}|h_{n,G}(X_r)|$ so that we use later in the proof. From Lemma~\ref{lem:donsker}, $|s_G|\le F_0$, hence $|h_{n,G}(x)|\le MF_0(x)$ uniformly over the local set. In particular, we have the bound $$P_*h_{n,G}^2\le M^2\|F_0\|^2_*.$$
For $\eps>0$, the union bound and Markov's inequality of order $2+\eta$ give
\[
\begin{aligned}
\PP\left\{\max_{1\le r\le n}F_0(X_r)>\eps\sqrt n\right\}
&\le n\,\PP_{X\sim P_*}\{F_0(X)>\eps\sqrt n\}\\
&\le \frac{n\,P_*F_0^{2+\eta}} {\eps^{2+\eta}n^{1+\eta/2}}
=\frac{P_*F_0^{2+\eta}}{\eps^{2+\eta}n^{\eta/2}}
\longrightarrow0.
\end{aligned}
\]
Consequently,
\begin{equation}\label{eq:uniformObservedPerturbationDetailed}
\sup_{\sqrt n\rho(G)\le M}
\max_{r\le n}\left|\frac{h_{n,G}(X_r)}{\sqrt n}\right|
\le \frac{M}{\sqrt n}\max_{r\le n}F_0(X_r)
=o_P(1).
\end{equation}

\textit{Step 3}: 
We expand the likelihood difference $\ell_n(G)-\ell_n(G_*)$ using Taylor expansion and estimate the first, second, and remainder components. Consider the event $$\Xi_n = \left\{\sup_{\sqrt n\rho(G)\le M}
\max_{r\le n}\left|\frac{h_{n,G}(X_r)}{\sqrt n}\right| < \frac{1}{2}\right\},$$
whose probability tends to $1$ as $n\to \infty$ by \eqref{eq:uniformObservedPerturbationDetailed}. In this event, by writing $\frac{p_G(X_r)}{p_*(X_r)}
=1+\frac{h_{n,G}(X_r)}{\sqrt n}$, and applying the following estimation
for $|x| < 1/2$
\[
\log(1+x)=x-\frac12x^2+R(x),
\qquad |R(x)|\le C|x|^3,
\]
we have from estimation \eqref{eq:uniformObservedPerturbationDetailed} that in the set whose probability tending to one when $n\to \infty$,
\begin{align}
\ell_n(G)-\ell_n(G_*)
&=\frac1{\sqrt n}\sum_{r=1}^n h_{n,G}(X_r)
-\frac1{2n}\sum_{r=1}^n h_{n,G}(X_r)^2
+R_{n,G}\nonumber \\
&= \underbrace{\sqrt{n}P_nh_{n,G}}_{\text{first-order}} - \underbrace{\frac{1}{2}P_nh^2_{n,G}}_{\text{second-order}} + \underbrace{R_{n,G}}_{\text{remainder}}.
\label{eq:LANdecompDetailed}
\end{align}
where the remainder $R_{n,G}$ can be estimated as  $$|R_{n,G}|
\le C\sum_{r=1}^n
\left|\frac{h_{n,G}(X_r)}{\sqrt n}\right|^3.$$ 
The next step is to give a representation or estimation for each first, second-order components and remainder $R_{n,G}$. 

\textit{Step 3.1: Representation for first-order component ${\sqrt n}P_nh_{n,G}$.} For every attainable local function, we have  $P_*h_{n,G}
=\sqrt n\int(p_G-p_*)\,d\mu=0$.
Therefore, the first-order term in equation \eqref{eq:LANdecompDetailed} can be written as 
\begin{equation}
\label{eqn:dung_taylor_expansion_first_order}
    {\sqrt n}P_nh_{n,G}
= \sqrt{n}(P_nh_{n,G}-P_*h_{n,G}) = \mathbb G_nh_{n,G}. 
\end{equation}

\textit{Step 3.2: Estimation for second-order component $\frac{1}{2}P_nh^2_{n,G}$}.
Taking into account that $h_{n,G}=a_Gs_G$ with $a_G\le M$, using the Glivenko-Cantelli of the class of function $\{s^2:s\in \mathcal{S}\}$ in Lemma~\ref{lem:donsker} , we have 
\begin{equation}
\label{eqn:dung_sup_P_n_minus_P_*_over_small_set}
    \sup_{\sqrt n\rho(G)\le M}
|(P_n-P_*)h_{n,G}^2|
\le M^2\sup_{s\in\cS}|(P_n-P_*)s^2|
=o_P(1).
\end{equation}

\textit{Step 3.3: Estimation for remainder component $R_{n,G}$}.
Since $P_*h_{n,G}^2\le M^2\|F_0\|^2_*$, estimation \eqref{eqn:dung_sup_P_n_minus_P_*_over_small_set} yields 
\begin{equation}
\label{eqn:dung_sup_P_n_over_small_set}
    \sup_{\sqrt n\rho(G)\le M}|P_nh_{n,G}^2| = O_P(1).
\end{equation}
Using $\sum|x_r|^3\le(\max|x_r|)\sum x_r^2$ with $x_r=h_{n,G}(X_r)/\sqrt n$,
\[
|R_{n,G}|
\le C
\left\{\max_{r\le n}\frac{|h_{n,G}(X_r)|}{\sqrt n}\right\}
P_nh_{n,G}^2.
\]
By plugging in the estimation in \eqref{eq:uniformObservedPerturbationDetailed} and \eqref{eqn:dung_sup_P_n_over_small_set} to the estimation above, we have 
\begin{equation}
\label{eqn:dung_taylor_expansion_remainder_final_estimation}
    \sup_{\sqrt n\rho(G)\le M}|R_{n,G}|=o_P(1).
\end{equation}

\textit{Step 3.4: Putting everything together.} Substituting the estimations in \eqref{eqn:dung_taylor_expansion_first_order}, \eqref{eqn:dung_sup_P_n_minus_P_*_over_small_set}, and \eqref{eqn:dung_taylor_expansion_remainder_final_estimation} into \eqref{eq:LANdecompDetailed}, we obtain
\[
\sup_{\sqrt n\rho(G)\le M}
\left|
\ell_n(G)-\ell_n(G_*)
-\mathbb G_nh_{n,G}
+\frac12\|h_{n,G}\|_*^2
\right|\pto0,
\]
which is \eqref{eq:uniformLAN}. This completes our proof. 
\end{proof}

\begin{lemma}[Existence and almost-sure uniqueness of the Gaussian argmax]\label{lem:unique}
The maximization problem~\eqref{eq:Hdagger} admits an almost surely unique
maximizer. In other words, there exists an almost surely unique
$H^\dagger\in \mathcal{T}_k$ such that
\[
Q(H^\dagger)=\sup_{h\in \mathcal T_k}Q(h).
\]
Moreover, for every fixed $M>0$, the restriction of $Q$ to
$\mathcal T_{k,M}:=\mathcal T_k\cap\{h:\|h\|_*\le M\}$ also admits an almost surely
unique maximizer $H_M^\dagger$.
\end{lemma}

\begin{proof}
\emph{Step 1: Existence.} We prove the almost sure existence of maximizer of $Q$ over $\cT_k$. 
By Lemma~\ref{lem:donsker}, the class $\mathcal S$ is
$P_*$-Donsker. Therefore, there exists a version of Gaussian process $\mathbb{W}$ in $\mathcal{S}$ with bounded and uniformly continuous sample path. This version extends continuously to the
$L^2(P_*)$-closure $\overline{\mathcal S}$. Since $\{h\in \mathcal T_k:\|h\|_*=1\}\subseteq\overline{\mathcal S}$, we obtain
\[
B:=\sup_{\substack{h\in \mathcal T_k\\\|h\|_*=1}}|\mathbb{W}(h)|<\infty
\qquad\text{almost surely}.
\]
% Using Lemma~\ref{lem:tangentcompact} and Example 1.5.10 in \citet{vdVW1996}, there exists a version of Gaussian process $\mathbb{W}$ with bounded and uniformly continuous sample path. We define the following positive random variable  
% We prove the existence almost surely of maximizer for the problem \eqref{eq:Hdagger}.  By Lemma~\ref{lem:donsker}, the compact closure $K_1$ of the normalized secants carries a tight, separable, uniformly $L^2(P_*)$-continuous version of the limiting isonormal process. Restrict to a single probability-one event on which that version is continuous and bounded on $K_1$. Positive homogeneity then extends it consistently to each bounded tangent section. In particular,
% \[
% B:=\sup\{|\mathbb W(h)|:h\in\cT_k,\ \|h\|_*=1\}<\infty
% \quad\text{almost surely},
% \]
Then $|\mathbb W(h)|\le B\|h\|_*$ for every $h\in \cT_k$. Since
$Q(h)=\mathbb W(h)-\|h\|_*^2/2\le B\|h\|_*-\|h\|_*^2/2 \to -\infty$ as $\|h\|_* \to \infty$, the supremum of $Q(h)$ is attained.

% The maximum and a measurable argmax can be defined using the countable dense subset of each compact section and the measurable maximum theorem.

\emph{Step 2: Uniqueness.} For any two almost surely measurable maximizer $h_1$ and $h_2$ of $Q$ over $\cT_k$, we show that $h_1=h_2$ almost surely.  Let $(e_r)_{r\ge1}$ be an orthonormal basis of $\cH=\overline{\operatorname{span}(\cT_k)}$, put $Z_r=\mathbb W(e_r)$, and set
$\mathbb W_{-r}(h)=\mathbb W(h)-Z_r\langle h,e_r\rangle_*$. For fixed $r$, $Z_r\sim N(0,1)$ is independent of the separable remainder process $\mathbb W_{-r}$. Conditional on this remainder, the function
\[
V_r(z)=\sup_{h\in\cT_k}
\left\{z\langle h,e_r\rangle_*+\mathbb W_{-r}(h)-\frac12\|h\|_*^2\right\}
\]

Joint Gaussianity and orthogonality imply that $Z_r\sim N(0,1)$
is independent of the entire process $\mathbb W_{-r}$.
Fix a realization on which
$|\mathbb W(h)|\le B\|h\|_*$ for every $h\in\cT_k$.
For $z\in\RR$, define
\[
V_r(z)
=
\sup_{h\in\cT_k}
\left\{
z\langle h,e_r\rangle_*
+\mathbb W_{-r}(h)
-\frac12\|h\|_*^2
\right\}.
\]
For each fixed $h$, the expression inside the supremum is
affine in $z$. Hence $V_r$ is convex. Moreover,
\[
\begin{aligned}
z\langle h,e_r\rangle_* +\mathbb W_{-r}(h) -\frac12\|h\|_*^2
&\le |z-Z_r|\,|\langle h,e_r\rangle_*| +B\|h\|_*-\frac12\|h\|_*^2\\
&\le \bigl(|z-Z_r|+B\bigr)\|h\|_* -\frac12\|h\|_*^2
\rightarrow -\infty
\end{aligned}
\]
as $\|h\|_*\to\infty$. Since bounded tangent sections are
compact and the criterion is continuous on these sections,
$V_r(z)$ is finite and its supremum is attained.

A finite convex function on $\RR$ is differentiable outside
a Lebesgue-null set. At a differentiability point $z$,
suppose that $h_1$ and $h_2$ both attain the supremum.
Their slopes $a_j=\langle h_j,e_r\rangle_*$, for $j=1,2$,
belong to the subdifferential $\partial V_r(z)$, because
the corresponding affine functions support $V_r$ at $z$.
Differentiability implies
$\partial V_r(z)=\{V_r'(z)\}$, so $a_1=a_2$.
Thus, at every differentiability point, all maximizers
have the same $e_r$ coordinate.

The exceptional nondifferentiability set may depend on
the realization of $\mathbb W_{-r}$, but has Lebesgue
measure zero for almost every such realization.
Conditional on $\mathbb W_{-r}$, the independent variable
$Z_r$ has a standard normal density. Consequently,
\[
\PP\left\{
Z_r\text{ belongs to the exceptional set}
\,\middle|\,\mathbb W_{-r}
\right\}
=0
\qquad\text{almost surely}.
\]
Taking expectations shows that, with probability one,
all maximizers $h$ of the original random criterion have
the same coordinate $\langle h,e_r\rangle_*$. Taking the countable intersection of these probability-one
events over the orthonormal basis, we have $\langle h_1-h_2,e_r\rangle_*=0$
for every $r$. Therefore $h_1=h_2$ in $\cH$.

\textit{Step 4: Restricted problem}. Lastly, we prove the existence and uniqueness for the solution of $Q$ in $\mathcal{T}_{k,M}$.  Existence follows from compactness of $\mathcal{T}_{k,M}$ and almost sure continuity of $Q$ on this set.

For uniqueness, we use the same argument as in Step 2. Indeed, we use the same orthonormal basis $(e_r)_{r\geq 1}$
and decomposition
\[
\mathbb{W}(h)=Z_r\langle h,e_r\rangle_*+\mathbb{W}_{-r}(h)
\]
as above. Conditional on $\mathbb{W}_{-r}$, define
\[
V_{r,M}(z):=
\sup_{h\in \mathcal{T}_{k,M}}
\left\{
z\langle h,e_r\rangle_*+\mathbb{W}_{-r}(h)
-\frac12\|h\|_*^2
\right\}.
\]
This function is finite and convex on $\mathbb R$, and
the supremum is attained. At each differentiability point
of $V_{r,M}$, all maximizing points have the same
$e_r$ coordinate, because their coordinates belong to
the subdifferential of $V_{r,M}$.
Since $Z_r$ is independent of $\mathbb{W}_{-r}$ and has a density,
it almost surely avoids the Lebesgue-null set of
nondifferentiability points. Taking a countable intersection
over $r$ shows that all maximizers have identical
coordinates in the orthonormal basis, and hence coincide.
% Hence the maximizer exists and is almost surely unique.
% is finite, convex and measurable in $z$. Finiteness follows from the same quadratic coercivity bound, now using the almost surely bounded remainder on the unit tangent section. A finite convex function on $\RR$ is differentiable outside a Lebesgue-null set. At a differentiability point, all maximizing affine functions have the same slope, hence all maximizers share the same $e_r$ coordinate. This exceptional set is measurable, since $V_r$ is a separable supremum of measurable affine functions. Because the independent $Z_r$ has a density, the realized $Z_r$ lies outside that set with conditional probability one.

% Take the countable intersection of these probability-one events over $r$. Any two maximizers then have equal coordinates along every $e_r$, and therefore agree as elements of $\cH$. No convexity of $\cT_k$ is used.
\end{proof}

% \begin{proposition}[Convergence of the likelihood maximizer]\label{prop:argmax}
% Under the standing assumptions,
% \[
% \widehat h_n=\sqrt n\frac{p_{\widehat G_n}-p_*}{p_*}\weakto H^\dagger
% \quad\text{in }L^2(P_*).
% \]
% \end{proposition}

\begin{proof}[Proof of Theorem~\ref{thm:density}]
We have established the convergence rate \eqref{eq:rhorate} of $\rho(\widehat G_n)$ in Lemma~\ref{lem:rootn}. For the weak limit of scaled relative density discrepancy $\widehat h_n$, we divide our proof into following steps. 

\emph{Step 1.} For $M>0$, define
\[
\cH_{n,M}
=
\left\{
h_{n,G}:=\sqrt n\,\frac{p_G-p_*}{p_*}:
G\in\cO_k(\Theta),\ \sqrt n\,\rho(G)\le M
\right\}.
\]
We can verify that $\cH_{n,M}$ is a compact subset of $L^2(P_*)$. Using Lemma~\ref{lem:Haussdorf}, we show that $\cH_{n,M}$ converges in Hausdorff distance
to $\cT_{k,M}$ in $L^2(P_*)$, defined in Lemma~\ref{lem:tangentcompact}.
Let $\overline{\cS}$ denote the closure of $\cS$ in
$L^2(P_*)$. By Lemma~\ref{lem:donsker}(ii), $\cS$ is
totally bounded, so $\overline{\cS}$ is compact.
Define
\[
K_M
=
\{as:0\le a\le M,\ s\in\overline{\cS}\}.
\]
The map $(a,s)\mapsto as$ is continuous, so $K_M$ is
compact. In addition, for $n$ sufficiently large, both
$\cH_{n,M}$ and $\cT_{k,M}$ are contained in $K_M$.

Suppose that $h_{n_j}\in\cH_{n_j,M}$ and
$h_{n_j}\to h$ in $L^2(P_*)$. There exist
$G_{n_j}\in\cO_k(\Theta)$ such that
\[
h_{n_j}
=
\frac{p_{G_{n_j}}-p_*}{t_{n_j}p_*},
\qquad t_{n_j}=n_j^{-1/2}.
\]
Proposition~\ref{prop:tangent} gives $h\in\cT_k$.
Norm continuity also gives $\|h\|_*\le M$, hence
$h\in\cT_{k,M}$.

Since $T_k$ is a cone, the subset
$\{h\in T_k:\|h\|_*<M\}$ is dense in $T_{k,M}$:
for every $h\in T_{k,M}$, the points $(1-j^{-1})h$
belong to this subset and converge to $h$. Fix $h\in\cT_{k,M}$ with $\|h\|_*<M$.
The recovery construction in
Proposition~\ref{prop:tangent} provides mixing measures
$G_t$ for all sufficiently small $t>0$ such that ${(p_{G_t}-p_*)}/(tp_*)\to h$ in $L^2(P_*)$. For $n$ sufficiently large, taking $t=t_n=n^{-1/2}$ and setting
\[
h_n
=
\sqrt n\,\frac{p_{G_{t_n}}-p_*}{p_*}
\]
gives $h_n\to h$. Since $\|h\|_*<M$, we have
$h_n\in\cH_{n,M}$ for all sufficiently large $n$. 
Combining all the above results with Lemma~\ref{lem:Haussdorf}, we have $d_{\mathrm H} \left(\cH_{n,M},\cT_{k,M}\right)\rightarrow 0$, where the closures and Hausdorff distance are taken in $L^2(P_*)$.

\textit{Step 2.}
Fix $M>0$. On the common compact index set $K_M$, define
\[
Z_n(h):=\mathbb G_n h-\frac12\|h\|_*^2,
\qquad h\in K_M.
\]
The Donsker property in Lemma~\ref{lem:donsker} extends to $K_M$ by
taking the closure of $\mathcal S$ and then allowing bounded
scalar multiplication. Consequently, $Z_n\rightsquigarrow Q$ in $\ell^\infty(K_M)$, where $Q(h)=W(h)-\frac12\|h\|_*^2$ has almost surely continuous sample paths under the $L^2(P_*)$ metric.

By Lemma~\ref{lem:LAN}, the local likelihood criterion satisfies $Q_n(h)=Z_n(h)+r_n(h)$ for $h\in\mathcal H_{n,M}$, where $\Delta_{n,M}:=
\sup_{h\in\mathcal H_{n,M}}|r_n(h)|=o_P(1)$. Let $\widehat h_{n,M}$ be a measurable maximizer of $Q_n$
over $\mathcal H_{n,M}$. Then
\[
\begin{aligned}
Z_n(\widehat h_{n,M}) \ge Q_n(\widehat h_{n,M})-\Delta_{n,M}
=\sup_{h\in\mathcal H_{n,M}}Q_n(h)-\Delta_{n,M}\ge\sup_{h\in\mathcal H_{n,M}}Z_n(h)-2\Delta_{n,M}.
\end{aligned}
\]
In addition, by Step 1, $d_H(\mathcal H_{n,M},\mathcal{T}_{k,M})\longrightarrow0$,
and both sets lie in $K_M$ for all sufficiently large $n$. The existence and uniqueness argument in Lemma~\ref{lem:unique} applies also to the compact set $T_{k,M}$, so $H_M^\dagger =
\operatorname*{arg\,max}_{h\in T_{k,M}}Q(h)$ exists and is almost surely unique. Theorem 1 of \citep{cox2022generalized} yields
\[
\widehat h_{n,M}\rightsquigarrow H_M^\dagger
\qquad\text{in }L^2(P_*).
\]

\textit{Step 3.} Lastly, we prove that $h_n\rightsquigarrow H^\dagger$. It is possible to choose the restricted maximizer $\widehat h_{n,M}$ to agree
with $\widehat h_n$ whenever $\|\widehat h_n\|_*\le M$ as $\widehat h_n$ is still a maximizer for restrained problem in $\mathcal H_{n,M}$.
Hence
\[
P\{\widehat h_{n,M}\ne\widehat h_n\}
\le P\{\|\widehat h_n\|_*>M\}.
\]
Similarly, almost sure uniqueness of the global maximizer
implies $H_M^\dagger=H^\dagger$ on $\{\|H^\dagger\|_*\le M\}$.

Let $f:L^2(P_*)\to\mathbb R$ be bounded and Lipschitz,
with $\|f\|_\infty\le1$. The preceding observations give
\[
\begin{aligned}
\bigl|E f(\widehat h_n)-E f(H^\dagger)\bigr|
\le
2P\{\|\widehat h_n\|_*>M\} +
\bigl|E f(\widehat h_{n,M})-E f(H_M^\dagger)\bigr| +
2P\{\|H^\dagger\|_*>M\}.
\end{aligned}
\]
For each fixed $M$, the middle term converges to zero as $n\to \infty$
by Step 2. Moreover, Lemma~\ref{lem:rootn} gives
$\|\widehat h_n\|_*=O_P(1)$, and Lemma~\ref{lem:unique} gives
$\|H^\dagger\|_*<\infty$ almost surely. Therefore,
taking $n\to\infty$ and then $M\to\infty$ yields $E f(\widehat h_n)\longrightarrow E f(H^\dagger)$. Since bounded Lipschitz functions determine weak convergence
on $L^2(P_*)$, we conclude that
\[
\widehat h_n\rightsquigarrow H^\dagger
\qquad\text{in }L^2(P_*).
\]
This completes the proof.

\end{proof}

% \begin{proof}[Proof of Theorem~\ref{thm:density}]
% Lemma~\ref{lem:rootn} gives the localization statement \eqref{eq:rhorate}.  Proposition~\ref{prop:argmax} proves the weak convergence \eqref{eq:densitylimit}.  These two conclusions are exactly the assertions of Theorem~\ref{thm:density}.
% \end{proof}

\section{Proofs for Section~\ref{sec:wassersteinlimit}}
\label{app:wassersteinproofs}
% The likelihood limit alone does not determine the Wasserstein limit.  We therefore prove a separate sequential derivative theorem for $W_2^2$.  The proof is constructive: local atoms are collapsed onto their true centers, cross-cluster flows are shown to have only $O(t)$ mass, and the remaining first-order problem becomes a finite transportation linear program.  We use only finite-support optimal-transport arguments, consistent with the general frameworks of \citet{Villani2009,Santambrogio2015,PeyreCuturi2019} and the variational convergence language of \citet{RockafellarWets1998}.

\subsection{Sequential expansion and continuity of the transportation functional}
\label{app:sequential_expansion}
In this subsection, we prove the sequential first-order
expansion in Proposition~\ref{prop:transportmain} and the
continuity of $\Omega$. Lemma~\ref{lem:collapse} separates
the squared Wasserstein cost into local collapse and
residual mass transport, while
Lemma~\ref{lem:transportstability} establishes stability
of the finite transportation problem governing the
residual cost. Together, these results yield the
sequential expansion. We then prove continuity of
$\Omega$ in Lemma~\ref{lem:Omegacont}, providing the
additional ingredient needed to combine this expansion
with the density weak limit.
% \begingroup\color{red}
\begin{lemma}[Local collapse and first-order squared transport]\label{lem:collapse}
Let $t_m\downarrow0$ and $G_m\in\cO_k(\Theta)$ satisfy
\[
\frac{p_{G_m}-p_*}{t_mp_*}\longrightarrow h_\gamma\in\cT_k
\quad\text{in }L^2(P_*).
\]
On a classified subsequence as in proof of Proposition~\ref{prop:tangent}, write $L_i$ for labels approaching $\theta_i^*$ and $L_0$ for the other labels, and set
\[
\Delta_{j,m}=\theta_{j,m}-\theta_i^*\quad(j\in L_i),\qquad
A_m=\sum_{i=1}^p\sum_{j\in L_i}\pi_{j,m}\|\Delta_{j,m}\|^2.
\]
Let $M_{i,m}=\sum_{j\in L_i}\pi_{j,m}$ and define the collapsed measure
\[
\overline G_m=\sum_{i=1}^pM_{i,m}\delta_{\theta_i^*}
+\sum_{j\in L_0}\pi_{j,m}\delta_{\theta_{j,m}}.
\]
Then
\begin{align}
A_m/t_m&\longrightarrow\sum_{i=1}^p\operatorname{tr}(C_i),\label{eq:AtraceDetailed}\\
W_2^2(G_m,G_*)&=A_m+W_2^2(\bar G_m,G_*)+o(t_m).\label{eq:collapseIdentityDetailed}
\end{align}
\end{lemma}

\begin{proof}
The trace limit follows immediately from the local second
moments in the of necessity part, proof of Proposition~\ref{prop:tangent}:
\[
\frac{1}{t_m}\sum_{j\in L_i}
\pi_{j,m}\Delta_{j,m}\Delta_{j,m}^{\top}\to C_i.
\]
Taking traces and summing over $i$ gives \eqref{eq:AtraceDetailed}.
It remains to estimate $W_2^2(G_m,G_*)$ through
$W_2^2(\overline G_m,G_*)$ and $A_m$.

\textit{Step 1. } We prove that $W_2^2(G_m,G_*)=O(t_m)$. Using the same argument as in Step 1, Part 1 in the proof
of Proposition~\ref{prop:tangent}, we have $\rho(G_m)=O(t_m)$ and $G_m$
belongs to the neighborhood $\mathcal U$ defined in
Lemma~\ref{lem:inverse} such that we can apply the Lemma~\ref{lem:w2bound} to achieve
$W_2^2(G_m,G_*)=O(t_m)$. 

\textit{Step 2.} We first relate transportation plans from $G_m$ to $G_*$
to those from $\overline G_m$ to $G_*$, and conversely.
Comparing their costs then allows us to estimate
$W_2^2(G_m,G_*)$ in terms of the local collapse cost
$A_m$ and the residual transportation cost
$W_2^2(\overline G_m,G_*)$ in the next step. Let $\Gamma=(\gamma_{jr})_{j,r}$ be a coupling from
$G_m$ to $G_*$. We can check that the transportation
plan $\Gamma$ induces a transportation plan
$\overline\Gamma=(\overline\gamma_{ir})_{i,r}$
between $\overline G_m$ and $G$ as
\begin{equation*}
\overline\gamma_{i,r}
=
\sum_{j\in L_i}\gamma_{jr}.
\end{equation*}
Conversely, let
$\overline\Gamma=(\overline\gamma_{ir})_{i,r}$
be a coupling between $\overline G_m$ and $G$.
We construct a transportation plan
$\Gamma=(\gamma_{jr})_{j,r}$ between $G_m$ to $G_*$
satisfying relation
$\overline\gamma_{i,r}=\sum_{j\in L_i}\gamma_{jr}$.
In fact, this transportation plan can be constructed as
\begin{equation*}
\gamma_{jr}
=
\begin{cases}
\overline\gamma_{jr},
& \text{if }j\in L_0,\\[0.5em]
\displaystyle
\frac{\pi_{j,m}}{\sum_{k\in L_i}\pi_{k,m}}
\overline\gamma_{jr},
& \text{if }j\in L_i\text{ and }i\ne0,
\end{cases}
\end{equation*}
and we can check the condition (75) for this
transportation plan. Henceforth, unless otherwise
stated, given a transportation plan $\Gamma$ from
$G_m$ to $G_*$, we write $\overline\Gamma$ for the
transportation plan from $\overline G_m$ to $G$
constructed via (75). Conversely, given a transportation
plan $\overline\Gamma$ from $\overline G_m$ to $G$,
we write $\Gamma$ for the transportation plan from
$G_m$ to $G_*$ constructed via (76).
We also write $\Pi_\Gamma(G_m,G_*)$ and
$\Pi_{\overline\Gamma}(\overline G_m,G_*)$ be the
total squared costs corresponding to transportation
plan $\Gamma$ and $\overline\Gamma$ as
\[
\Pi_\Gamma(G_m,G_*)
=
\sum_{j,r}\gamma_{jr}
\|\theta_{j,m}-\theta_r^*\|^2,
\qquad
\Pi_{\overline\Gamma}(\overline G_m,G_*)
=
\sum_{i,r}\overline\gamma_{ir}
\|\theta_i^*-\theta_r^*\|^2
+
\sum_{i\in L_0}\sum_r\overline\gamma_{ir}
\|\theta_{i,m}-\theta_r^*\|^2.
\]

For $j\in L_i$, we consider the expansion
\begin{equation*}
\|\theta_{j,m}-\theta_r^*\|^2
=
\|\Delta_{j,m}\|^2
+
\|\theta_i^*-\theta_r^*\|^2
+
2\Delta_{j,m}^{\top}(\theta_i^*-\theta_r^*).
\end{equation*}
For a transportation plan $\Gamma=(\gamma_{jr})_{j,r}$,
by multiplying with weight $\gamma_{jr}$, we have
\[
\begin{aligned}
\Pi_\Gamma(G_m,G_*)
&=
\sum_{i=1}^{k}\sum_{j\in L_i}\sum_r
\gamma_{jr}
\left(
\|\Delta_{j,m}\|^2
+\|\theta_i^*-\theta_r^*\|^2
+2\Delta_{j,m}^{\top}(\theta_i^*-\theta_r^*)
\right)
+
\sum_{j\in L_0}\sum_r
\gamma_{jr}\|\theta_{j,m}-\theta_r^*\|^2\\
&=
\underbrace{
\left(
\sum_{i=1}^{k}\sum_{j\in L_i}\sum_r
\gamma_{jr}\|\Delta_{j,m}\|^2
\right)
}_{=A_m}
+
\underbrace{
\left(
\sum_{i=1}^{k}\sum_{j\in L_i}\sum_r
\gamma_{jr}\|\theta_i^*-\theta_r^*\|^2
+
\sum_{j\in L_0}\sum_r
\gamma_{jr}\|\theta_{j,m}-\theta_r^*\|^2
\right)
}_{\Pi_{\overline\Gamma}(\overline G_m,G_*)}\\
&\quad
+
2\sum_{j\in L_i}\sum_r
\gamma_{jr}\Delta_{j,m}^{\top}
(\theta_i^*-\theta_r^*).
\end{aligned}
\]
Therefore, we have
\begin{equation}
\Pi_\Gamma(G_m,G_*)
=
A_m+\Pi_{\overline\Gamma}(\overline G_m,G_*)
+
2\underbrace{
\sum_{j\in L_i}\sum_{r\ne i}
\gamma_{jr}\Delta_{j,m}^{\top}
(\theta_i^*-\theta_r^*)
}_{=:I_m}.
\label{eqn:dung_main_estimation_wasserstein_based_on_grouping}
\end{equation}

\textit{Step 3. } We prove the estimation \eqref{eq:collapseIdentityDetailed}.
Let $\Gamma=(\gamma_{jr})_{j,r}$ be the optimal coupling
between $G_m$ and $G_*$. For
$s_*:=\min_{i\ne r}\|\theta_i^*-\theta_r^*\|>0$
be the minimal separation among true center,
since the local displacements vanish, for all large $m$, we have
$\max_i\max_{j\in L_i}\|\Delta_{j,m}\|\le s_*/4$.
Hence a local atom from cluster $i$ transported to
a different true atom $\theta_r^*$ travels at least
$3s_*/4$. In the optimal plan $\Gamma$, we therefore
achieve that $\gamma_{jr}=O(t_m)$ for $j\in L_i$
and $i\ne r$. Consequently, we have the following
estimation for $I_m$ as
\begin{equation}
|I_m|
\le
\sum_{j\in L_i}\sum_{r\ne i}
\underbrace{\gamma_{jr}}_{O(t_m)}
\underbrace{\|\Delta_{j,m}\|}_{o(1)}
\cdot
\underbrace{\|\theta_i^*-\theta_r^*\|}_{\le\operatorname{Diam}(\Theta)}
=
o(t_m).
\label{eqn:estimation_for_cross_term}
\end{equation}
Thus, by plugging in the estimation \eqref{eqn:estimation_for_cross_term} for $I_m$
and the inequality
$\Pi_{\overline\Gamma}(\overline G_m,G_*)
\ge W_2^2(\overline G_m,G_*)$
to \eqref{eqn:dung_main_estimation_wasserstein_based_on_grouping}, we achieve
\begin{equation}
W_2^2(G_m,G_*)
\ge
A_m+W_2^2(\overline G_m,G_*)+o(t_m).
\label{eqn:lower_bound_of_W_2^2(G_m,G_*)}
\end{equation}

For the inverse direction, let
$\overline\Gamma=(\overline\gamma_{ir})_{i,r}$
be the optimal transportation plan from
$\overline G_m$ to $G_*$, which corresponds to the
transportation plan $\Gamma$ from $G_m$ to $G_*$.
As $W_2^2(G_m,G_*)=O(t_m)$, from the estimation \eqref{eqn:lower_bound_of_W_2^2(G_m,G_*)}
and the fact that $A_m\ge0$, we have
$W_2^2(\overline G_m,G_*)=O(t_m)$.
Therefore, the off-diagonal mass
$\overline\gamma_{ir}$ transported between distinct
true atoms $i,r$ is of order $O(t_m)$, since every such unit
of mass costs at least $s_*^2$.
As a consequence, we also achieve the estimation
$|I_m|=o(t_m)$ in \eqref{eqn:estimation_for_cross_term}. Moreover, $\Pi_\Gamma(G_m,G_*)\ge W_2^2(G_m,G_*)$,
we plug in these estimations above to \eqref{eqn:dung_main_estimation_wasserstein_based_on_grouping} and achieve
\[
W_2^2(G_m,G_*)
\le
A_m+W_2^2(\overline G_m,G_*)+o(t_m).
\]
This completes our proof.
\end{proof}

\begin{lemma}[Stability of a finite transportation problem]\label{lem:transportstability}
Fix distinct locations $\theta_1^*,\ldots,\theta_p^*$ and a nonnegative integer $q$. Suppose $\alpha_m\to\alpha$ in $\mathbb R^p$, $c_m\to c$ in $[0,\infty)^q$, and $\eta_{\ell,m}\to\eta_\ell$ in $\Theta$ such that $\eta_\ell \neq \theta_i^*$ for $1\leq i\leq p$ and 
\[\sum_i\alpha_{i,m}+\sum_{\ell=1}^q c_{\ell,m}=0.\]
For each $m$ define $V_m$ as the minimum of
\[
\Psi_m(x,y) = \sum_{i\ne j}\|\theta_i^*-\theta_j^*\|^2x_{ij}
+\sum_{\ell=1}^q\sum_{i=1}^p\|\eta_{\ell,m}-\theta_i^*\|^2y_{\ell i}
\]
over $x_{ij}$ and $y_{\ell i}$ belong to the constraint set 
\begin{equation}\label{eq:stableLPbalance}
\mathcal{P}_m = \left\{x,y \geq 0,\quad \sum_{i}y_{\ell i}=c_{\ell,m}, \quad \sum_{j\neq i}(x_{ij}-x_{ji}) = \alpha_{i,m} + \sum_{\ell}y_{\ell i}\right\}
\end{equation}
Let $V$ denote the analogous minimum at the limiting masses and locations:
\begin{equation*}
    V = \min_{(x,y) \in \mathcal{P}}\Psi(x,y):=  \sum_{i\ne j}\|\theta_i^*-\theta_j^*\|^2x_{ij}
+\sum_{\ell=1}^q\sum_{i=1}^p\|\eta_{\ell}-\theta_i^*\|^2y_{\ell i}
\end{equation*}
where \begin{equation*}
    \mathcal{P} = \left\{x,y \geq 0,\quad \sum_{i}y_{\ell i}=c_{\ell}, \quad \sum_{j\neq i}(x_{ij}-x_{ji}) = \alpha_{i} + \sum_{\ell}y_{\ell i}\right\}
\end{equation*}
Then $V_m\to V$. This conclusion allows $c_\ell=0$ and arbitrary changes in the active edges of the minimizing flows.
\end{lemma}

\begin{proof}
\textit{Step 1.} First, we show that $\sup_m V_m<\infty$ and there exists a constant $M>0$ such that for $m$ sufficiently large, there exists an optimal solution $(x^m,y^m)$ of $\Psi_m$ over $\mathcal{P}_m$ such that $\|x^m\|_{\infty},\|y^m\|_{\infty}\leq M$. In fact, as $\theta_1^*,\ldots,\theta_p^*$ are distinct and $\eta_{\ell} \neq \theta_i^*$ for all $1\leq i\leq p$, there exists a positive constant $s_*$ such that for $m$ sufficiently large, 
\begin{equation*}
    \|\theta_i^* - \theta_j^*\|\geq s_*, \quad \|\eta_{\ell,m} - \theta_i^*\| \geq s_*. 
\end{equation*}
Therefore, for $m$ sufficiently large, $\Psi_m(x,y) \geq s_*^2(\sum_{i\neq j}x_{ij} + \sum_{\ell = 1}^q\sum_{i=1}^py_{\ell i})$. Meanwhile, we can see that as $c_{\ell,m} \to c$ and $\alpha_m\to \alpha$ and bounded, there exists a constant $C$ such that $\mathcal{P}_m \cap \{\|x\|_{\infty},\|y\|_{\infty} \leq C\}$ and $\mathcal{P} \cap \{\|x\|_{\infty},\|y\|_{\infty} \leq C\}$ are non-empty. Moreover, $\Psi_m \leq Cp(p+q-1)\cdot\mathrm{Diam}(\Theta)$ in $\mathcal{P}_m\cap \{\|x\|_{\infty},\|y\|_{\infty} \leq C\}$ and $\Psi\leq Cp(p+q-1)\cdot\mathrm{Diam}^2(\Theta)$ in $\mathcal{P}\cap \{\|x\|_{\infty},\|y\|_{\infty} \leq C\}$. Combining these above observations, we have $\sup_{m}V_m \leq Cp(p+q-1)\cdot\mathrm{Diam}^2(\Theta) < \infty$ and there exists $M > 0$ such that some optimal solution $(x^m,y^m)$ of $\Psi_m$ over $\mathcal{P}_m$ satisfying $\|x^m\|_{\infty},\|y^m\|_{\infty}\leq M$.

\textit{Step 2. } We prove that $\lim \inf_{m}V_m \geq V$. Consider a sequence $\{m_j\}_j$ such that $V_{m_j} \to \lim \inf_{m}V_m$. By extracting a subsequence, we can suppose that the optimal solution $(x^{m_j},y^{m_j})$ constructed in Step 1 converges to $(x^{\star},y^{\star}) \in \mathcal{P}$. Therefore, 
\begin{equation*}
    \lim \inf_{m}V_m = \lim_{j\to \infty} V_{m_j} = \lim_{j\to \infty} \Psi_{m_j}(x^{m_j},y^{m_j})\geq \lim_{j\to \infty} \Psi_{m_j}(x^{\star},y^{\star}) = \Psi(x^{\star},y^{\star}) \geq V. 
\end{equation*}

\textit{Step 3.} We prove that $\limsup_m V_m\le V$
Fix an optimal solution $(x,y)\in\mathcal P$, we construct a sequence $(x^{(m)},y^{(m)}) \in \mathcal{P}_m$ based on $(x,y)$.  For every $\ell=1,\ldots,q$ and $i=1,\ldots,p$, define
\[
y_{\ell i}^{(m)}
=
\begin{cases}
\displaystyle
\frac{c_{\ell,m}}{c_\ell}y_{\ell i},
& c_\ell>0,\\[0.5em]
c_{\ell,m}\mathbf 1_{\{i=1\}},
& c_\ell=0.
\end{cases}
\]
Then $y^{(m)}\ge0$ and
$\sum_i y_{\ell i}^{(m)}=c_{\ell,m}$.
Moreover, $y^{(m)}\to y$: when $c_\ell=0$,
nonnegativity and $\sum_i y_{\ell i}=0$ imply
$y_{\ell i}=0$ for every $i$. For $x^{(m)}$, we define the balance discrepancies
\[
\delta_{i,m}
=
\alpha_{i,m}+\sum_{\ell=1}^q y_{\ell i}^{(m)}
-\sum_{j\ne i}(x_{ij}-x_{ji}),
\qquad i=1,\ldots,p.
\]
Since $(x,y)\in\mathcal P$ and $y^{(m)}\to y$,
we have $\delta_{i,m}\to0$ for every $i$.
Furthermore, mass conservation gives
\[
\sum_{i=1}^p\delta_{i,m}
=
\sum_{i=1}^p\alpha_{i,m}
+\sum_{\ell=1}^q c_{\ell,m}
=0.
\]
Write $\delta_{i,m}^+=\max\{\delta_{i,m},0\}$ and
$\delta_{i,m}^-=\max\{-\delta_{i,m},0\}$, and put $D_m
:=
\sum_{i=1}^p\delta_{i,m}^+
=
\sum_{i=1}^p\delta_{i,m}^-
=
\frac12\sum_{i=1}^p|\delta_{i,m}|$.
For $i\ne j$, define the correcting flows explicitly by
\[
z_{ij}^{(m)}
=
\begin{cases}
\displaystyle
\frac{\delta_{i,m}^+\delta_{j,m}^-}{D_m},
& D_m>0,\\[0.5em]
0,
& D_m=0.
\end{cases}
\]
These flows are nonnegative. If $D_m>0$, the identity
$\delta_{i,m}^+\delta_{i,m}^-=0$ implies
$\sum_{j\ne i}z_{ij}^{(m)}=\delta_{i,m}^+$ and
$\sum_{j\ne i}z_{ji}^{(m)}=\delta_{i,m}^-$. The same identities hold when $D_m=0$, because all $\delta_{i,m}$ then vanish. Consequently,
$\sum_{j\ne i}
\bigl(z_{ij}^{(m)}-z_{ji}^{(m)}\bigr)
=
\delta_{i,m}$. Define $x_{ij}^{(m)}:=x_{ij}+z_{ij}^{(m)}$.
Then
\[
\sum_{j\ne i}
\bigl(x_{ij}^{(m)}-x_{ji}^{(m)}\bigr)
=
\alpha_{i,m}+\sum_{\ell=1}^q y_{\ell i}^{(m)},
\]
so $(x^{(m)},y^{(m)})\in\mathcal P_m$. Finally, $\sum_{i\ne j}z_{ij}^{(m)}=D_m\to 0$, and hence $x^{(m)}\to x$.
Together with $y^{(m)}\to y$ and
$\eta_{\ell,m}\to\eta_\ell$, this gives
\[
\Psi_m(x^{(m)},y^{(m)})
\longrightarrow
\Psi(x,y)=V.
\]
Since $V_m\le\Psi_m(x^{(m)},y^{(m)})$, we conclude that
$\limsup_m V_m\le V$.
Combined with Step 2, this proves $V_m\to V$.

\end{proof}

% \endgroup
\begin{proof}[Proof of Proposition~\ref{prop:transportmain}]

% \begin{proposition}[Sequential directional derivative of $W_2^2$]\label{prop:transport}
% Let $t_m\downarrow0$ and $G_m\in\cO_k(\Theta)$ satisfy
% \[
% \left\|\frac{p_{G_m}-p_*}{t_mp_*}-h_\gamma\right\|_*\to0
% \]
% for some $h_\gamma\in\cT_k$. Then
% \begin{equation}\label{eq:directionalW2}
% \frac{W_2^2(G_m,G_*)}{t_m}\to \Omega(h_\gamma).
% \end{equation}
% \end{proposition}

We show that every subsequence $\{G_{m_i}\}_i$ of $\{G_m\}_m$ admits a further subsequence along which the convergence in \eqref{eq:directionalmain} holds. For notational simplicity, we relabel the subsequence $\{G_{m_i}\}_i$ as $\{G_m\}_m$. 

Using Proposition~\ref{prop:tangent}, we extract a further subsequence of $G_m$ on which all fitted labels have limiting locations and are classified into local sets $L_i$ and a nonlocal set $L_0$. For notational convenience, we also relabel this subsequence as $G_m$. Along this classified subsequence the scaled local moments converge to the canonical coefficients $(\alpha,\beta,C)$ and the positive scaled nonlocal masses converge to the atoms of $\nu$, which contributes to the tangent configuration $\overline{\gamma} = (\alpha,\beta,C,\nu)$. By uniqueness of the canonical tangent representation, every classified subsequence gives the same $\overline{\gamma} = \gamma = (\alpha,\beta,C,\nu)$.

In this proof, we use the same notation as in Lemma~\ref{lem:collapse}. Using this lemma, we reduce the transportation problem $W_2^2(G_m,G_*)$ to the transportation problem $W_2^2(\overline{G}_m,G_*)$, which we study in the next step. In addition, consider the constant $\mathcal{C}_0$ such that $$W_2^2(G_m,G_*)\leq \mathcal{C}_0t_m, \quad W_2^2(\overline{G}_m,G_*) \leq \mathcal{C}_0t_m.$$

\textit{Step 1: Solution of the transportation problem between $\overline{G}_m$ and $G_*$.}  % We prove that for $m$ sufficiently large, the transportation problem $W_2^2(\overline{G}_m,G_*)$ is equivalent to the quadratic optimization 
% $$\Psi_m(x,y) = \sum_{i\ne j}d_{ij}^2x_{ij}
% +\sum_{\ell,i}\|\eta_{\ell,m}-\theta_i^*\|^2y_{\ell i},$$
% over the set of constraint
% \begin{equation}\label{eq:finiteYrowsDetailed}
% \mathcal{P}_m = \left\{x,y\geq 0, \quad \sum_i y_{\ell i}=c_{\ell,m},\quad \sum_{j\ne i}(x_{ij}-x_{ji})
% =\alpha_{i,m}+\sum_\ell y_{\ell i}\right\}. 
% \end{equation}
% In other words, for $m$ sufficiently large, 
% \begin{equation}
% \label{eqn:equivalent_of_W_2}
%     \min_{(x,y)\in \mathcal{P}_m}\Psi_m(x,y) = \frac{W_2^2(\overline{G}_m,G)}{t_m},
% \end{equation}
% and it is possible to construct the solution of one problem from the other. 
We retain all nonlocal labels in the classified subsequence, while including labels whose scaled weights converge to zero. Denote these nonlocal labels as $\theta_{\ell,m}$, where $1\leq \ell \leq q_0$.  Set
\[
c_{\ell,m}=\frac{\pi_{\ell,m}}{t_m},
\qquad
\eta_{\ell,m}=\theta_{\ell,m},
\qquad
\alpha_{i,m}=\frac{M_{i,m}-\pi_i^*}{t_m}.
\]
After the chosen subsequence,
\[
c_{\ell,m}\to c_\ell\ge0,
\qquad
\eta_{\ell,m}\to\eta_\ell,
\qquad
\alpha_{i,m}\to\alpha_i.
\]
The labels with $c_\ell>0$ are exactly the atoms of $\nu$. Exact mass conservation gives
\begin{equation}\label{eq:finiteMassBalanceDetailed}
\sum_i\alpha_{i,m}+\sum_{\ell=1}^{q_0}c_{\ell,m}=0.
\end{equation}

Consider the following function $$\Psi_m(x,y) = \sum_{i\ne j}d_{ij}^2x_{ij}
+\sum_{\ell,i}\|\eta_{\ell,m}-\theta_i^*\|^2y_{\ell i},$$ 
along with the constraint
\begin{equation}\label{eq:finiteYrowsDetailed}
\mathcal{P}_m = \left\{x,y\geq 0, \quad \sum_i y_{\ell i}=c_{\ell,m},\quad \sum_{j\ne i}(x_{ij}-x_{ji})
=\alpha_{i,m}+\sum_\ell y_{\ell i}\right\},
\end{equation}
where the constraint $x,y\ge 0$ is understood as $x,y$ are element-wise non-negative. We prove that for $m$ sufficiently large, the optimal transportation problem $W_2^2(\overline{G}_m,G)$ and the problem $\min_{(x,y)\in \mathcal{P}_m}\Psi_m(x,y)$ under constraint \eqref{eq:finiteYrowsDetailed} are equivalent in the sense that
\begin{equation}
\label{eqn:equivalent_of_W_2}
    \min_{(x,y)\in \mathcal{P}_m}\Psi_m(x,y) = \frac{W_2^2(\overline{G}_m,G)}{t_m},
\end{equation}
and it is possible to construct the optimal solution for each problem from the other problem. Noting that $\min_{(x,y)\in \mathcal{P}_m}\Psi_m(x,y)$ is well-defined, as $\Psi_m(x,y)$ is a continuous function satisfying $\lim_{x,y\to \infty}\Psi_m(x,y) = \infty$, and $\mathcal{P}_m$ is a closed set. 

Consider the optimal coupling from $\overline G_m$ to $G_*$, let $\gamma_{ij}^{(m)}$ be the flow from true source $i$ to true target $j$ and let $\zeta_{\ell i}^{(m)}$ be the flow from nonlocal source $\ell$ to true target $i$. Since $W_2^2(\overline G_m,G_*)=O(t_m)$ and distinct true atoms have positive separation $s_*>0$, each $\gamma^{(m)}_{ij}$ for $i\neq j$ or $\zeta^{(m)}_{\ell i}$ are $O(t_m)$.  Define
\[
x_{ij}^{(m)}=\frac{\gamma_{ij}^{(m)}}{t_m},\quad i\ne j,
\qquad
y_{\ell i}^{(m)}=\frac{\zeta_{\ell i}^{(m)}}{t_m}.
\]
The nonlocal row constraints are $\sum_i y_{\ell i}^{(m)}=c_{\ell,m}$ for all $1\leq \ell\leq q_0$. In addition, the flow constraint for each center $\theta^*_i$ as source and target can be expressed as 
\begin{equation*}
    \sum_{1\le j \leq k}\gamma^{(m)}_{ji} +   \sum_{1\le \ell \leq q_0}\zeta^{(m)}_{\ell i}  = \pi_i^*, \quad \sum_{1\leq j\leq k} \gamma^{(m)}_{ij} = M_{i,m}. 
\end{equation*}
By subtracting these two constraints, we have
\begin{equation*}
    \sum_{j\ne i}(x_{ij}^{(m)}-x_{ji}^{(m)})
=\alpha_{i,m}+\sum_\ell y_{\ell i}^{(m)}.
\end{equation*}
Thus, the optimal transportation problem $W_2^2(\overline G_m,G)$ is equivalent to minimize $\Psi_m$ over a subset of $\mathcal{P}_m$, which implies 
\begin{equation*}
    \min_{(x,y)\in \mathcal{P}_m}\Psi_m(x,y) \leq  \frac{W_2^2(\overline{G}_m,G)}{t_m}.
\end{equation*}

Furthermore, as ${W_2^2(\overline{G}_m,G)} \leq \mathcal{C}_0{t_m}$, we have $\min_{(x,y)\in \mathcal{P}_m}\Psi_m(x,y) \leq \mathcal{C}_0$, and for $m$ sufficiently large, there exists a constant $s_0$ such that $\|\eta_{l,m} - \theta_i^*\| \geq s_0$ for all $\ell,i$. It follows that the optimal solution $(x^{(m)},y^{(m)})$ of  $\min_{(x,y)\in \mathcal{P}_m}\Psi_m(x,y)$ satisfies 
\begin{equation*}
    x^{(m)}_{ij} \leq  \frac{\mathcal{C}_0}{\min_{i\neq j}d^2_{ij}}, \quad y^{(m)}_{\ell i} \leq \frac{\mathcal{C}_0}{s_0^2}. 
\end{equation*}
We construct the coupling between $\overline{G}_m$ and $G_*$ from $(x^{(m)},y^{(m)})$ as following. For each $i\neq j$ denote the different indices of true centers, let $\gamma^{(m)}_{ij} = t_mx ^{(m)}_{ij}$ denote the total mass that moves from $\theta^*_{i}$ of $\overline{G}_m$ to $\theta^*_{j}$ of $G_*$, and $\zeta^{(m)}_{\ell i} = t_my ^{(m)}_{\ell i}$ denote the total mass that moves from $\eta_{\ell,m}$ of $\overline{G}_m$ to $\theta^*_{i}$ of $G_*$. The total mass that remains in $\theta_i^*$ is defined as $\gamma_{ii}^{(m)} =\pi_i^*+t_m\alpha_{i,m}-t_m\sum_{j\ne i}x_{ij}$. As $x^{(m)}_{ij}$ is bounded for $m$ sufficiently large and $t_m \downarrow 0$, we have $\gamma_{ii}^{(m)} \geq 0$ for $m$ large enough. We can check that the transportation plan consisting of $\gamma^{(m)}_{ij}$ and $\zeta^{(m)}_{\ell i}$ is actually a coupling between $\overline{G}_m$ and $G_*$. Therefore, we have 
\begin{equation*}
    \min_{(x,y)\in \mathcal{P}_m}\Psi_m(x,y) \geq  \frac{W_2^2(\overline{G}_m,G)}{t_m}.
\end{equation*}
Thus, we establish the equivalence in \eqref{eqn:equivalent_of_W_2} and show how to construct an optimal solution to either problem from an optimal solution to the other.

\medskip\noindent\emph{Step 2: Stability of the finite transportation problem.}
We retain the finitely the nonlocal rows whose limiting masses may be zero. The parameter sequences $(\alpha_{i,m},c_{\ell,m},\eta_{\ell,m})$ satisfy the exact mass balance \eqref{eq:finiteMassBalanceDetailed} and converge to $(\alpha_i,c_\ell,\eta_\ell)$. Lemma~\ref{lem:transportstability} therefore applies to the feasible sets \eqref{eq:finiteYrowsDetailed}, and yields
\begin{equation}\label{eq:collapsedLPlimitDetailed}
\min_{(x,y)\in\mathcal P_m}\Psi_m(x,y)
\longrightarrow\min_{(x,y)\in\mathcal P}\Psi(x,y)
=\mathfrak M_*(\alpha,\nu).
\end{equation}
% In particular the result does not require fixed optimal edges or positive limiting weights for all nonlocal labels. To return from a bounded scaled flow to an actual coupling, set the off-diagonal true-to-true and nonlocal-to-true masses equal to $t_mx_{ij}$ and $t_my_{\ell i}$, and place the remaining mass on the true-node diagonal. The diagonal entries are nonnegative for sufficiently large $m$, since $\pi_i^*>0$ and all scaled flows are $O(1)$.

Finally, Step 1, Lemma~\ref{lem:collapse} and \eqref{eq:collapsedLPlimitDetailed} give
\[
\frac{W_2^2(G_m,G_*)}{t_m}
\longrightarrow\sum_{i=1}^p\operatorname{tr}(C_i)+\mathfrak M_*(\alpha,\nu)
=\Omega(h_\gamma).
\]
The same argument applies to a further subsequence of every original subsequence; uniqueness of the canonical coordinates identifies the same limit throughout. Thus the full sequential assertion holds.
\end{proof}

\begin{remark}[Why $\beta$ does not enter $\Omega$]
The local first-moment coordinate $\beta_i$ represents a common displacement of order $t$ in the recovery path. Its squared transportation cost is $O(t^2)$ and therefore disappears after division by $t$. By contrast, the splitting coordinate $C_i$ is generated by displacements of order $\sqrt t$ and contributes at order $t$. See equation~\eqref{eqn:expansion_of_fitted_center} for a construction. Root-$t$ mass discrepancies between separated true atoms and $O(t)$ nonlocal masses also contribute at order $t$ to $W_2^2$.
\end{remark}

\begin{lemma}[Continuity of the transportation functional]\label{lem:Omegacont}
The map $\Omega: \mathcal{T}_k \to \mathbb{R}_{\geq 0}$ is continuous. In other words, if $h_m,h\in\cT_k$ and $h_m\to h$ in $L^2(P_*)$, then
\[
\Omega(h_m)\to\Omega(h).
\]
\end{lemma}

\begin{proof}

For every fixed $m$, from Proposition~\ref{prop:tangent} we have a family $G_{m,t}\in\cO_k(\Theta)$, $t\downarrow0$, such that
\[
\frac{p_{G_{m,t}}-p_*}{tp_*}\to h_m
\quad\text{in }L^2(P_*).
\]
Applying Proposition~\ref{prop:transportmain} to $G_{m,t}$, we have 
\[
\frac{W_2^2(G_{m,t},G_*)}{t}\to\Omega(h_m).
\]
For each $m$, choose $t_m\in(0,1/m)$ sufficiently small that
\begin{equation}\label{eq:OmegaDiagDensityDetailed}
\left\|
\frac{p_{G_{m,t_m}}-p_*}{t_mp_*}-h_m
\right\|_*\le m^{-1}, \quad 
\left|
\frac{W_2^2(G_{m,t_m},G_*)}{t_m}-\Omega(h_m)
\right|\le m^{-1}.
\end{equation}
Given that $h_m\to h$,
\[
\left\|
\frac{p_{G_{m,t_m}}-p_*}{t_mp_*}-h
\right\|_*
\le m^{-1}+\|h_m-h\|_*\to0.
\]
Proposition~\ref{prop:transportmain} therefore gives
\[
\frac{W_2^2(G_{m,t_m},G_*)}{t_m}\to\Omega(h).
\]
Combining this convergence with the second part of \eqref{eq:OmegaDiagDensityDetailed} yields
\[
|\Omega(h_m)-\Omega(h)|
\le m^{-1}
+\left|
\frac{W_2^2(G_{m,t_m},G_*)}{t_m}-\Omega(h)
\right|
\to0.
\]
Consequently, $\Omega$ is continuous on $\cT_k$ in the $L^2(P_*)$ topology.
\end{proof}

% \section{Proof of the joint Wasserstein weak-limit theorem}\label{sec:maincompletion}
\subsection{Proof of the joint density and Wasserstein weak limit}
\label{sec:maincompletion}
Once we understand the geometry of the sequential first-order
expansion, we prove Theorem~\ref{thm:wasserstein} on the joint
density and Wasserstein weak limit. The argument combines
the density convergence in Theorem~\ref{thm:density} with
the sequential expansion and continuity of $\Omega$ through
the extended continuous mapping theorem.
\begin{proof}[Proof of Theorem~\ref{thm:wasserstein}]
By Theorem~\ref{thm:density}, $\widehat h_n\rightsquigarrow H^\dagger$ in $L^2(P_*)$. For fixed $n$, write
\[
\cH_n=\left\{ h_{n,G} = \sqrt n(p_G-p_*)/p_*:G\in\cO_k(\Theta)\right\}.
\]
Assumption~\ref{ass:strongid} about ordinary mixture identifiability makes the map $$(\mathcal{O}_k(\Theta),W_2)\to (\mathcal{H}_n,L^2(P_*)),\quad G \mapsto h_{n,G}$$ injective. In addition, we also have the inverse map $h_{n,G} \to G$ is continuous. Indeed, if $h_{n,G_m}\to h_{n,G}$ in $L^2(P_*)$, then $p_{G_m}\to p_G$ in $L^1(\mu)$ by Cauchy--Schwarz with respect to $P_*$. Every subsequence of $G_m$ has a further $W_2$-convergent subsequence by compactness of $\cO_k(\Theta)$. Continuity of the kernel mixture map identifies its density limit as $p_G$, and ordinary identifiability identifies the measure limit as $G$. Therefore $G_m\to G$ in $W_2$. In particular, the finite-$n$ statistic below is Borel measurable on $\cH_n$.

Define $F_n:\cH_n\to L^2(P_*)\times\mathbb R$ by
\[
F_n(h_{n,G})=\bigl(h_{n,G},\sqrt n\,W_2^2(G,G_*)\bigr)
\]
and define $F:\cT_k\to L^2(P_*)\times\mathbb R$ by $F(h)=(h,\Omega(h))$. For any deterministic $h_n\in\cH_n$ converging to $h\in\cT_k$, write $h_n=(p_{G_n}-p_*)/(t_np_*)$ with $t_n=n^{-1/2}$. Proposition~\ref{prop:transportmain} gives
\[
\sqrt n W_2^2(G_n,G_*)=\frac{W_2^2(G_n,G_*)}{t_n}
\longrightarrow\Omega(h),
\]
so $F_n(h_n)\to F(h)$. The support of $H^\dagger$ is contained in $\cT_k$ and $F$ is continuous there by Lemma~\ref{lem:Omegacont}. The extended continuous mapping theorem for maps on moving domains therefore implies
\[
F_n(\widehat h_n)\rightsquigarrow F(H^\dagger),
\]
which is precisely the joint convergence in the theorem. Applying the continuous square-root map on $[0,\infty)$, we have 
\begin{equation*}
    n^{1/4}W_2(\widehat G_n,G_*) \rightsquigarrow \{\Omega(H^{\dagger})\}^{1/2}. 
\end{equation*}
As a result, $W_2(\widehat G_n,G_*) = O_{\mathbb{P}}(n^{-1/4})$. This completes the proof. 
\end{proof}

% \section{Proofs for bootstrap inference}\label{sec:bootproofs}
\section{Proofs for Section~\ref{sec:bootstrap}}
\label{app:bootstrapproofs}
% \subsection{Proof of Theorem~\ref{thm:bootstrap}}
\begin{proof}[Proof of Theorem~\ref{thm:bootstrap}] We divide the proof into following steps.

\textit{Step 1: Consistency of $Q_n^{\otimes m}$.} 
Recall that $Q_n=P_{\widehat G_n}$, we prove that under the assumption $m = o(n)$, we have $\|Q_n^{\otimes m}-P_*^{\otimes m}\|_{\mathrm{TV}}
\xrightarrow{P}0$. Assumption~\ref{ass:smooth} includes positivity of $p_*$, so $Q_n=P_{\widehat G_n}$ is absolutely continuous with respect to $P_*$. Therefore, the following density ration is well-defined
\[
R_n(x)
:=
\frac{dQ_n}{dP_*}(x)
=
\frac{p_{\widehat G_n}(x)}{p_*(x)}.
\]
As $Q_n$ is a probability measure, $P_*R_n=1$. Pearson's $\chi^2$ divergence therefore can be expressed as
\[ \chi^2(Q_n\|P_*) = P_*(R_n-1)^2= \int\left(\frac{p_{\widehat G_n}-p_*}{p_*}\right)^2 p_*\,d\mu= \rho(\widehat G_n)^2.\]
The tensor-product identity for chi-square divergence gives
\begin{equation*}
    1+\chi^2(Q_n^{\otimes m}\Vert P_*^{\otimes m})
=\{1+\rho(\widehat G_n)^2\}^{m}.
\end{equation*}
Therefore, we have the following comparison
\begin{equation}\label{eqn:dung_TV_comparison}
2\|Q_n^{\otimes m}-P_*^{\otimes m}\|_{\mathrm{TV}} \leq \{\chi^2(Q_n^{\otimes m}\|P_*^{\otimes m})\}^{1/2}
\le\left[\{1+\rho(\widehat G_n)^2\}^m-1\right]^{1/2}.
\end{equation}
The root-$n$ density bound in Lemma~\ref{lem:rootn} yields $n\rho(\widehat G_n)^2=O_P(1)$. Thus $m\rho(\widehat G_n)^2=o_P(1)$ whenever $m/n\to0$, and
\begin{equation}\label{eqn:dung_1_+_distance_to_0}
\{1+\rho(\widehat G_n)^2\}^m-1
\le\exp\{m\rho(\widehat G_n)^2\}-1=o_P(1).
\end{equation}
Combining this result with equation~\eqref{eqn:dung_TV_comparison}, we have $\|Q_n^{\otimes m}-P_*^{\otimes m}\|_{\mathrm{TV}} \pto 0$. 

\textit{Step 2: Prove equation~\eqref{eq:explicitBLbound}}. Put
\[
S_{n,m}^*=m^{1/4}W_2(\widehat G^*_m,G_*),
\quad T_m^0=m^{1/4}W_2(\widetilde G_m,G_*),
\]
We estimate $T^*_{n,m}$ through $S^*_{n,m}$. Using the triangle inequality, we have
\[
|T_{n,m}^*-S_{n,m}^*|\le m^{1/4}W_2(\widehat G_n,G_*).
\]

Denote $E^*f:=E[f\mid X_1,\ldots,X_n]$, $P^*$ and $d^*_{\mathrm{BL}}$ be the
conditional expectation of $f$, conditional probability, and the conditional $d_{\mathrm{BL}}$-distance given sample
$X_1,\ldots,X_n$. For any $f\in\mathrm{BL}_1$, from estimation above,
we have
\begin{equation}
\label{eqn:dung_first_equation_of_f_T}
|E^*f(T_{n,m}^*)-E^*f(S_{n,m}^*)|
\le
E^*[|T_{n,m}^*-S_{n,m}^*|]
\le
m^{1/4}W_2(\widehat G_n,G_*).
\end{equation}
% Theorem~4.2 gives
% $T_n=n^{1/4}W_2(\widehat G_n,G_*)=O_P(1)$,
% and from the hypothesis $m/n\to0$, we have
% \begin{equation}
% m^{1/4}W_2(\widehat G_n,G_*)
% =
% \left(\frac mn\right)^{1/4}T_n
% =
% o_P(1).
% \tag{89}
% \end{equation}
Consider the measurable statistic on $\mathcal X^m$
\[
\Phi_m(x_1,\ldots,x_m)
=
m^{1/4}W_2\{\mathcal A_m(x_1,\ldots,x_m),G_*\}.
\]
Measurability follows from the assumed measurable MLE
selection and continuity of $W_2$ on the compact
mixing-measure space.
Under $Q_n^{\otimes m}$, $\Phi_m$ has conditional law
equal to the law of $S_{n,m}^*$; under $P_*^{\otimes m}$
it has the law of $T_m^0$. By applying following estimation for any measurable $g$
with $\|g\|_\infty\le1$,
\[
\left|
\int g\,dQ-\int g\,dP
\right|
\le
2\|Q-P\|_{\mathrm{TV}},
\]
to $g=f\circ\Phi_m$ for $f\in\mathrm{BL}_1$,
we obtain
\begin{equation}
    \label{eqn:dung_second_equation_of_f_T}
|E^*f(S_{n,m}^*)-Ef(T_m^0)|
\le
2\|Q_n^{\otimes m}-P_*^{\otimes m}\|_{\mathrm{TV}}.
\end{equation}
% Using the total-variation comparison in \eqref{eqn:dung_TV_comparison}, we have

%     |E^*f(S_{n,m}^*)-Ef(T_m^0)| \le \left[\{1+\rho(\widehat G_n)^2\}^m-1\right]^{1/2}.

By combining the estimations from \eqref{eqn:dung_first_equation_of_f_T} and \eqref{eqn:dung_second_equation_of_f_T}, we have for all $f \in \mathrm{BL}_1$: 
\begin{equation*}
    |E^*f(T_{n,m}^*)-Ef(T_m^0)| \leq m^{1/4}W_2(\widehat G_n,G_*) + 2\|Q_n^{\otimes m}-P_*^{\otimes m}\|_{\mathrm{TV}}. 
\end{equation*}
Therefore we achieve estimation \eqref{eq:explicitBLbound}. 

\textit{Step 3: Conclusion}. We give the proof of the bounds \eqref{eq:bootlimit} and \eqref{eq:bootfiniteNtarget}. Firstly, we observe that  Theorem~\ref{thm:wasserstein} gives
$T_n=n^{1/4}W_2(\widehat G_n,G_*)=O_P(1)$,
and from the hypothesis $m/n\to0$, we have
\begin{equation*}
m^{1/4}W_2(\widehat G_n,G_*)
=
\left(\frac mn\right)^{1/4}T_n
=
o_P(1).
\end{equation*}
Combining this estimation with Step 1, we have two terms in the RHS of \eqref{eq:explicitBLbound} converge in probability to 0. Hence $d_{\mathrm{BL}}
        \{\mathcal L^*(T_{n,m}^*),\mathcal L(T_m^0)
    \} \pto 0$.

% The bounded-Lipschitz distance between the conditional laws of these two statistics is bounded by the same right-hand side. Meanwhile, pushing the two product laws through the common measurable statistic and applying \eqref{eqn:dung_TV_comparison} gives

Theorem~\ref{thm:wasserstein} applied to an independent size-$m$ sample gives $T_m^0\rightsquigarrow T_\infty$. Since bounded-Lipschitz distance metrizes weak convergence on $\mathbb R$, $d_{\mathrm{BL}}\{\mathcal L(T_m^0),\mathcal L(T_\infty)\}\to0$. Using triangle inequality for $d_{\mathrm{BL}}$, we achieve immediately $d_{\mathrm{BL}}\{\mathcal L^*(T_{n,m}^*),\mathcal L(T_\infty)\} \pto 0$, which proves the inequality~\eqref{eq:bootlimit}. In addition, a second application of Theorem~\ref{thm:wasserstein} gives
$d_{\mathrm{BL}}\{\mathcal L(T_n),\mathcal L(T_\infty)\}\pto0$, which leads to $d_{\mathrm{BL}}\{\mathcal L^*(T_{n,m}^*),\mathcal L(T_n)\} \pto 0$
and therefore establishes \eqref{eq:bootfiniteNtarget}.

% The triangle inequality proves the explicit bound \eqref{eq:explicitBLbound}. Its first term is
% $(m/n)^{1/4}T_n=o_P(1)$ by Theorem~\ref{thm:wasserstein}; its second term is $o_P(1)$ by \eqref{eqn:dung_1_+_distance_to_0}. Hence
% \begin{equation}\label{eq:proofBootIntermediate}
% d_{\mathrm{BL}}\{\mathcal L^*(T_{n,m}^*),\mathcal L(T_m^0)\}\xrightarrow{P}0.
% \end{equation}
% Theorem~\ref{thm:wasserstein} applied to an independent size-$m$ sample gives $T_m^0\rightsquigarrow T_\infty$. Since bounded-Lipschitz distance metrizes weak convergence on $\mathbb R$,
% \[
% d_{\mathrm{BL}}\{\mathcal L(T_m^0),\mathcal L(T_\infty)\}\to0.
% \]
% Combining this deterministic convergence with \eqref{eq:proofBootIntermediate} establishes \eqref{eq:bootlimit}.. All displayed distances are between \emph{probability laws}, not random variables.
\end{proof}

\begin{proof}[Proof of Corollary~\ref{cor:bootsquared}]
Recall from equation~\eqref{eq:bootlimit} that
$d_{\mathrm{BL}}\{\mathcal L^*(T_{n,m}^*),
\mathcal L(T_\infty)\}\xrightarrow{P}0$.
As $L(T_\infty)$ and $H^\dagger$ follow the same
distribution, with the choose of $H^\dagger$ such that
$H^\dagger$ is independent of sample $X_1,\ldots,X_n$.
Using Lemma~\ref{lem:dung_BL_distance_compare}, we have
\[
d_{\mathrm{BL}}^*
\bigl((T_{n,m}^*)^2,\Omega(H^\dagger)\bigr)
\le
(2R+4)d_{\mathrm{BL}}
\{\mathcal L^*(T_{n,m}^*),\mathcal L(T_\infty)\}
+
4\PP\bigl(\Omega(H^\dagger)>R^2\bigr).
\]
By first letting $n\to\infty$, then $R\to\infty$,
we achieve that $d_{\mathrm{BL}}^*
\bigl((T_{n,m}^*)^2,\Omega(H^\dagger)\bigr)
\xrightarrow{P}0$.
As $d_{\mathrm{BL}}$ metricizes the weak convergence,
we achieve that
\[
\mathcal L^*
\left\{
\sqrt m\,W_2^2(\widehat G_m^*,\widehat G_n)
\right\}
\rightsquigarrow_P
\mathcal L\{\Omega(H^\dagger)\}.
\]
This proves equation~\eqref{eq:bootsquared}.
\end{proof}
\begin{proof}[Proof of Corollary~\ref{cor:bootquantile}] We divide the proof into following steps. 

\emph{Step 1: Convergence of quantile $q_{n,m}^*(p_0)$.}
We prove that $q_{n,m}^*(p_0)$ converges in probability to $q_{p_0}$. According to Theorem~\ref{thm:bootstrap}, we have
$d_{\mathrm{BL}}\{T_{n,m}^*,T_\infty\}\xrightarrow{P}0$.
Let $t$ be any continuity point of $F_\infty$, we show
that $F_{n,m}^*(t)\xrightarrow{P}F_\infty(t)$.
For $\varepsilon\in(0,1)$, using Lemma~\ref{lem:compare_kolmogorov_BL}, we have
\[
|F_{n,m}^*(t)-F_\infty(t)|
\le
\varepsilon^{-1}d_{\mathrm{BL}}^*(T_{n,m}^*,T_\infty)
+
\bigl(F_Q(t+\varepsilon)-F_Q(t-\varepsilon)\bigr).
\]
By first letting $n\to\infty$, then $\varepsilon\to0$,
we achieve that
$|F_{n,m}^*(t)-F_\infty(t)|\xrightarrow{P}0$.

Fix $\varepsilon>0$ and write $p_0=1-\alpha$.
By hypothesis \eqref{eq:quantreg},
$F_\infty(q_{p_0}-\varepsilon)<p_0
<F_\infty(q_{p_0}+\varepsilon)$.
As distribution function has at most countably many
discontinuities, we can choose continuity points $a,b$
satisfying $q_{p_0}-\varepsilon<a<q_{p_0}<b<q_{p_0}+\varepsilon$ and still $F_\infty(a)<p_0<F_\infty(b)$.
Let
$\delta=\frac12\min\{p_0-F_\infty(a),F_\infty(b)-p_0\}>0$,
by the convergence of CDF at continuous point of
$F_\infty$ just established, we have
\[
F_{n,m}^*(a)\xrightarrow{P}F_\infty(a),
\qquad
F_{n,m}^*(b)\xrightarrow{P}F_\infty(b).
\]
Hence with probability tending to one,
\[
F_{n,m}^*(a)\le F_\infty(a)+\delta<p_0
\quad\text{and}\quad
F_{n,m}^*(b)\ge F_\infty(b)-\delta>p_0.
\]
By the definition of the generalized inverse,
$F_{n,m}^*(a)<p_0$ implies $q_{n,m}^*(p_0)>a$,
while $F_{n,m}^*(b)>p_0$ leads to
$q_{n,m}^*(p_0)\le b$.
Thus, with probability tending to one,
\[
q_{p_0}-\varepsilon<a<q_{n,m}^*(p_0)
\le b<q_{p_0}+\varepsilon.
\]
Since $\varepsilon$ is arbitrary, this proves $q_{n,m}^*(p_0)\pto q_{p_0}$.

\emph{Step 2: Consistency of random Wasserstein ball.} We prove the equation~\eqref{eq:coverage} and~\eqref{eq:confball}.
Theorem~\ref{thm:wasserstein} gives $T_n\rightsquigarrow T_\infty$,
and the quantile consistency established above gives
$q_{n,m}^*(p_0)\xrightarrow{P}q_{p_0}$.
Noting that $q_{p_0}$ is a constant, we employ
Slutsky's theorem to achieve
\[
(T_n,q_{n,m}^*(p_0))
\rightsquigarrow
(T_\infty,q_{p_0}).
\]
Consider the closed half-space
$A=\{(x,y)\in\RR^2:x\le y\}$, whose boundary is
$\{x=y\}$.
Under the limiting law, the boundary event
$(T_\infty,q_{p_0})\in A$ is exactly
$\{T_\infty=q_{p_0}\}$, which has probability zero
by hypothesis~\eqref{eq:quantreg}. The Portmanteau theorem therefore gives
\[
\PP\{T_n\le q_{n,m}^*(p_0)\}
\to
\PP\{T_\infty\le q_{p_0}\}
=
F_\infty(q_{p_0})
=
p_0.
\]
This proves equation~\eqref{eq:coverage}.
To prove equation~\eqref{eq:confball}, from the definition of $T_n$, we have
\[
G_*\in\mathcal C_{n,1-\alpha}
\iff
W_2(G_*,\widehat G_n)
\le n^{-1/4}q_{n,m}^*(1-\alpha)
\iff
T_n\le q_{n,m}^*(1-\alpha).
\]
Therefore, equation~\eqref{eq:confball} follows from
\[
\PP(G_*\in\mathcal C_{n,1-\alpha})
=
\PP(T_n\le q_{n,m}^*(1-\alpha))
\to p_0.
\]
This completes our proof. 
% Write $p_0=1-\alpha$, $q=q_{p_0}$ and $F=F_\infty$. Under \eqref{eq:quantcross}, for each $\epsilon>0$ one can choose continuity points $a,b$ of $F$ satisfying
% \[
% q-\epsilon<a<q<b<q+\epsilon,
% \qquad F(a)<p_0<F(b).
% \]
% Theorem~\ref{thm:bootstrap} gives
% $F_{n,m}^*(a)\xrightarrow{P}F(a)$ and $F_{n,m}^*(b)\xrightarrow{P}F(b)$. Thus with probability tending to one,
% $F_{n,m}^*(a)<p_0<F_{n,m}^*(b)$, and the generalized inverse obeys
% $a<q_{n,m}^*(p_0)\le b$. Since $\epsilon$ is arbitrary, this proves \eqref{eq:qconsistency}.

% By Theorem~\ref{thm:wasserstein} and Slutsky's theorem,
% $(T_n,q_{n,m}^*(p_0))\rightsquigarrow(T_\infty,q)$.
% The boundary of $\{(t,z):t\le z\}$ has limit probability $\Pr(T_\infty=q)=0$ by \eqref{eq:quantreg}. The Portmanteau theorem therefore yields
% \[
% \Pr\{T_n\le q_{n,m}^*(p_0)\}\to F(q)=p_0.
% \]
% This is \eqref{eq:coverage}. The event $\{T_n\le q_{n,m}^*(p_0)\}$ is exactly $\{G_*\in\mathcal C_{n,1-\alpha}\}$, proving \eqref{eq:confcoverage}. The argument uses \emph{both} the quantile-crossing condition and the no-atom condition; unconditional weak convergence alone does not imply exact coverage.
\end{proof}

\begin{proof}[Proof of Proposition~\ref{prop:MCbootstrap}]
% Conditionally on the observed data, the $B$ bootstrap statistics are i.i.d. with distribution function $F_{n,m}^*$. The Dvoretzky--Kiefer--Wolfowitz inequality gives, for every $\epsilon>0$,
% \begin{equation}\label{eqn:dung_DKW_bound}
% \Pr^*\left\{\sup_t|\widehat F_{n,m,B}^*(t)-F_{n,m}^*(t)|>\epsilon\right\}
% \le2\exp(-2B\epsilon^2).
% \end{equation}
% Taking expectations yields \eqref{eq:DKWconv} when $B\to\infty$. For $a,b$ as in the preceding proof, combine \eqref{eqn:dung_DKW_bound} with $F_{n,m}^*(a)\to_P F(a)$ and $F_{n,m}^*(b)\to_P F(b)$ to obtain
% \[
% \Pr\{\widehat F_{n,m,B}^*(a)<p_0<\widehat F_{n,m,B}^*(b)\}\to1.
% \]
% Hence $a<\widehat q_{n,m,B}^*(p_0)\le b$ with probability tending to one, proving \eqref{eq:MCquantile}. Under \eqref{eq:quantreg}, Slutsky and Portmanteau give
% \[
% \Pr\{T_n\le\widehat q_{n,m,B}^*(1-\alpha)\}\longrightarrow1-\alpha.
% \]
% This is exactly the event $G_*\in\widehat{\mathcal C}_{n,1-\alpha}$, proving \eqref{eq:MCball}.
We divide the proof into following steps. 

\emph{Step 1: Consistency of Monte Carlo approximation.} We prove equation~\eqref{eq:DKWconv}. Given $\mathcal F_n$, the $B=B_n$ bootstrap repetitions
are independent and identically distributed, and each
statistic $T_{n,m}^{*(b)}$ has CDF $F_{n,m}^*$.
Thus, conditional on $\mathcal F_n$,
$\widehat F_{n,m,B}^*$ is the empirical CDF based on
$B$ i.i.d. observations from $F_{n,m}^*$.
The DKW inequality (see \citet{DvoretzkyKieferWolfowitz1956,Massart1990}) gives
\begin{equation}
\label{eqn:DKW_bound}
P^*\left\{
\sup_t
\left|
\widehat F_{n,m,B}^*(t)-F_{n,m}^*(t)
\right|>\varepsilon
\right\}
\le 2e^{-2B\varepsilon^2},
\end{equation}
where the right side is deterministic and tends to zero
because $B_n\to\infty$. Taking unconditional expectations gives
\[
\PP\left\{
\sup_t
\left|
\widehat F_{n,m,B}^*(t)-F_{n,m}^*(t)
\right|>\varepsilon
\right\}
\le 2e^{-2B_n\varepsilon^2}\to0,
\]
which proves equation~\eqref{eq:DKWconv}.

\emph{Step 2: Asymptotic behavior of quantiles and coverage.} We study the asymptotic behavior of $\widehat q_{n,m,B}^*(p_0)$ and $\widehat{\mathcal C}_{n,1-\alpha}$.
Fix $p_0=1-\alpha$ and $\varepsilon>0$.
Choose continuity points $a,b$ exactly as in the proof
of Corollary~\ref{cor:bootquantile} so that
\[
q_{p_0}-\varepsilon<a<q_{p_0}<b<q_{p_0}+\varepsilon,
\qquad
F_\infty(a)<p_0<F_\infty(b).
\]
Theorem~\ref{thm:bootstrap} yields
$F_{n,m}^*(a)\xrightarrow{P}F_\infty(a)$ and
$F_{n,m}^*(b)\xrightarrow{P}F_\infty(b)$,
while the DKW bound \eqref{eqn:DKW_bound} gives
\[
\widehat F_{n,m,B}^*(a)-F_{n,m}^*(a)
\xrightarrow{P}0,
\qquad
\widehat F_{n,m,B}^*(b)-F_{n,m}^*(b)
\xrightarrow{P}0.
\]
Therefore
$\widehat F_{n,m,B}^*(a)<p_0<
\widehat F_{n,m,B}^*(b)$
with probability tending to one.
Moreover, the generalized-inverse definition then forces
$a<\widehat q_{n,m,B}^*(p_0)\le b$
with probability tending to one.
Since $a,b$ lie within $\varepsilon$ of $q_{p_0}$
and $\varepsilon$ is arbitrary, we achieve the limit~\eqref{eq:MCquantile} as
\[
\widehat q_{n,m,B}^*(p_0)\xrightarrow{P}q_{p_0}.
\]
Using Slutsky theorem for
$T_n\rightsquigarrow T_\infty$ and the equation \eqref{eq:MCquantile},
we have
$(T_n,\widehat q_{n,m,B}^*(p_0))
\rightsquigarrow(T_\infty,q_{p_0})$.
By considering $A=\{(x,y)\in\RR^2:x\le y\}$
and using the same argument as in Step 2 in the proof
of Corollary~\ref{cor:bootquantile}, we have
\[
\PP\{T_n\le\widehat q_{n,m,B}^*(1-\alpha)\}
\to
\PP\{T_\infty\le q_{p_0}\}
=
1-\alpha.
\]
This event is exactly
$\{G_*\in\widehat{\mathcal C}_{n,1-\alpha}\}$
by the definition of the ball, proving the equation~\eqref{eq:MCball}.
\end{proof}
\section{Auxiliary Results}
\label{sec:auxiliary_results}
\begin{proof}[Proof of Corollary~\ref{cor:gaussian}]

In this section, we denote 
\[
f_\theta(x)
=(2\pi)^{-d/2}|\Sigma|^{-1/2}
\exp\left\{-\frac12(x-\theta)^\top\Sigma^{-1}(x-\theta)\right\}.
\]

\textit{Step 1: Verify Assumption~\ref{ass:smooth}.}
It is obvious to check that the gradient of Gaussian density with respect to $\theta$ is 
\[
\nabla_\theta f_\theta(x)
=\Sigma^{-1}(x-\theta)f_\theta(x),
\]
and every higher derivative has the form
\[
\nabla_\theta^r f_\theta(x)=P_r(x-\theta)f_\theta(x),
\qquad r=0,1,2,3,
\]
where $P_r$ is a tensor-valued polynomial of degree at most $r$ whose coefficients depend only on $\Sigma^{-1}$.

Fix any true atom, without loss of generality, w can choose $\theta_1^*$. Since $p_*(x)\ge\pi_{\min}^*f_{\theta_1^*}(x)$,
\[
\frac{f_\theta(x)}{p_*(x)}
\le(\pi_{\min}^*)^{-1}\frac{f_\theta(x)}{f_{\theta_1^*}(x)}.
\]
A direct calculation gives
\[
\log\frac{f_\theta(x)}{f_{\theta_1^*}(x)}
=x^\top\Sigma^{-1}(\theta-\theta_1^*)
-\frac12\{\theta^\top\Sigma^{-1}\theta-(\theta_1^*)^\top\Sigma^{-1}\theta_1^*\}.
\]
As $\Theta$ is compact, there exist constants $C_1,C_2$ such that uniformly in $\theta\in\Theta$, ${f_\theta(x)}/{p_*(x)}\le C_1e^{C_2\|x\|}$.
Combining this with the polynomial derivative representation and compactness of $\Theta$, there exists constant $C>0$ such that 
\begin{equation}\label{eq:GaussianEnvelopeDetailed}
\sup_{\theta\in\Theta}
\frac{\|\nabla_\theta^r f_\theta(x)\|}{p_*(x)}
\le C(1+\|x\|^3)e^{C\|x\|},
\qquad r\le3.
\end{equation}
Fix finite $a>0$ and an integer $m\ge1$, and let $h(x)=(1+\|x\|^3)^m e^{a\|x\|}$. Compactness of $\Theta$ and the Gaussian moment-generating bound imply
\[
\sup_{\theta\in\Theta}\int_{\mathbb R^d}h(x)f_\theta(x)\,dx<\infty,
\qquad
P_*h=\sum_{i=1}^p\pi_i^*\int_{\mathbb R^d}h(x)f_{\theta_i^*}(x)\,dx<\infty.
\]
Consequently the right side of \eqref{eq:GaussianEnvelopeDetailed} belongs to $L^q(P_*)$ for every finite $q$. Assumption~\ref{ass:smooth} follows for any desired $\eta>0$ with $F(x) = C(1+\|x\|^3)e^{C\|x\|}$.

\textit{Step 2: Verify Assumption~\ref{ass:consistency}.}
First we construct an integrable envelope for the log-likelihood class for the class $\{\log p_G:G\in \mathcal{O}_k(\Theta)\}.$ For the upper bound $\ell$, every Gaussian density with covariance $\Sigma$ is bounded by the same finite constant ${(2\pi)^{-d/2}(\det \Sigma)^{-1/2}}$, so
\[
\log p_G(x)\le -\frac{d}{2}\log(2\pi) - \frac{1}{2}\log(\det \Sigma). 
\]
uniformly in $G\in\cO_k(\Theta)$. For the lower bound, in any $k$-component representation, there exists one weight, whose mass is at least $1/k$. If that component has location $\theta\in\Theta$, then $p_G(x)\ge k^{-1}f_\theta(x)$. Furthermore, compactness of $\Theta$ gives $\inf_{\theta\in\Theta}\log f_\theta(x)
\ge -C(1+\|x\|^2)$. Therefore, there exists a constant $C$ depending only on $\Theta$ such that 
\[
-C(1+\|x\|^2)\le\log p_G(x)\le C
\]
uniformly in $G\in\cO_k(\Theta)$. This envelope is integrable under the finite Gaussian mixture $P_*$.

Moreover, the parameter space of $G$, given as  $\Delta_{k-1}\times\Theta^k$, is a compact set. For each $x$, the map
\[
(\pi,\theta_1,\ldots,\theta_k)
\mapsto
\log\left\{\sum_{j=1}^k\pi_jf_{\theta_j}(x)\right\}
\]
is continuous for every $x$. The compact parametrization, pointwise continuity, and integrable envelope establish the $P_*$-Glivenko--Cantelli property. 

Now suppose that $G,H$ produce the same Gaussian location-mixture density. In this case, Fourier transformation gives
\[
e^{-\frac12t^\top\Sigma t}\int e^{-it^\top\theta}\,dG(\theta)
=e^{-\frac12t^\top\Sigma t}\int e^{-it^\top\theta}\,dH(\theta).
\]
The common Gaussian factor is never zero, so the characteristic functions of $G$ and $H$ coincide
\begin{equation*}
    \int e^{-it^\top\theta}\,dG(\theta) = \int e^{-it^\top\theta}\,dH(\theta).
\end{equation*}
Uniqueness of characteristic functions gives $G=H$. Thus the population likelihood has the unique maximizer required in Assumption~\ref{ass:consistency}.

\textit{Step 3: Verify Assumption~\ref{ass:strongid}. } Suppose that for pairwise distinct $\xi_1,\ldots,\xi_m$,
\begin{equation}
\label{eqn:dung_equal_identifiabitily_identity}
    \sum_{j=1}^m\left[
 a_jf_{\xi_j}
+b_j^\top\nabla_\theta f_{\xi_j}
+\frac12\operatorname{tr}(C_j\nabla_\theta^2f_{\xi_j})
\right]=0,
\end{equation}
where $C_1,\ldots,C_M$ are symmetric matrices. We show that all the coefficients $a_j$, $b_j$, and $C_j$ are equal to 0. In this proof, we use the Fourier transform convention $\widehat g(t)=\int e^{-it^\top x}g(x)\,dx$. For a Gaussian location kernel $f_{\theta}(x)$, the Fourier transformation with respect to $x$ of $f$ and its derivatives are  
\[
\widehat f_{\xi_j}(t)=e^{-it^\top\xi_j}e^{-\frac12t^\top\Sigma t},
\qquad
\widehat{\nabla_\theta f_{\xi_j}}(t)=-it\,\widehat f_{\xi_j}(t),
\qquad
\widehat{\nabla_\theta^2 f_{\xi_j}}(t)=-tt^\top\widehat f_{\xi_j}(t).
\]
By taking the Fourier transformation to both side of identity~\eqref{eqn:dung_equal_identifiabitily_identity} and dividing by the everywhere nonzero Gaussian factor $e^{-\frac12t^\top\Sigma t}$, we have 
\begin{equation}\label{eq:GaussianExpPolynomialDetailed}
\sum_{j=1}^m e^{-it^\top\xi_j}P_j(t)=0,
\qquad t\in\RR^d,
\end{equation}
with $P_j(t)=a_j-i b_j^\top t-\frac12 t^\top C_jt.$  Now we consider $\mathcal V
=\RR^d\setminus\bigcup_{j\ne\ell}
\{v:v^\top(\xi_j-\xi_\ell)=0\}.$
The complement is a finite union of proper hyperplanes, so $\mathcal V$ is open and dense. Fix an arbitrary $v\in\mathcal V$ and put $\lambda_j(v)=v^\top\xi_j$. These numbers are pairwise distinct. Restricting \eqref{eq:GaussianExpPolynomialDetailed} to $t=sv$ gives
\begin{equation}
\label{eqn:short_form_of_identity}
    \sum_{j=1}^m e^{-is\lambda_j(v)}q_{j,v}(s)=0,
\end{equation}
where every $q_{j,v}:=P_j(sv)$ is a polynomial of degree at most two. For the fixed $v\in\mathcal V$ and a fixed index $j$, apply the differential operator
$\mathcal L_{j,v}=\prod_{\ell\ne j}\left(\frac{d}{ds}+i\lambda_\ell(v)\right)^3$,
to the preceding identity. Note that we have the following identity $$\left(\frac{d}{ds}+ib\right)(e^{ia}f) = e^{ia}\left(\left(\frac{d}{ds} + i(a+b)\right)f\right),$$ therefore for $\ell\ne j$,
\[
\left(\frac{d}{ds}+i\lambda_\ell(v))\right)^3
\{e^{-is\lambda_\ell(v)}q_{\ell,v}(s)\}
=e^{-is\lambda_\ell(v)}q_{\ell,v}^{(3)}(s)=0,
\]
because $\deg(q_{\ell,v})\le2$. Therefore,
\[
e^{-is\lambda_j(v)}
\prod_{\ell\ne j}\left(\frac{d}{ds}+i(\lambda_\ell(v)- \lambda_j(v))\right)^3q_{j,v}(s)=0.
\]
Furthermore, if $c\ne0$, the linear operator $d/ds+ic$ is injective on the set of polynomials: the highest-degree term of $icq$ cannot be cancelled by $q'$. Hence the product operator is injective and $q_{j,v}\equiv0$. As a result, the identity $q_{j,v}\equiv0$ holds for $j$ and for every $v\in\mathcal V$.

Taking $s=1$ shows that $P_j(v)=0$ for every $v\in\mathcal V$. Since $\mathcal V$ is dense and $P_j$ is a polynomial, $P_j$ vanishes identically on $\RR^d$. Consequently, the formula $P_j(t)=a_j-i b_j^\top t-\frac12t^\top C_jt$ and the fact that $C_j$ is symmetric 
then gives $a_j=0$, $b_j=0$, and $C_j=0$ for each $j$ (the matrices $C_j$ are symmetric). This establishes second-order strong identifiability.
\end{proof}

\begin{lemma}
[Criterion for Hausdorff convergence]
\label{lem:Haussdorf}
Let $(X,d)$ be a metric space, and let $A,A_n\subset X$,
$n\geq 1$, be nonempty compact sets. Suppose that:
\begin{enumerate}
    \item There exists a compact set $K\subset X$ such that
    $A_n\subset K$ for every $n$.
    
    \item For every strictly increasing sequence $(n_m)$ and
    every sequence $y_m\in A_{n_m}$, if $y_m\to y$ in $X$,
    then $y\in A$.
    
    \item There exists $A'\subset A$ with $\overline{A'}=A$
    such that, for every $x\in A'$, there is a sequence
    $x_n\in A_n$ satisfying $x_n\to x$.
\end{enumerate}
Then
\[
    d_H(A_n,A)\longrightarrow 0,
\]
where
\[
    d_H(B,C)
    :=
    \max\left\{
        \sup_{b\in B}\inf_{c\in C}d(b,c),
        \sup_{c\in C}\inf_{b\in B}d(c,b)
    \right\}.
\]
\end{lemma}

\begin{proof}
Write $d(x,B):=\inf_{b\in B}d(x,b)$. We first show that
\[
    \sup_{y\in A_n}d(y,A)\longrightarrow 0.
\]
Otherwise, there exist $\varepsilon>0$, a strictly increasing
sequence $(n_m)$, and points $y_m\in A_{n_m}$ such that
$d(y_m,A)\geq\varepsilon$ for every $m$.
Since $y_m\in K$ and $K$ is compact, a subsequence of $(y_m)$
converges to some $y\in K$. By the second assumption, $y\in A$.
Along this subsequence,
\[
    d(y_m,A)\leq d(y_m,y)\longrightarrow 0,
\]
which is a contradiction.

Next, fix $\varepsilon>0$. Since $A'$ is dense in $A$ and
$A$ is compact, there exist $z_1,\ldots,z_r\in A'$ such that
\[
    A\subset \bigcup_{i=1}^r B(z_i,\varepsilon/2).
\]
For each $i$, choose $z_{i,n}\in A_n$ such that
$z_{i,n}\to z_i$. Since there are finitely many such sequences,
for all sufficiently large $n$,
\[
    \max_{1\leq i\leq r}d(z_i,z_{i,n})<\varepsilon/2.
\]
For any $x\in A$, choose $i$ with $d(x,z_i)<\varepsilon/2$.
Then
\[
    d(x,A_n)
    \leq d(x,z_{i,n})
    \leq d(x,z_i)+d(z_i,z_{i,n})
    <\varepsilon.
\]
Thus $\sup_{x\in A}d(x,A_n)\to 0$. Combining the two
directed distance bounds gives $d_H(A_n,A)\to 0$.
\end{proof}

% \begin{lemma}
% \label{lem:compact_convergence}
%     Let $q$ a continuous on $K_M$. Consider a family of compact sets $A_n,A\subset K_M$ satisfy $d_H(A_n,A)\to0$, and let $q_n$ are bounded continuous functions on $A_n$ satisfying $\sup_{x\in A_n}|q_n-q|\to0$. Then,
%     \begin{enumerate}
%         \item [(i)] The maximum value of $q_n$ over $A_n$ converges to maximum value of $q$ over $A$: 
%         \begin{equation*}
%            \lim_{n\to \infty} \max_{x \in A_n} q_n(x) = \max_{x \in A} q(x).
%         \end{equation*}
%         \item [(ii)] Let $h_n \in \arg\max_{x\in A_n}q_n(x)$, and suppose furthermore that $q$ admits unique maximizer $h^{\star}$ over $A$. Then, $h_n\to h^\star$. 
%     \end{enumerate} 
% \end{lemma}
% \begin{proof}
%     $(i)$. Using notation as in the statement of part $(ii)$, as $h_n \in A_n \subset K_m$ where $K_m$ is a compact set, we can extract a further subsequence $h_{n_j} \to \overline{h}$. As $A_n$ converge to $A$ in Hausdorff sense, $h \in A$. The approximate-maximizer inequality and uniform convergence imply
% \[
% q(\bar h)
% = \lim_{j}q(h_{n_j})=\lim_j q_{n_j}(h_{n_j})
% \ge \lim_j q_{n_j}(h^\star)
% =q(h^\star).
% \]
% This only happens when $\lim_{n\to \infty} \max_{x \in A_n} q_n(x) = \max_{x \in A} q(x)$, which implies part $(i)$ of Lemma \ref{lem:compact_convergence}.

% $(ii)$. In case of $h^\star$ is the unique maximizer of $q$ in $A$, we have $q(\overline h) = q(h^\star)$, which implies $\overline h = h^\star$. 
% \end{proof}

\begin{lemma}
\label{lem:dung_BL_distance_compare}
Let $X$ and $Y$ be real-valued random variables. Then, for any $R>0$, we have 
\begin{equation*}
    d_{\mathrm{BL}}(X^2,Y^2) \leq (2R+2)d_{\mathrm{BL}}(X,Y)+4\mathbb P(|X|>R-1).
\end{equation*}
\end{lemma}

\begin{proof}
Write $\Delta:=d_{\mathrm{BL}}(X,Y)$. 
For $R\ge1$, define
$g_R(x)=\min\{x^2,R^2\}$. Since $g_R$ is $2R$-Lipschitz,
for every $f\in\mathrm{BL}_1(\mathbb R)$ the function
$(f\circ g_R)/(2R)$ belongs to $\mathrm{BL}_1(\mathbb R)$.
Hence
\[
\left|
\mathbb E[f(g_R(X))]
-\mathbb E[f(g_R(Y))]
\right|
\le 2R\Delta.
\]
Moreover, $g_R(x)=x^2$ whenever $|x|\le R$, so
\[
\begin{aligned}
\left|
\mathbb E[f(X^2)]-\mathbb E[f(Y^2)]
\right|
&\le 2R\Delta
+2\mathbb P(|X|>R)
+2\mathbb P(|Y|>R).
\end{aligned}
\]

To bound the tail probability $\mathbb P(|Y|>R)$, define
$q_R(x)=\min\{1,(|x|-R+1)_+\}$.
Then $q_R\in\mathrm{BL}_1(\mathbb R)$ and
\[
\mathbf 1_{\{|x|>R\}}
\le q_R(x)
\le \mathbf 1_{\{|x|>R-1\}}.
\]
It follows that
\[
\mathbb P(|Y|>R)
\le \mathbb E[q_R(Y)]
\le \mathbb E[q_R(X)]+\Delta
\le \mathbb P(|X|>R-1)+\Delta.
\]
Taking the supremum over $f\in\mathrm{BL}_1(\mathbb R)$ gives
\[
d_{\mathrm{BL}}(X_n^2,X^2)
\le (2R+2)\Delta+4\mathbb P(|X|>R-1).
\]
\end{proof}

\begin{lemma}
\label{lem:compare_kolmogorov_BL}
Let $P,Q$ be probability measures on $\mathbb R$, with distribution
functions $F_P,F_Q$, then, for every $t\in\mathbb R$ and $\varepsilon>0$,
\[
|F_P(t)-F_Q(t)|
\le
\max\{1,\varepsilon^{-1}\}\,d_{\mathrm{BL}}(P,Q)
+
(F_Q(t+\varepsilon)-F_Q(t-\varepsilon))
.
\]
\end{lemma}

\begin{proof}
Fix $\varepsilon>0$ and define
\[
u(x)=
\begin{cases}
1, & x\le t,\\
1-(x-t)/\varepsilon, & t<x<t+\varepsilon,\\
0, & x\ge t+\varepsilon,
\end{cases}
\qquad
\ell(x)=u(x+\varepsilon).
\]
Both functions have sup norm at most $1$ and Lipschitz constant
at most $\varepsilon^{-1}$. Thus, with
$c_\varepsilon=\max\{1,\varepsilon^{-1}\}$,
\[
\left|\int u\,dP-\int u\,dQ\right|
\le c_\varepsilon d_{\mathrm{BL}}(P,Q),
\qquad
\left|\int \ell\,dP-\int \ell\,dQ\right|
\le c_\varepsilon d_{\mathrm{BL}}(P,Q).
\]
The pointwise inequalities $\mathbf 1_{(-\infty,t-\varepsilon]}
\le \ell
\le \mathbf 1_{(-\infty,t]}
\le u
\le \mathbf 1_{(-\infty,t+\varepsilon]}$
therefore imply
\[
F_Q(t-\varepsilon)-c_\varepsilon d_{\mathrm{BL}}(P,Q)
\le F_P(t)
\le F_Q(t+\varepsilon)+c_\varepsilon d_{\mathrm{BL}}(P,Q).
\]
Subtracting $F_Q(t)$ from both sides gives 
\begin{equation*}
    |F_P(t)-F_Q(t)|
\le
\max\{1,\varepsilon^{-1}\}\,d_{\mathrm{BL}}(P,Q)
+
\max\{F_Q(t+\varepsilon)-F_Q(t),\ F_Q(t)-F_Q(t-\varepsilon)\},
\end{equation*}
completing the proof. 
\end{proof}

\FloatBarrier
\bibliographystyle{plainnat}
\bibliography{references}

@article{Chen1995,
  author  = {Chen, Jiahua},
  title   = {Optimal rate of convergence for finite mixture models},
  journal = {Annals of Statistics},
  year    = {1995},
  volume  = {23},
  number  = {1},
  pages   = {221--233},
  doi     = {10.1214/aos/1176324464}
}

@article{DacunhaGassiat1999,
  author  = {Dacunha-Castelle, Didier and Gassiat, Elisabeth},
  title   = {Testing the order of a model using locally conic parametrization: Population mixtures and stationary {ARMA} processes},
  journal = {Annals of Statistics},
  year    = {1999},
  volume  = {27},
  number  = {4},
  pages   = {1178--1209},
  doi     = {10.1214/aos/1017938921}
}

@article{dempster1977maximum,
  author  = {Dempster, Arthur P. and Laird, Nan M. and Rubin, Donald B.},
  title   = {Maximum likelihood from incomplete data via the {EM} algorithm},
  journal = {Journal of the Royal Statistical Society, Series B},
  year    = {1977},
  volume  = {39},
  number  = {1},
  pages   = {1--38}
}

@article{Geyer1994,
  author  = {Geyer, Charles J.},
  title   = {On the asymptotics of constrained {M}-estimation},
  journal = {Annals of Statistics},
  year    = {1994},
  volume  = {22},
  number  = {4},
  pages   = {1993--2010},
  doi     = {10.1214/aos/1176325768}
}

@article{HeinrichKahn2018,
  author  = {Heinrich, Philippe and Kahn, Jonas},
  title   = {Strong identifiability and optimal minimax rates for finite mixture estimation},
  journal = {Annals of Statistics},
  year    = {2018},
  volume  = {46},
  number  = {6A},
  pages   = {2844--2870},
  doi     = {10.1214/17-AOS1641}
}

@article{HoNguyen2016,
  author  = {Ho, Nhat and Nguyen, XuanLong},
  title   = {On strong identifiability and convergence rates of parameter estimation in finite mixtures},
  journal = {Electronic Journal of Statistics},
  year    = {2016},
  volume  = {10},
  number  = {1},
  pages   = {271--307},
  doi     = {10.1214/16-EJS1105}
}

@article{HoNguyen2016Weak,
  author  = {Ho, Nhat and Nguyen, XuanLong},
  title   = {Convergence rates of parameter estimation for some weakly identifiable finite mixtures},
  journal = {Annals of Statistics},
  year    = {2016},
  volume  = {44},
  number  = {6},
  pages   = {2726--2755},
  doi     = {10.1214/16-AOS1444}
}

@article{Keribin2000,
  author  = {Keribin, Christine},
  title   = {Consistent estimation of the order of mixture models},
  journal = {Sankhya: The Indian Journal of Statistics, Series A},
  year    = {2000},
  volume  = {62},
  number  = {1},
  pages   = {49--66}
}

@article{KieferWolfowitz1956,
  author  = {Kiefer, Jack and Wolfowitz, Jacob},
  title   = {Consistency of the maximum likelihood estimator in the presence of infinitely many incidental parameters},
  journal = {Annals of Mathematical Statistics},
  year    = {1956},
  volume  = {27},
  number  = {4},
  pages   = {887--906}
}

@article{Leroux1992,
  author  = {Leroux, Brian G.},
  title   = {Consistent estimation of a mixing distribution},
  journal = {Annals of Statistics},
  year    = {1992},
  volume  = {20},
  number  = {3},
  pages   = {1350--1360}
}

@book{Lindsay1995,
  author    = {Lindsay, Bruce G.},
  title     = {Mixture Models: Theory, Geometry and Applications},
  series    = {NSF-CBMS Regional Conference Series in Probability and Statistics},
  volume    = {5},
  publisher = {Institute of Mathematical Statistics},
  year      = {1995}
}

@inproceedings{ManoleHo2022,
  author    = {Manole, Tudor and Ho, Nhat},
  title     = {Refined convergence rates for maximum likelihood estimation under finite mixture models},
  booktitle = {Proceedings of the 39th International Conference on Machine Learning},
  series    = {Proceedings of Machine Learning Research},
  volume    = {162},
  pages     = {14979--15006},
  year      = {2022}
}

@book{McLachlanPeel2000,
  author    = {McLachlan, Geoffrey J. and Peel, David},
  title     = {Finite Mixture Models},
  publisher = {Wiley},
  address   = {New York},
  year      = {2000}
}

@article{Nguyen2013,
  author  = {Nguyen, XuanLong},
  title   = {Convergence of latent mixing measures in finite and infinite mixture models},
  journal = {Annals of Statistics},
  year    = {2013},
  volume  = {41},
  number  = {1},
  pages   = {370--400},
  doi     = {10.1214/12-AOS1065}
}

@article{RednerWalker1984,
  author  = {Redner, Richard A. and Walker, Homer F.},
  title   = {Mixture densities, maximum likelihood and the {EM} algorithm},
  journal = {SIAM Review},
  year    = {1984},
  volume  = {26},
  number  = {2},
  pages   = {195--239},
  doi     = {10.1137/1026034}
}

@book{vdV1998,
  author    = {van der Vaart, Aad W.},
  title     = {Asymptotic Statistics},
  publisher = {Cambridge University Press},
  address   = {Cambridge},
  year      = {1998}
}

@book{vanDeGeer2000,
  author    = {van de Geer, Sara A.},
  title     = {Empirical Processes in {M}-Estimation},
  year      = {2000},
  publisher = {Cambridge University Press},
  address   = {Cambridge}
}

@book{vdVW1996,
  author    = {van der Vaart, Aad W. and Wellner, Jon A.},
  title     = {Weak Convergence and Empirical Processes: With Applications to Statistics},
  publisher = {Springer},
  address   = {New York},
  year      = {1996}
}

@book{Villani2009,
  author    = {Villani, C{\'e}dric},
  title     = {Optimal Transport: Old and New},
  publisher = {Springer},
  address   = {Berlin},
  year      = {2009}
}

@article{WuYang2020,
  author  = {Wu, Yihong and Yang, Pengkun},
  title   = {Optimal estimation of {Gaussian} mixtures via denoised method of moments},
  journal = {Annals of Statistics},
  year    = {2020},
  volume  = {48},
  number  = {4},
  pages   = {1981--2007},
  doi     = {10.1214/19-AOS1873}
}

@article{DvoretzkyKieferWolfowitz1956,
  author  = {Dvoretzky, Aryeh and Kiefer, Jack and Wolfowitz, Jacob},
  title   = {Asymptotic Minimax Character of the Sample Distribution Function and of the Classical Multinomial Estimator},
  journal = {The Annals of Mathematical Statistics},
  volume  = {27},
  number  = {3},
  pages   = {642--669},
  year    = {1956}
}

@article{Massart1990,
  author  = {Massart, Pascal},
  title   = {The Tight Constant in the {D}voretzky--{K}iefer--{W}olfowitz Inequality},
  journal = {The Annals of Probability},
  volume  = {18},
  number  = {3},
  pages   = {1269--1283},
  year    = {1990}
}

@article{Gasparetto2004,
  author = {Gasparetto, Maura and Gentry, Tracy and Sebti, Said and
            O'Bryan, Erica and Nimmanapalli, Ramadevi and Blaskovich, Michelle A. and
            Bhalla, Kapil and Rizzieri, David and Haaland, Perry and Dunne, Jack and Smith, Clay},
  title = {Identification of compounds that enhance the anti-lymphoma activity of rituximab using flow cytometric high-content screening},
  journal = {Journal of Immunological Methods},
  year = {2004}, volume = {292}, number = {1--2}, pages = {59--71},
  doi = {10.1016/j.jim.2004.06.003}
}

@article{LoHahne2009,
  author = {Lo, Kenneth and Hahne, Florian and Brinkman, Ryan R. and Gottardo, Raphael},
  title = {{flowClust}: a {Bioconductor} package for automated gating of flow cytometry data},
  journal = {BMC Bioinformatics},
  year = {2009}, volume = {10}, pages = {145},
  doi = {10.1186/1471-2105-10-145}
}

@article{NguyenLeRinaldoHo2026,
  author  = {Nguyen, Huy and Le, Dung and Rinaldo, Alessandro and Ho, Nhat},
  title   = {On the Geometry of Separation in Finite {Gaussian} Mixtures},
  journal = {arXiv preprint arXiv:2606.16179},
  year    = {2026},
  doi     = {10.48550/arXiv.2606.16179}
}

@article{DossEtAl2023,
  author  = {Doss, Natalie and Wu, Yihong and Yang, Pengkun and Zhou, Harrison H.},
  title   = {Optimal Estimation of High-Dimensional {Gaussian} Location Mixtures},
  journal = {The Annals of Statistics},
  year    = {2023},
  volume  = {51},
  number  = {1},
  pages   = {62--95},
  doi     = {10.1214/22-AOS2207}
}

@book{rockafellar,
  author    = {Rockafellar, R. Tyrrell and Wets, Roger J.-B.},
  title     = {Variational Analysis},
  year      = {1998},
  publisher = {Springer},
  address   = {Berlin, Heidelberg},
  series    = {Grundlehren der mathematischen Wissenschaften},
  volume    = {317},
  doi       = {10.1007/978-3-642-02431-3},
  isbn      = {978-3-540-62772-2}
}

@article{EmbrechtsHofert2013,
  author  = {Embrechts, Paul and Hofert, Marius},
  title   = {A note on generalized inverses},
  journal = {Mathematical Methods of Operations Research},
  year    = {2013},
  volume  = {77},
  number  = {3},
  pages   = {423--432},
  doi     = {10.1007/s00186-013-0436-7}
}

@article{HyndmanFan1996,
  author  = {Hyndman, Rob J. and Fan, Yanan},
  title   = {Sample Quantiles in Statistical Packages},
  journal = {The American Statistician},
  year    = {1996},
  volume  = {50},
  number  = {4},
  pages   = {361--365},
  doi     = {10.1080/00031305.1996.10473566}
}

@article{BarilettoLeRinaldoHo2026,
  author  = {Bariletto, Nicola and Le, Dung and
             Rinaldo, Alessandro and Ho, Nhat},
  title   = {Convergence rates for latent mixing measures in
             infinite homoscedastic location-scale mixture models},
  journal = {arXiv preprint arXiv:2605.06845},
  year    = {2026},
  doi     = {10.48550/arXiv.2605.06845}
}

@article{ManoleKhalili2021,
  author  = {Manole, Tudor and Khalili, Abbas},
  title   = {Estimating the Number of Components in Finite Mixture
             Models via the {Group-Sort-Fuse} Procedure},
  journal = {The Annals of Statistics},
  year    = {2021},
  volume  = {49},
  number  = {6},
  pages   = {3043--3069},
  doi     = {10.1214/21-AOS2072}
}

@article{DrtonPlummer2017,
  author  = {Drton, Mathias and Plummer, Martyn},
  title   = {A {Bayesian} Information Criterion for Singular Models},
  journal = {Journal of the Royal Statistical Society:
             Series B (Statistical Methodology)},
  year    = {2017},
  volume  = {79},
  number  = {2},
  pages   = {323--380},
  doi     = {10.1111/rssb.12187}
}

@article{NguyenEtAl2023Approximation,
  author  = {Nguyen, TrungTin and Chamroukhi, Faicel
             and Nguyen, Hien D. and McLachlan, Geoffrey J.},
  title   = {Approximation of probability density functions via
             location-scale finite mixtures in {Lebesgue} spaces},
  journal = {Communications in Statistics---Theory and Methods},
  year    = {2023},
  volume  = {52},
  number  = {14},
  pages   = {5048--5059},
  doi     = {10.1080/03610926.2021.2002360}
}

@article{cox2022generalized,
  author  = {Cox, Gregory},
  title   = {A generalized argmax theorem with applications},
  journal = {arXiv preprint arXiv:2209.08793},
  year    = {2022},
  doi     = {10.48550/arXiv.2209.08793}
}

\end{document}